\documentclass[final,3p,times,timeauthor]{elsarticle}

\usepackage{graphicx,amsmath,amsthm,amsfonts,latexsym,amssymb,hyperref,comment,subeqnarray,bm}
\usepackage{natbib}
\usepackage{hyperref,mathrsfs}
\usepackage[mathscr]{euscript}
\usepackage{colortbl}
\usepackage{todonotes}
\usepackage{amsmath}
\usepackage{caption}
\usepackage{subcaption}

\usepackage{algorithmic}
\let\chapter\section               
\usepackage{algorithm}  

\usepackage{lineno}

\newtheorem{remark}{Remark}[section]
\newtheorem{theorem}{Theorem}[section]

\newtheorem{lemma}{Lemma}[section]
\newtheorem{assumption}{Assumption}[section]
\numberwithin{equation}{section}

\newcommand{\vertiii}[1]{{\left\vert\kern-0.25ex\left\vert\kern-0.25ex\left\vert #1
    \right\vert\kern-0.25ex\right\vert\kern-0.25ex\right\vert}}
\newcommand{\verti}[1]{{\left\vert #1
    \right\vert}}

\newcommand{\aver}[1]{\left\{\!\!\left\{#1\right\}\!\!\right\}}

\newcommand{\jump}[1]{\left[\!\left[#1\right]\!\right]}

\newcommand{\dd}{\,{\rm d}}
\newcommand{\bfn}{{\bf n}}

\newcommand{\bfalpha}{{\bf \alpha}}
\newcommand{\bfx}{{\bf x}}
\newcommand{\bfy}{{\bf y}}
\newcommand{\bfA}{{\bf A}}
\newcommand{\bfu}{{\bf u}}
\newcommand{\bff}{{\bf f}}
\newcommand{\bfB}{{\bf B}}

\journal{Computer Methods In Applied Mechanics And Engineering}

\begin{document}

\begin{frontmatter}

\title{An Enriched Immersed Finite Element Method for 3D Interface Problems}

\tnotetext[mytitlenote]{Xu Zhang is partially supported by National Science Foundation DMS-2110833. 
Ruchi Guo was partially supported by National Science Foundation DMS-2309777.}

\author{Ruchi Guo}
\address{School of Mathematics, Sichuan University, Chengdu, China (ruchiguo@scu.edu.cn)}
\author{Xu Zhang}
\address{Department of Mathematics, Oklahoma State University, Stillwater OK 74078 (xzhang@okstate.edu)}

\begin{abstract}
We introduce an enriched immersed finite element method for addressing interface problems characterized by general non-homogeneous jump conditions. Unlike many existing unfitted mesh methods, our approach incorporates a homogenization concept. The IFE trial function set is composed of two components: the standard homogeneous IFE space and additional enrichment IFE functions. These enrichment functions are directly determined by the jump data, without adding extra degrees of freedom to the system. Meanwhile, the homogeneous IFE space is isomorphic to the standard finite element space on the same mesh. 
This isomorphism remains stable regardless of interface location relative to the mesh, 
ensuring optimal $\mathcal{O}(h^2)$ conditioning that is independent of the interface location 
and facilitates an immediate development of a multigrid fast solver;
namely the iteration numbers are independent of not only the mesh size but also the relative interface location.
Theoretical analysis and extensive numerical experiments are carried out in the efforts to demonstrate these features.
\end{abstract}


\end{frontmatter}


\section{Introduction}

\begin{figure}[htbp]
\begin{center}
\includegraphics[width=.3\linewidth]{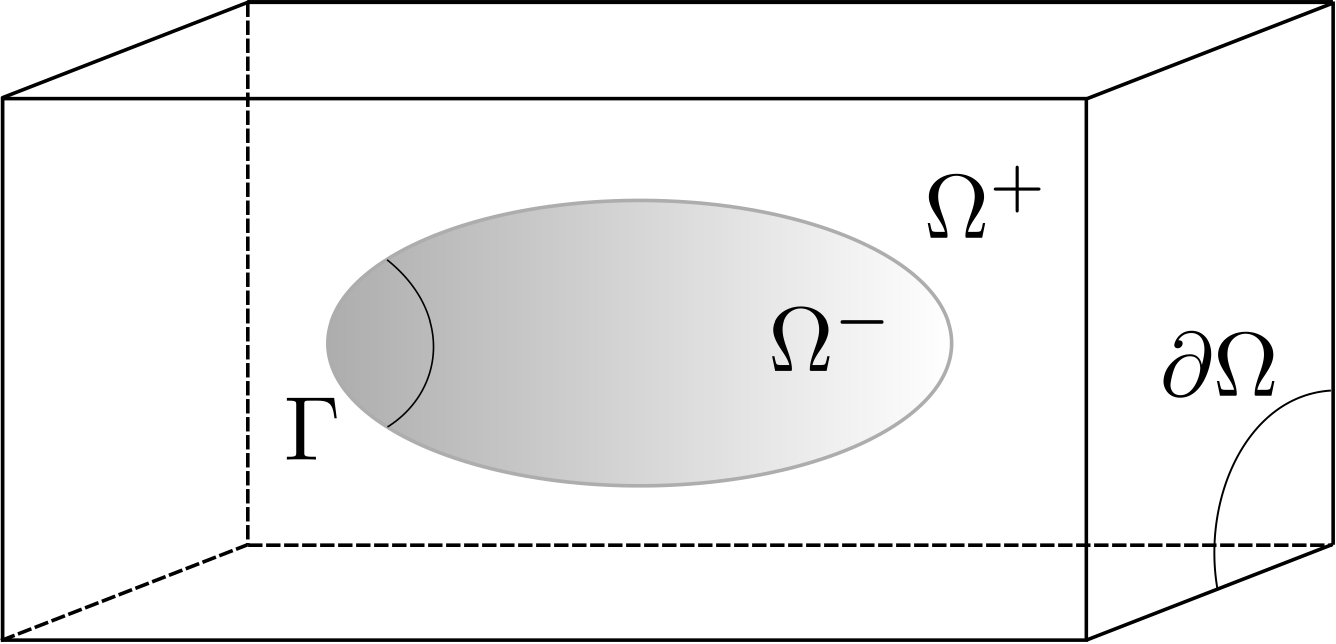}  
\caption{A three-dimensional domain with interface.}
\label{fig: domain}
\end{center}
\end{figure}

Let $\Omega\subseteq\mathbb{R}^3$ be an open bounded domain. 
We assume that $\Omega$ is separated into two subdomains $\Omega^-$ and $\Omega^+$ by a closed $C^2$ manifold $\Gamma\subseteq\Omega$ known as the interface.  
See Figure \ref{fig: domain} for an illustration. 
These subdomains contain different materials identified by a piecewise constant function 
$\beta(\mathbf{x})$ which is discontinuous across the interface $\Gamma$, i.e.,
\begin{equation*}
\beta(\mathbf{x})=
\left\{\begin{array}{cc}
\beta^- & \text{in} \; \Omega^- ,\\
\beta^+ & \text{in} \; \Omega^+,
\end{array}\right.
\end{equation*}
where $\beta^\pm>0$ and $\mathbf{x}= (x_1,x_2,x_3)$. 
We consider the following interface problem of the elliptic type on $\Omega$ 
coupled with non-homogeneous jump conditions:
\begin{subequations}
\label{model}
\begin{align}
\label{inter_PDE}
 -\nabla\cdot(\beta\nabla u)=f, &~~~~  \text{in} ~ \Omega^-  \cup \Omega^+, \\
 \jump{u}_{\Gamma} = q_1, &~~~~  \text{on} ~ \Gamma \label{jump_cond_1}, \\
 \jump{\beta \nabla u\cdot \mathbf{n}}_{\Gamma}  = q_2, & ~~~~\text{on}~\Gamma \label{jump_cond_2}, \\
  u=g, &~~~~\text{on} ~ \partial\Omega,
\end{align}
\end{subequations}
where $\mathbf{n}$ is the unit normal vector to $\Gamma$.
Here $\jump{v}_\gamma := (v|_{\Omega^+})|_\gamma-(v|_{\Omega^-})|_\gamma$ denotes the jump of $v$ on any manifold $\gamma\subseteq\Omega$, 
and if there is no danger of causing confusion, we shall drop $\gamma$ for simplicity. 
The functions $f$ and $g$ denote the source and boundary data, 
and $q_1$ and $q_2$ denote the solution and flux jump across the interface $\Gamma$. 
For simplicity, we denote $u^{\pm}=u|_{\Omega^{\pm}}$, in the rest of this article.


Interface problems as described in \eqref{model} are prevalent in physics and chemistry. 
For example, the non-zero flux jump $q_2 \neq 0$ occurs in models of electrical potentials across isotropic media with surface charge densities \cite{1975Cook}. 
Another instance is the Burton-Cabrera-Frank-type model for the epitaxial growth of thin films, discussed in \cite{2004BanschHauberLakkisLiVoigt, 2003CaflischLi}. 
Furthermore, the condition $q_1 \neq 0$ appears in the kinematic jump conditions of the Laplace-Young model for Hele-Shaw flows, as demonstrated in \cite{1997HouLiOsherZhao}. 
Additionally, cases where both $q_1 \neq 0$ and $q_2 \neq 0$ apply, such as in the Navier-Stokes equations, involve discontinuous pressures caused by surface tension and singular forces at interfaces, see \cite{2003LeeLeVeque,1997LevequeLi,2001LiLai} for example.

Many of the aforementioned problems are time-dependent for which \eqref{model} has to be solved at each time step with evolving interface. 
Classic numerical methods require fitted meshes for solving interface problems, 
such as finite element methods (FEM) \cite{1970Babuska,1998ChenZou}, discontinuous Galerkin (DG) methods \cite{2011CaiYeZhang} and so on;
otherwise their performance may be arbitrarily bad \cite{2000BabuskaOsborn}. 
It may make the simulation inefficient especially in the three-dimensional (3D) case
as generating high-quality meshes is usually very expensive.  
One remedy is to solve interface problems directly on unfitted meshes, 
which demands special modification to either the computation scheme or the approximation spaces. 
We refer readers to \cite{2015BurmanClausHansboLarsonMassing,2017HuangWuXiao,2016WangXiaoXu,2020LiuZhangZhangZheng} for penalty-type unfitted mesh FEMs
and \cite{1998Li,2003LiLinWu,2015LinLinZhang,2022JiWangChenLi} for immersed FEMs (IFEMs). 
In addition, a large class of finite-difference based methods have been also developed for structured Cartesian meshes 
such as immersed boundary methods \cite{2002Peskin}, immersed interface methods \cite{2006LiIto} and matched interface and boundary methods \cite{2006ZhouWei}.


Roughly speaking, in order to handle the non-homogeneous interface conditions on unfitted meshes, 
most of the aforementioned FEMs resort to some ``rich enough" approximation spaces on those elements intersecting with the interface. 
Let $T$ be a typical interface element, and take the linear case as an example.
The methods in \cite{2015ChenwuXiao,2017ChenWeiWen} generate a local triangulation of $T$ and employ standard FE space on this triangulation for approximation. 
Such a triangulation is also used by generalized FEM \cite{2012BabuskaBanerjee,2010FriesBelytschko,2017BabuskaBanerjeeKergrene} to construct enrichment functions which can be non-polynomials. 
Penalty-type methods \cite{2015BurmanClausHansboLarsonMassing,2017HuangWuXiao} employ two polynomial spaces $\mathcal{P}_1(T^-)$ and $\mathcal{P}_1(T^+)$
which are glued together through the Nitsche's penalty in the numerical schemes,
and the resulting space is $\mathcal{P}_1(T^-)\oplus\mathcal{P}_1(T^+)$.
While these spaces are richer than the basic $\mathcal{P}_1(T)$ and provide enhanced approximation capabilities, they also introduce additional degrees of freedom (DoFs), which in turn can deteriorate the conditioning of the resulting linear system. 
More precisely, this is due to that the condition number tends to be adversely affected by small cutting edges created by the interface. 
Some more complex penalties are needed to handle this issue \cite{2010Burman}. 

In this work, the proposed enriched IFEM enjoys a very different approach to handle the general non-homogeneous interface conditions as it does not introduce any additional unknown DoFs. 
On an interface element $T$, the approximation space is still $\mathcal{P}_1(T^-)\oplus\mathcal{P}_1(T^+)$ but decomposed into two subspaces:
\begin{equation}
\label{IFE_decomp}
\mathcal{P}_1(T^-)\oplus\mathcal{P}_1(T^+) = \mathcal{S}^0_h(T) \oplus \mathcal{S}^J_h(T),
\end{equation} 
where $\mathcal{S}^0_h(T)$ is a piecewise polynomial space determined by homogeneous jump conditions, 
i.e., the linear IFE space in the literature \cite{2005KafafyLinLinWang},
and $\mathcal{S}^J_h(T)$ is another piecewise polynomial space specifically constructed for the non-homogeneous jumps, called the enrichment IFE space.
See the precise definitions of these spaces in Section \ref{sec:IFEspaces}. Such a decomposition has the following two key features:
\begin{itemize}
\item In computation, functions in $\mathcal{S}^J_h(T)$ are determined immediately from the jump data $q_1$ and $q_2$. 
No extra unknowns are added to the global system. 

\item The local spaces $\mathcal{S}^0_h(T)$ can be put together through nodal continuity to form the global space $\mathcal{S}^0_h$ 
that is isomorphic to the standard linear FE space denoted by $\widetilde{\mathcal{S}}_h$, namely
\begin{equation}
\label{isomap0}
\mathbb{I}_h : \mathcal{S}^0_h \rightarrow \widetilde{\mathcal{S}}_h,
\end{equation}
see the precise definition in \eqref{isomap1}. 

\end{itemize}
The isomorphism in \eqref{isomap0} holds particular significance due to its ability to produce a highly structured linear system, which proves especially advantageous for addressing moving interface problems in both theoretical analysis and numerical computation \cite{2013HeLinLinZhang, 2020Guo, 2019AdjeridChaabaneLinYue}. 
Furthermore, \eqref{isomap0} is stable with respect to the interface location relative to the mesh, 
which is crucial for ensuring that the conditioning of the resulting matrix remains optimal, i.e., $\mathcal{O}(h^{-2})$, 
and also independent of the interface location.
The similar concept was also employed in \cite{2023Ji,2020AdjeridBabukaGuoLin} to show that the condition number is not affected by the interface location.
This property further enables the efficient application of multigrid methods (MG) for solving the linear system, thereby significantly enhancing computational efficiency.

The enrichment idea was originally introduced in \cite{2007GongLiLi,2010GongLi} for handling non-homogeneous interface conditions. 
In these works, the enrichment functions are typically constructed by level-set functions which exactly satisfy the jump conditions but may not be polynomials.
In \cite{2011HeLinLin}, the enrichment functions are constructed by piecewise polynomials weakly satisfying the interface conditions but restricted to the 2D case.
A similar idea can be found in \cite{2013HouSongWangZhao} that studies a Robin-type jump condition.
We also refer readers to \cite{2016HanWangHeLinWang} for the application to plasma simulation 
and \cite{2017WangHouShi} for the generalization to non-linear jump conditions. 
It is worthwhile to mention that the enrichment idea has been widely used in extended finite element methods (X-FEMs) \cite{2001DolbowMoesBelytschko,2000DolbowMoesBelytschko,2006VaughanSmithChopp}. 

However, the analysis for enriched IFE methods still remains largely unexplored, 
except for a 2D study in \cite{2020AdjeridBabukaGuoLin}. 
The primary objective of this work is, for the first time, to establish a comprehensive framework for 3D enriched IFE methods, encompassing the construction, convergence rates, and condition number. 
Remarkably, both the convergence rates and the condition number are shown to be optimal and independent of the interface location. 
It is worth noting that the theoretical techniques required to handle the enrichment functions are seldom discussed in the literature. 
One major challenge lies in developing an appropriate interpolation operator capable of incorporating the jump conditions.
Additionally, the regularity of the jump data further complicates the analysis, as those data are usually defined only on interface surfaces, i.e., the low-dimensional manifolds.

This article consists of 5 additional sections. 
In the next section, we present some preliminary results and assumptions. 
In Section \ref{sec:IFEspaces}, we construct IFE functions. 
In Section \ref{sec:analysis}, we analyze the interpolation errors. 
In Section \ref{sec:scheme}, we present the analysis of the IFE scheme including the solution errors and conditioning number.
In Section \ref{sec:num_sec}, we provide extensive numerical examples to demonstrate the effectiveness and efficiency of the proposed method.


\section{Preliminary}

\subsection{Mesh And Spaces}
Through this article, we assume that the unfitted mesh, denoted as $\mathcal{T}_h$, is generated by subdividing a background cubic mesh into regular tetrahedra, as shown in Figure \ref{fig:cubtet}. 
The collection of faces and nodes are denoted as $\mathcal{F}_h$ and $\mathcal{N}_h$, respectively. 
We also use $\mathcal{T}^i_h$ and $\mathcal{F}^i_h$ to denote the collection of elements and faces cut by the interface, respectively. 

Given a subregion $\omega\subseteq\Omega$, we let $H^s(\omega)$, $s>0$, be the standard Sobolev space on $\omega$.
If $\omega$ intersects with $\Gamma$ and is partitioned into two subregions $\omega^{\pm}$,
we introduce the piecewisely-defined space 
$$H^s(\omega^-\cup\omega^+) = \{ v \in L^2(\omega): v|_{\omega^{\pm}}\in H^s(\omega^{\pm}) \}.$$ 
In addition, we define $H^s_0(\omega)$ and $H^s_0(\omega^-\cup\omega^+)$ as their subspaces containing functions with zero traces on $\partial\omega$.
Furthermore, we let $\mathcal{P}_1(\omega)$ be the linear polynomial space defined on the region $\omega$. 
But, as polynomials are trivially well-defined everywhere, we will often use $\mathcal{P}_1$ by dropping $\omega$.

\subsection{Geometric Approximation}

\begin{figure}[h]
\centering
\begin{minipage}{.4\textwidth}
\begin{subfigure}{.4\textwidth}
    \includegraphics[width=1.2in]{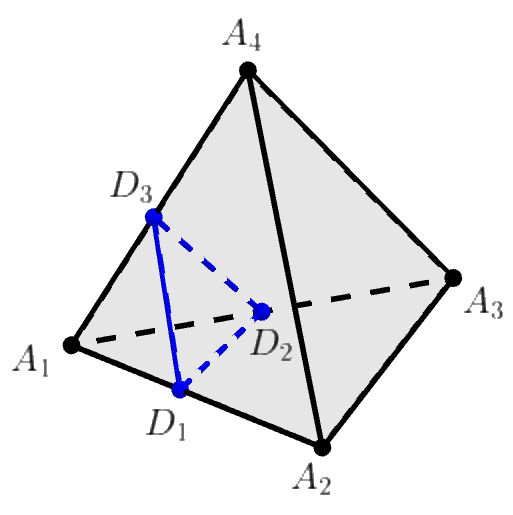}
\end{subfigure}
 \begin{subfigure}{.2\textwidth}
    \includegraphics[width=1.2in]{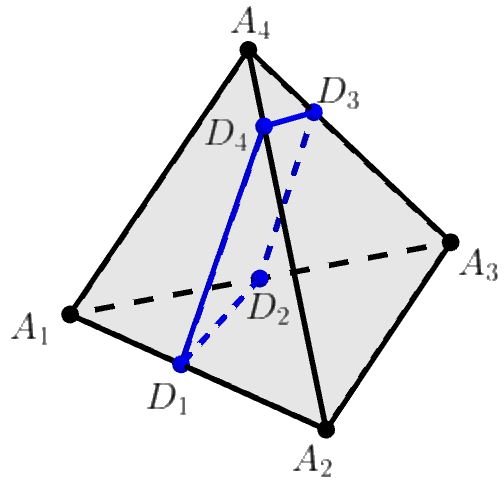}
\end{subfigure}
     \caption{Interface elements: 3 cutting points (left) and 4 cutting points (right).}
  \label{fig:tet} 
  \end{minipage}
  ~~~~~~~
  \begin{minipage}{.4\textwidth}
  \centering
  \includegraphics[width=1.1in]{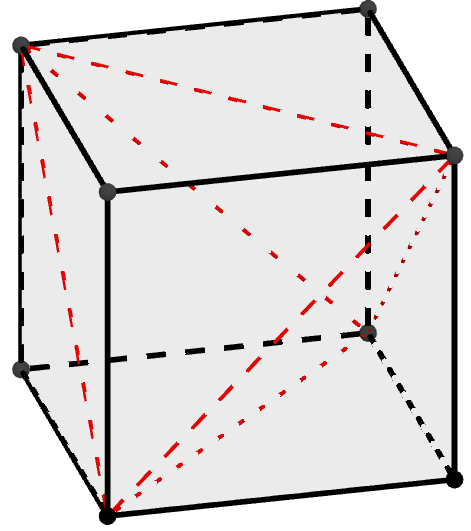}
    \caption{Partition a cuboid into 5 tetrahedra.}
  \label{fig:cubtet}
  \end{minipage}
\end{figure}

To begin with, we form a linear approximation to the interface that is then used in the construction of IFE functions.
For this purpose, we make the following assumption for the interface cutting meshes.
\begin{assumption}
\label{assump}
The interface can cut a tetrahedral element with only 3 or 4 points which locate at different edges.
\end{assumption}
By Assumption \ref{assump}, it is not hard to see that there exist only two possible interface element configuration shown in Figure \ref{fig:tet}. 
For the case of 3 cutting points, shown by the left plot of Figure \ref{fig:tet}, 
the plane determined by those points can be naturally used as the approximation to the actual surface. 
For the case of 4 cutting points, shown by the right plot of Figure \ref{fig:tet}, 
the points may not be coplanar, making the choice of an approximation plane more delicate.
Here, we follow \cite{2020GuoLin2} to employ the maximum angle condition \cite{1976BabuskaAziz},
which can theoretically guarantee the geometrical accuracy, see Lemma \ref{lem_geo_gamma} below.
Specifically, we select three points such that each interior angle is bounded above by $\theta < \pi$ and use these points to construct the approximation plane. 
The existence of such a selection is guaranteed by the following lemma.
We refer readers to \cite{2005KafafyLinLinWang,2005KafafyWangLin} for other geometrical approximation approaches. 
In particular, the level-set method, as employed in \cite{2022ChenGuoZou, 2023Ji}, ensures that all 3 or 4 intersection points lie on a common plane by construction.
The proposed analysis is also applicable to these methods.

\begin{lemma}
\label{max_angle}
Under Assumption \ref{assump}, there always exist 3 cutting points forming a triangle that satisfies the maximum angle condition.
\end{lemma}
\begin{proof}
The case of 3 cutting points is simple and immediately follows from the results in \cite{2017ChenWeiWen}. 
We focus on the case of 4 cutting points shown in the right plot of Figure \ref{fig:tet}. 
In fact, we only need to find one of the 4 angles $\theta_1:=\angle D_4D_1D_2$, $\theta_2:=\angle D_1D_2D_3$, $\theta_3:=\angle D_2D_3D_4$ and $\theta_4:=\angle D_3D_4D_1$ 
which is both bounded below from $0$ and above by $\pi$, i.e., such an angle is neither small nor large. 
Note that the two diagonals $D_1D_3$ and $D_2D_4$ both have the length of $\mathcal{O}(h_T)$, i.e., they never shrink, 
and every edge $D_1D_2$, $D_2D_3$, $D_3D_4$ and $D_4D_1$ is bounded above by $\mathcal{O}(h_T)$. 
Hence, the elementary geometry implies the existence of an angle $\underline{\theta}$ such that $\theta_1,\theta_2,\theta_3,\theta_4\ge \underline{\theta}$. 
Furthermore, as $D_1$, $D_2$, $D_3$ and $D_4$ may not be coplanar, we can conclude 
$$
\theta_1 + \theta_2 + \theta_3 + \theta_4 \le 2 \pi.
$$
Therefore, $\min\{\theta_1,\theta_2,\theta_3,\theta_4\}\in [\underline{\theta}, \pi/2]$ which finishes the proof.
\end{proof}
In the following discussion, without loss of generality, we assume that the 3 points given by Lemma \ref{max_angle} are always denoted as $D_1$, $D_2$ and $D_3$. 
We shall use the plane determined by these 3 points denoted as $\Gamma^T_h$ as the approximation to the surface $\Gamma^T:=\Gamma\cap T$. 

Furthermore, for each element $T$, we introduce its patch $\omega_T$ defined as
\begin{equation}
\label{patch}
\omega_T = \{ T'\in \mathcal{T}_h: \partial T' \cap \partial T \neq \emptyset \}.
\end{equation}
For analysis, we need the geometric estimate regarding between $\Gamma$ and $\Gamma^T_h$ within the patches of interface elements. 
Specifically, let $\Gamma^{\omega_T}=\Gamma\cap\omega_T$ and $\Gamma^{\omega_T}_h=\Gamma^T_h\cap\omega_T$. The geometric errors caused by $\Gamma^{\omega_T}$ and $\Gamma^{\omega_T}_h$ are described by the following lemma.
\begin{lemma}[Theorem 2.2,\cite{2020GuoLin2}]
\label{lem_geo_gamma}
Suppose $\Gamma$ is smooth enough and let the mesh size be sufficiently small. 
For every interface element $T$ and its patch $\omega_T$, denote $\bfn$ and $\bar{\bfn}$ as the normal vectors to $\Gamma^{\omega_T}$ and $\Gamma^{\omega_T}_h$ which have the same orientation. 
For each $\bfx\in\Gamma^{\omega_T}$, let $\bfx_{\bot}$ be the projection of $\bfx$ onto $\Gamma^{\omega_T}_h$.
Then, there hold that
\begin{subequations}
\label{lem_geo_gamma_eq0}
\begin{align}
 \| \bfx - \bfx_{\bot} \| \le c_{\Gamma,1} h^2_T,  \label{lem_geo_gamma_eq01}\\
 \| \bfn - \bar{\bfn} \| \le c_{\Gamma,2} h_T. \label{lem_geo_gamma_eq02}
\end{align}
\end{subequations}
\end{lemma}


\section{Then enriched IFE method}
\label{sec:IFEspaces}
In this section, we develop the enriched linear IFE spaces as well as the associated computation scheme.

\subsection{The enriched linear IFE space}
Let $T$ be an interface element with the vertices $A_1$, $A_2$, $A_3$, and $A_4$, see Figure \ref{fig:tet} for example. 
For any piecewise-defined function $\phi_T$ with $\phi_T^\pm = \phi_T|_{T^\pm}\in\mathcal{P}_1$,
we introduce the following linear functionals: 
\begin{subequations}
\begin{align}
& \mathcal{N}_i(\phi_T) = \phi_T(A_i), ~~~ i=1,2,3,4, \label{approx_jumps0} \\
& \mathcal{J}_i(\phi_T) = \jump{\phi_T}_{D_i},~~~i=1,2,3, ~~~~~~~\text{and}~~~~~   \mathcal{J}_4(\phi_T) = \jump{\beta\nabla \phi_T\cdot \bar{\mathbf{n}}}_{\Gamma^T_h}, \label{approx_jumps}
\end{align}
\end{subequations}
where $\jump{\cdot}_{D_i}$ is the jump evaluated at $D_i$. 
These linear functionals form the degrees of freedom (DoFs).
As $\nabla \phi_T$ is a constant vector, $\phi_T\cdot \bar{\mathbf{n}}$ becomes a constant.

\begin{figure}[ht]
  \centering
  \begin{subfigure}[b]{\textwidth}
    \centering
  \includegraphics[width=.2\linewidth]{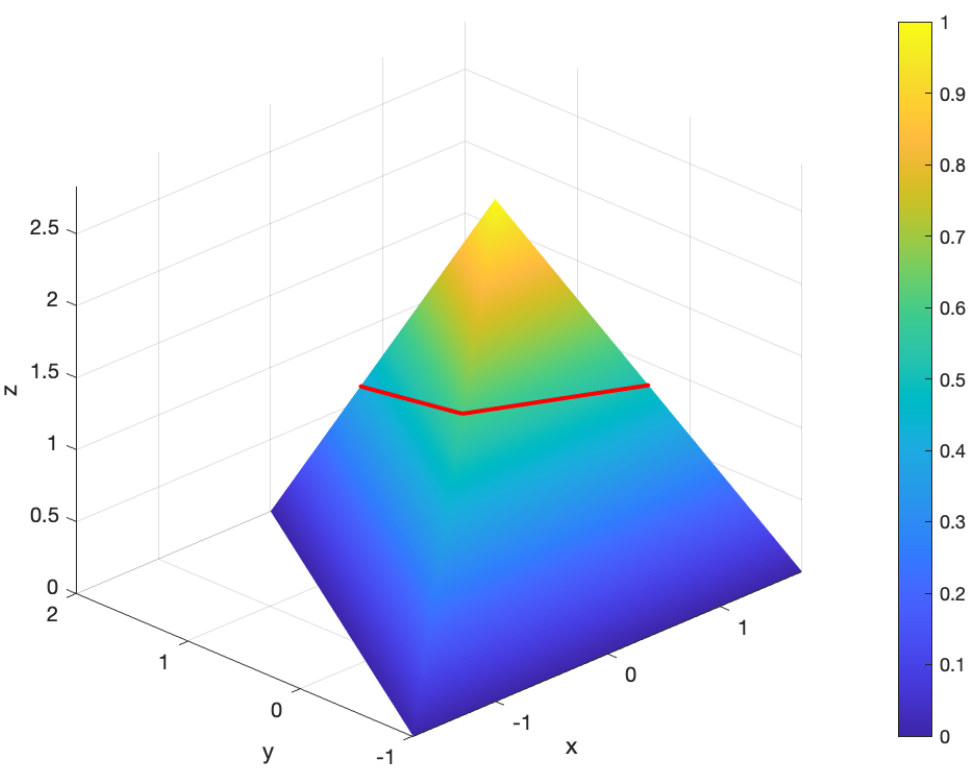}  
  \includegraphics[width=.2\linewidth]{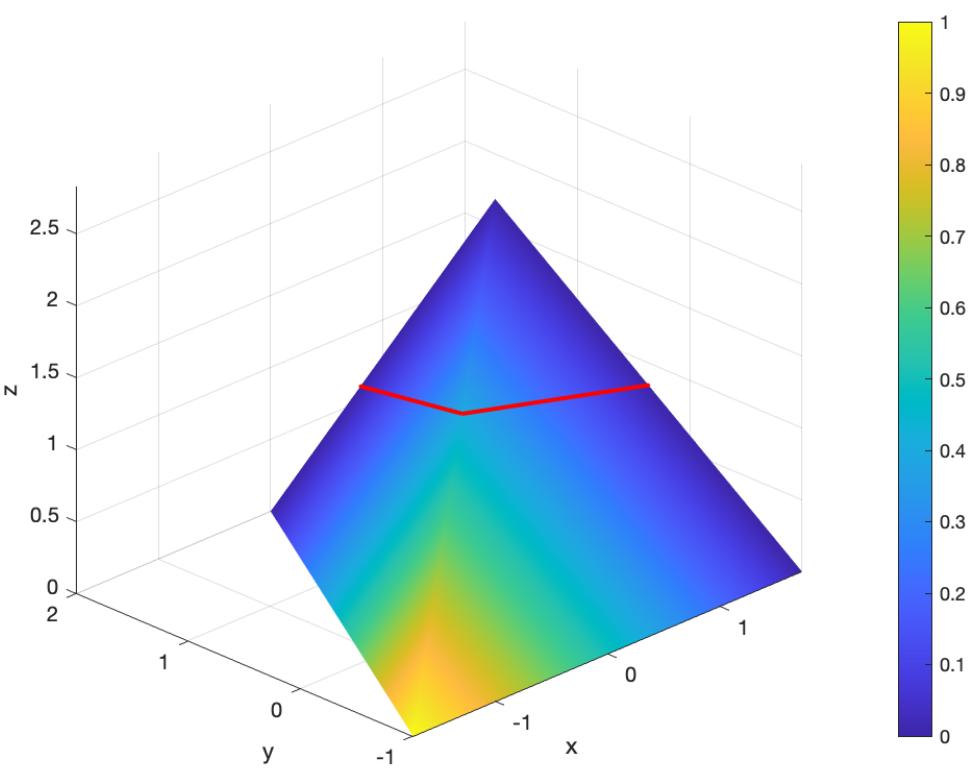}  
  \includegraphics[width=.2\linewidth]{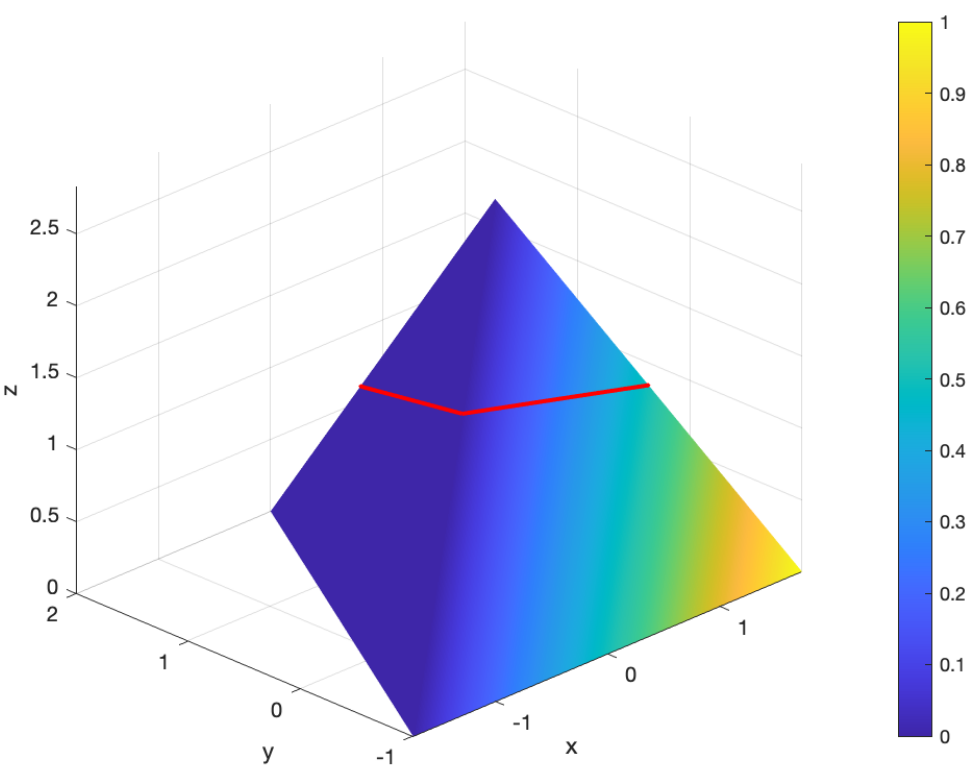}  
  \includegraphics[width=.2\linewidth]{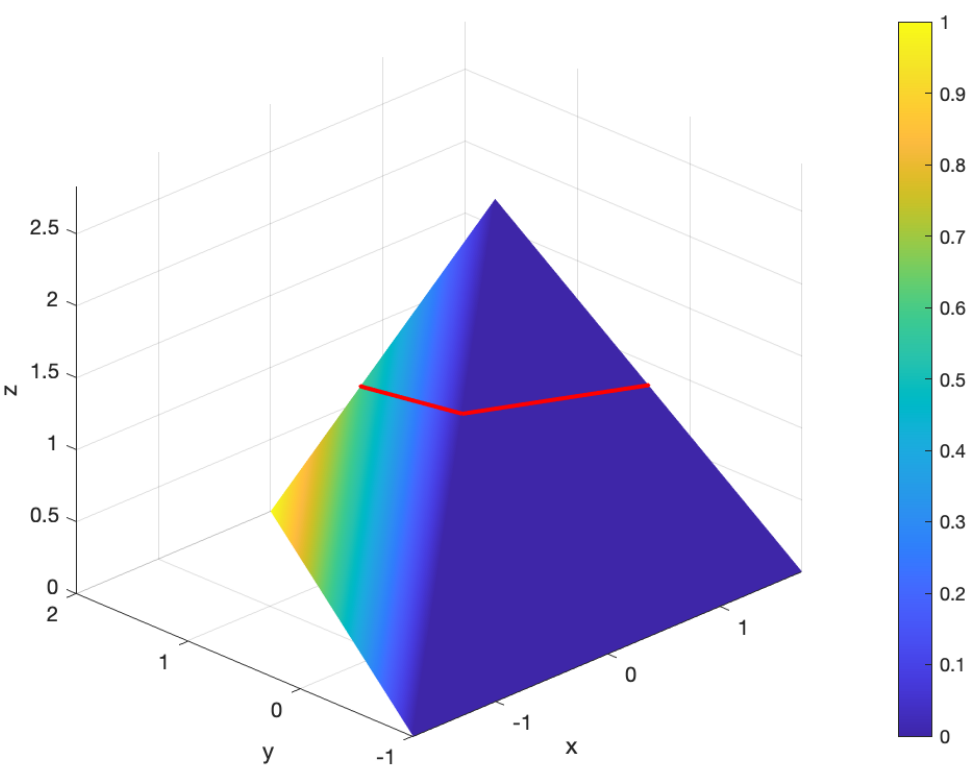}  
  \caption{Four $P_1$ FE shape functions $\psi_{j,T}$}
  \label{fig:sub-first}
\end{subfigure}
\centering
  \begin{subfigure}[b]{\textwidth}
    \centering
  \includegraphics[width=.2\linewidth]{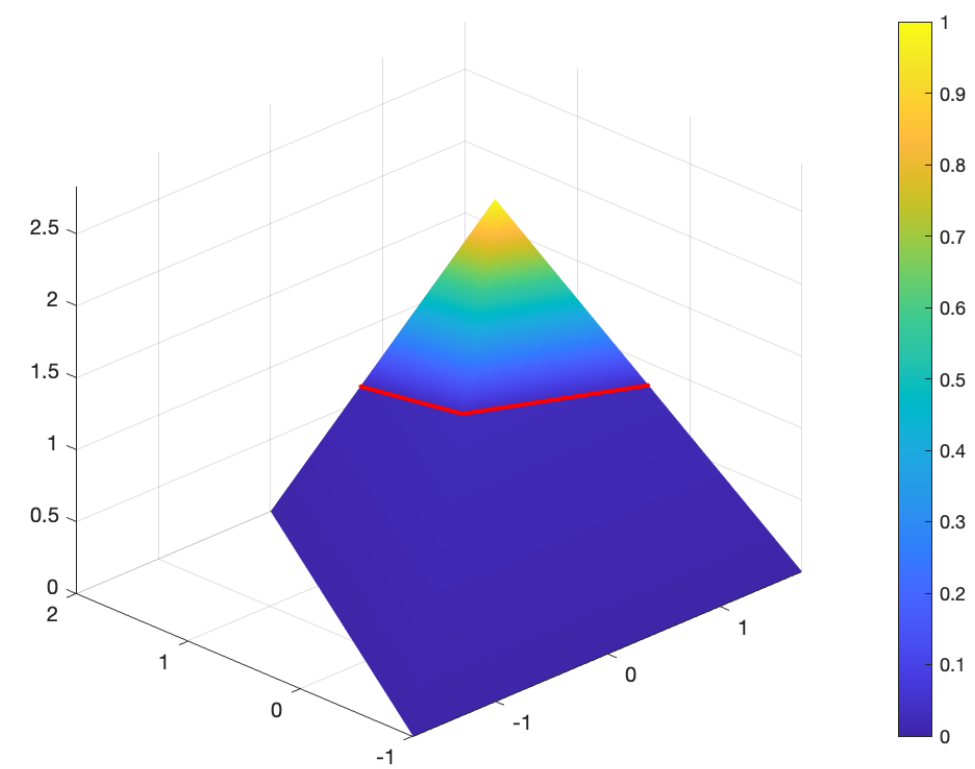}  
  \includegraphics[width=.2\linewidth]{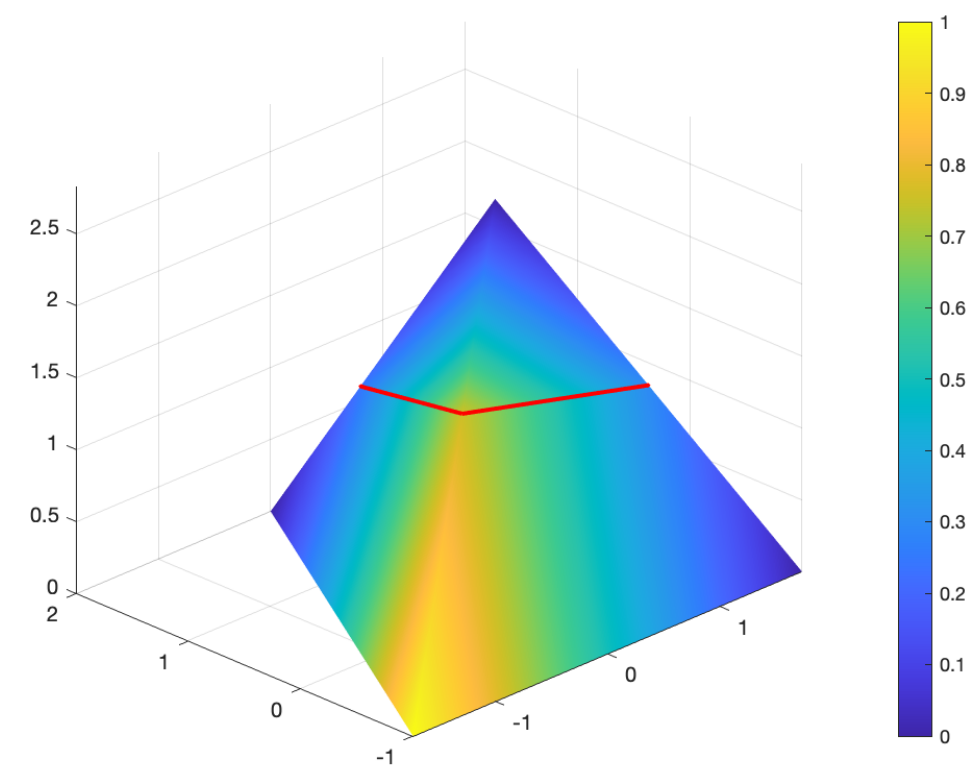}  
  \includegraphics[width=.2\linewidth]{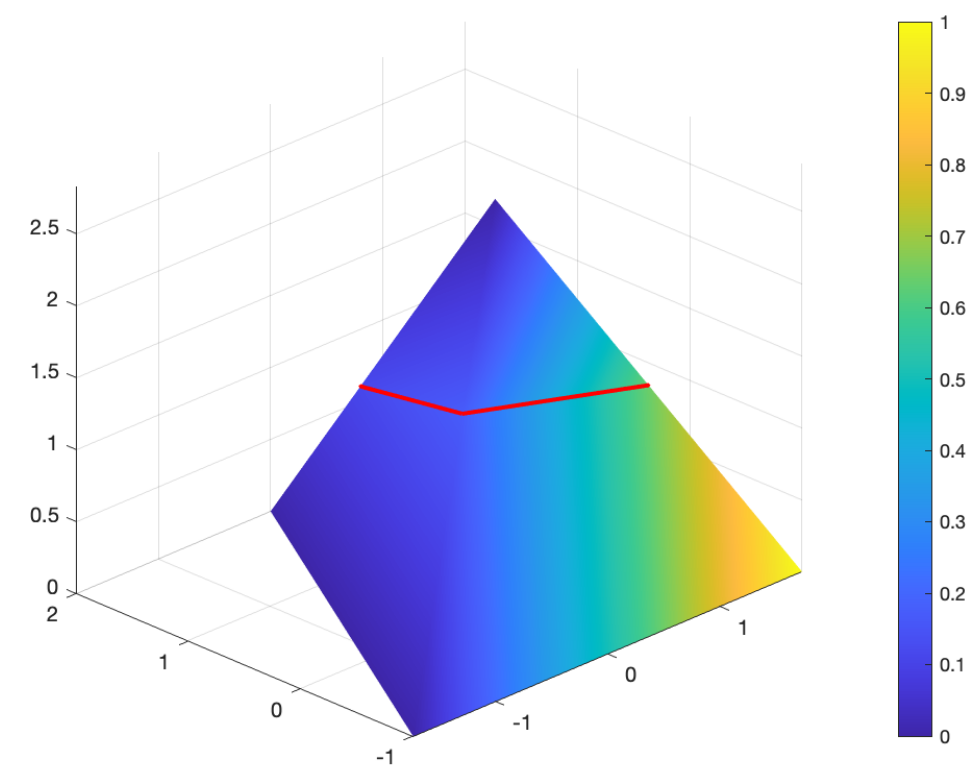}  
  \includegraphics[width=.2\linewidth]{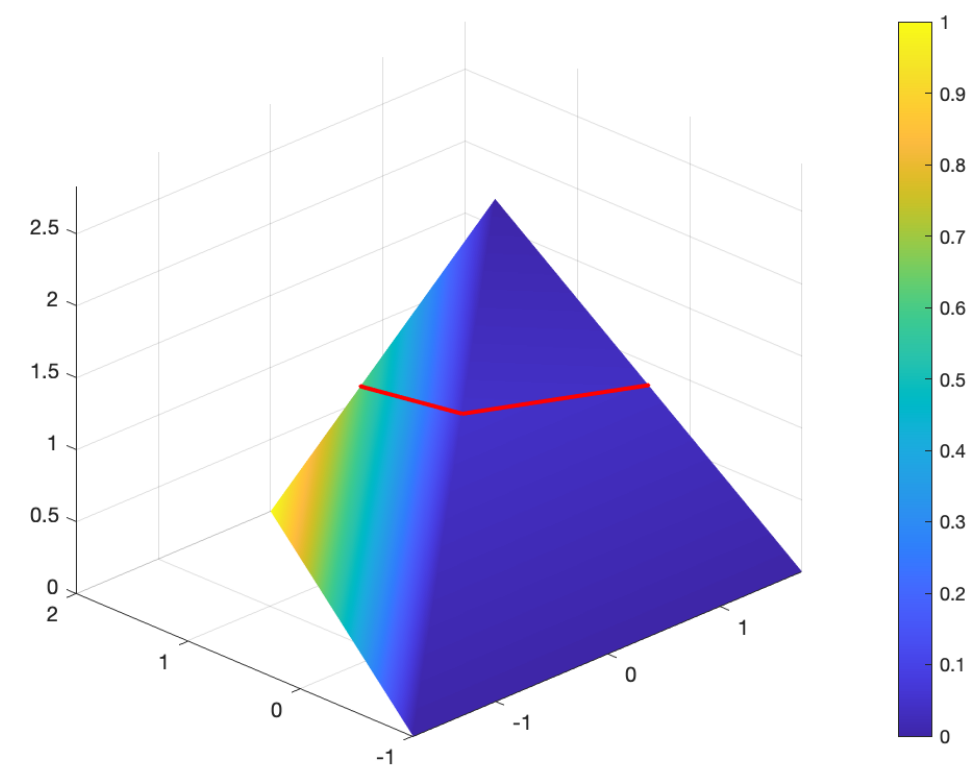}  
  \caption{Four $P_1$ IFE shape functions $\phi_{j,T}$ for homogeneous jump}
  \label{fig:sub-second}
\end{subfigure}
\centering
  \begin{subfigure}[b]{\textwidth}
    \centering
  \includegraphics[width=.2\linewidth]{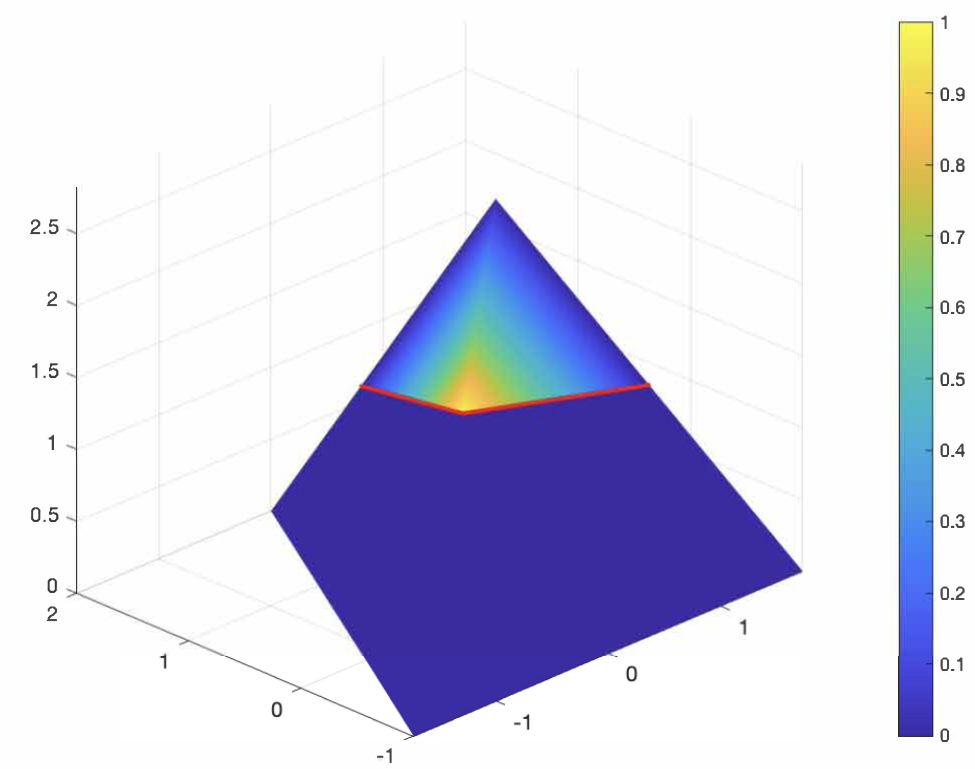}  
  \includegraphics[width=.2\linewidth]{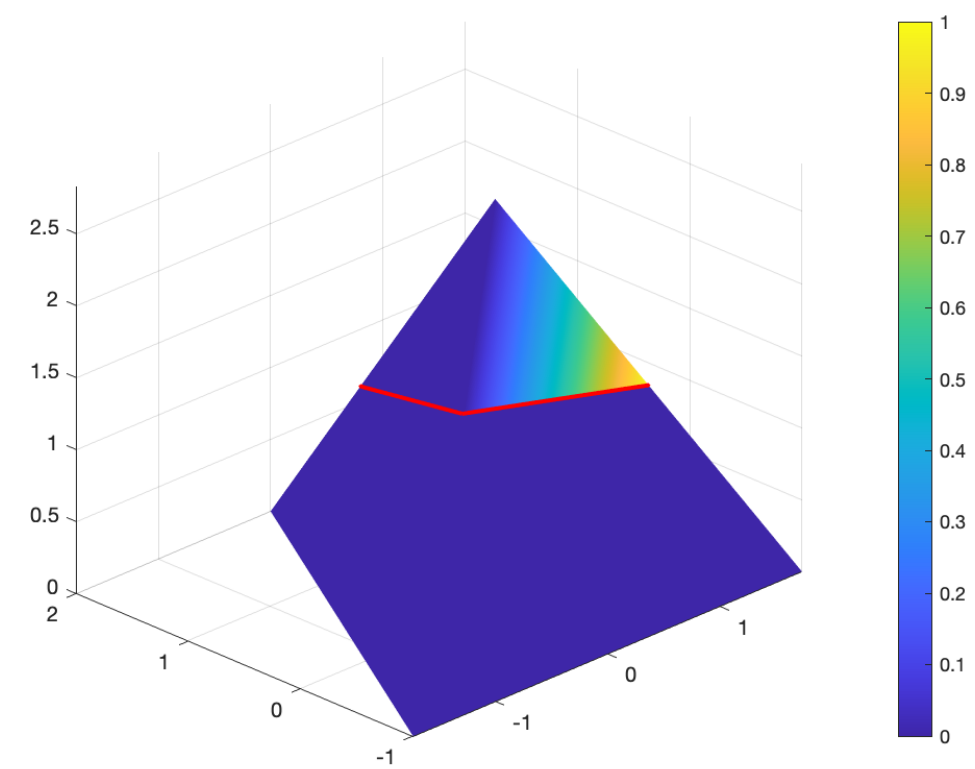}  
  \includegraphics[width=.2\linewidth]{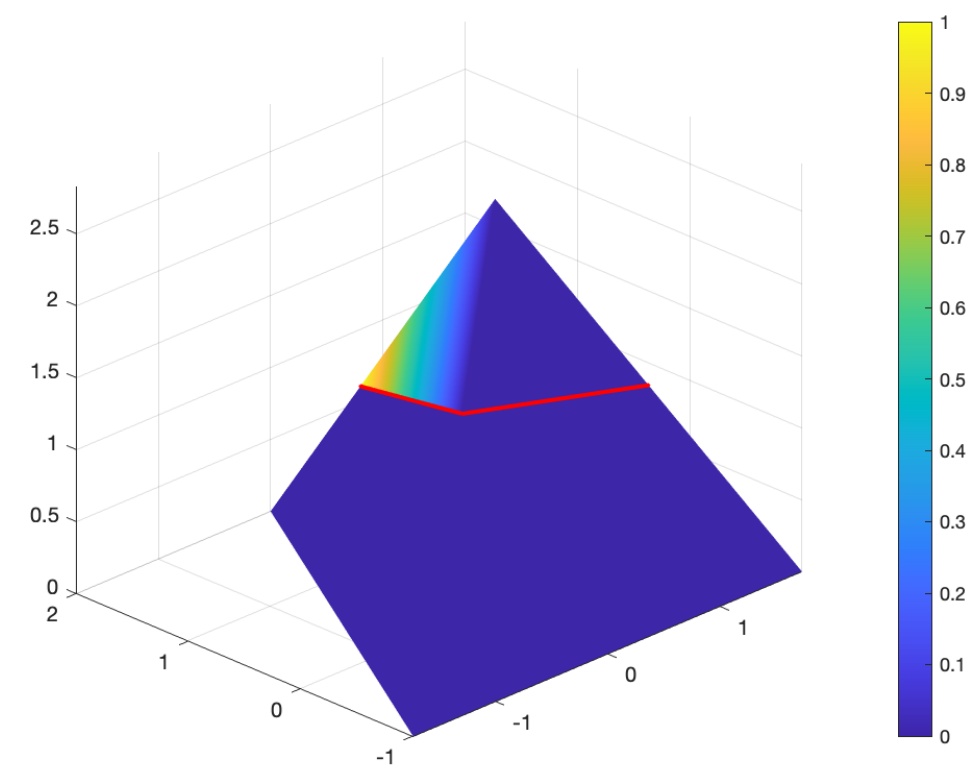}  
  \includegraphics[width=.2\linewidth]{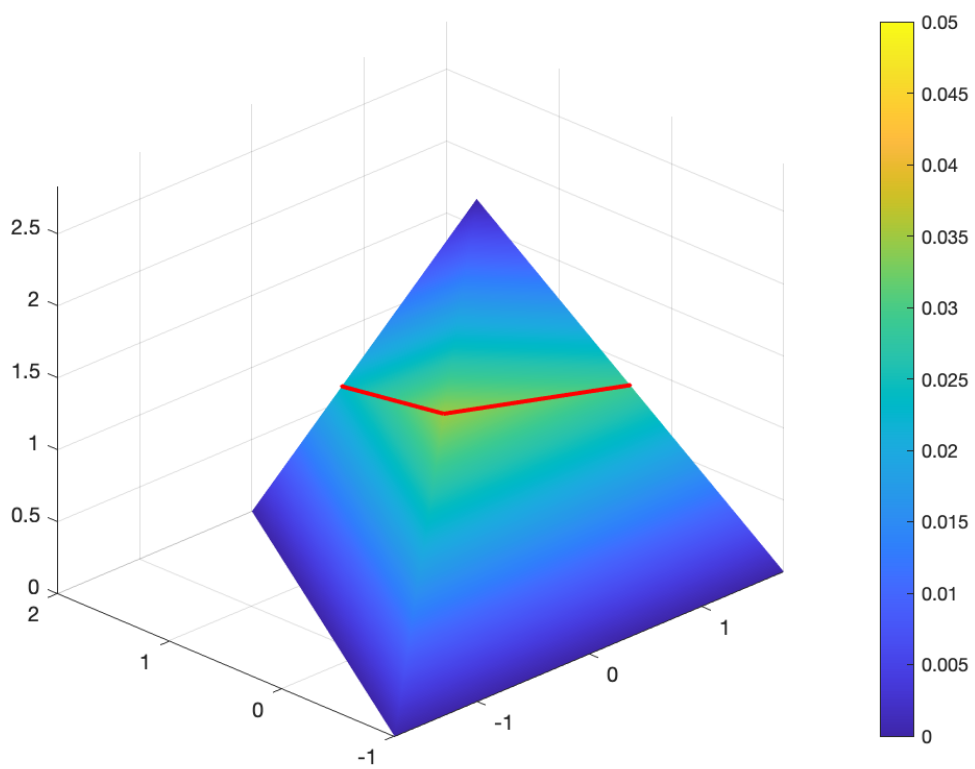}  
  \caption{Four $P_1$ IFE shape functions $\xi_{j,T}$ for non-homogeneous jump}
  \label{fig:sub-third}
\end{subfigure}  
\caption{Comparison of FE and IFE basis functions:
the homogeneous IFE shape functions $\phi_{j,T}$ are shown in second row, 
while non-homogeneous IFE shape functions $\xi_{j,T}$ are shown in the third row.}
\label{fig:basis}
\end{figure}

The homogeneous local linear IFE space $\mathcal{S}^0_h(T)$ contains the four shape functions $\phi_{j,T}$, $j=1,2,3,4$ such that 
\begin{equation}\label{eq: tmp1}
\mathcal{N}_i(\phi_{j,T}) = \delta_{ij}~~~~~\text{and}~~~~~~\mathcal{J}_i(\phi_{j,T}) = 0,~~~~\text{for}~~i,j=1,2,3,4.
\end{equation}
The non-homogeneous local linear IFE space $\mathcal{S}^J_h(T)$ contains another four shape functions $\xi_{j,T}$, $j=1,2,3,4$ such that 
\begin{equation}\label{eq: tmp2}
\mathcal{N}_i(\xi_{j,T}) = 0 ~~~~~\text{and}~~~~~~\mathcal{J}_i(\xi_{j,T}) = \delta_{ij},~~~~\text{for}~~i,j=1,2,3,4.
\end{equation}
Here, $\phi_{j,T}'s$ are basically the IFE version of the Lagrange-type nodal shape functions which have been well-studied in \cite{2005KafafyLinLinWang,2020GuoLin2,2022ChenGuoZou}. 
Let us briefly recall the construction procedure. We denote the coordinates of vertices $A_i = (x_i,y_i,z_i)$, $i=1,2,3,4$ and $D_j = (x_{dj},y_{dj},z_{dj})$, $j=1,2,3$. 
The equation \eqref{eq: tmp1} leads to a linear system
\begin{equation}\label{eq: system1}
\bfB\mathbf{c}_i =\mathbf{v}
\end{equation}
where 
\begin{equation}
\label{matA}
\bfB = 
\left(\begin{array}{cccccccc}
1 & x_{1} & y_{1} & z_{1} & 0 & 0 & 0 & 0 \\
1 & x_{2} & y_{2} & z_{2} & 0 & 0 & 0 & 0 \\
1 & x_{3} & y_{3} & z_{3}& 0 & 0 & 0 & 0 \\
0 & 0 & 0 & 0 & 1 & x_{4} & y_{4} & z_{4} \\
-1 & -x_{d1} & -y_{d1} & -z_{d1} & 1 & x_{d1} & y_{d1} & z_{d1} \\
-1 & -x_{d2} & -y_{d2} & -z_{d2} & 1 & x_{d2} & y_{d2} & z_{d2} \\
-1 & -x_{d3} & -y_{d3} & -z_{d3} & 1 & x_{d3} & y_{d3} & z_{d3} \\
0 & -\beta^-n_1& -\beta^-n_2 & -\beta^-n_3 & 0 & \beta^+n_1& \beta^+n_2 & \beta^+n_3 
\end{array}\right)
\end{equation}
and $\mathbf{v} = \mathbf{e}_i\in \mathbb{R}^8$, $i=1,2,3,4$ are the first four canonical unit vectors. 
Each $\mathbf{c}_i$ contains 8 unknowns, representing the coefficients of the terms $\{1, x, y, z\}$ on each side of an interface element.
The solvability of \eqref{eq: system1}, i.e., the non-singularity of \eqref{matA}, and the resulting unisolvence has been established in \cite{2005KafafyLinLinWang}.
Then the homogeneous local IFE space on $T$ can be written as 
\begin{equation}
\mathcal{S}^0_{h}(T) = \text{Span}\{\phi_{j,T}: j=1,2,3,4\}.
\end{equation}
To determine the non-homogeneous IFE shape functions, let us repeat \eqref{eq: system1} as
\begin{equation}\label{eq: system2}
\bfB{\mathbf{c}}_i =\mathbf{v},
\end{equation}
but modify the right-hand-side to be $\mathbf{v} =\mathbf{e}_i$, $i=5,6,7,8$ denotes the last four canonical vectors in  $\mathbb{R}^8$, 
and the unknowns ${\mathbf{c}}_i$ have the same meaning. 
Note that \eqref{eq: system1} and \eqref{eq: system2} share the same the coefficient matrix $\bfB$. 
Thus, the system \eqref{eq: system2}, with only a different right-hand side, is also uniquely solvable. 
Now, we define the non-homogeneous IFE spaces as
\begin{equation}
\label{IFE_nonhomo}
 \mathcal{S}_{h}^J(T) = \text{Span}\{\xi_{j,T}: j=1,2,3,4\}.
\end{equation}
All these shape functions are illustrated in Figure \ref{fig:basis}.
Therefore, we write down
\begin{equation}
\label{directSum}
\mathcal{S}_h(T) := \mathcal{S}^0_{h}(T) + \mathcal{S}_{h}^J(T) = \mathcal{P}_1(T^+)\oplus \mathcal{P}_1(T^-) := \{ v_h : v_h|_{T^{\pm}} \in \mathcal{P}_1(T^{\pm}) \},
\end{equation}
which will be shown to have the desired properties discussed in the introduction below \eqref{IFE_decomp}. 
It will be demonstrated that the functions in $\mathcal{S}_{h}^J(T)$ can be determined a priori directly from the data, rather than being solved from the large global linear system. 
Then, These ``enrichment terms" can then be straightforwardly transferred to the right-hand side of the computational scheme, 
resulting in a linear system that differs from the homogeneous case solely in its modified right-hand side, see the discussion in the next subsection.


\begin{remark}
\label{rem_approxi}
While the approximation plane $\Gamma^{T}_h$ is used in the construction procedure \eqref{approx_jumps}, 
the resulting IFE functions are still piecewisely defined with the actual interface. 
As the linear case is considered, the minor difference between $\Gamma^T_h$ and $\Gamma^T$ is not critical for computation. 
In fact, for each homogeneous IFE function $v_h\in \mathcal{S}^0_h(T)$, 
we can show the following estimate regarding their high-order difference on $\Gamma$ \cite[Theorem 4.8]{2020GuoZhang}:
\begin{equation}
\label{thm_trace_interface_eq0}
\| \jump{v_h} \|_{L^2(\Gamma^T)} \lesssim h^{3/2}_T \| \nabla v_h \|_{L^2(T)}.
\end{equation}
\end{remark}

\begin{remark}
\label{rem_unisol}
Usually, a finite element function can be uniquely determined by its DoFs. 
For IFE functions, given every $v_h\in \mathcal{S}_h(T)$, the unisolvence can yield the following identity
\begin{equation}
\label{interp_1}
v_h = \sum_{j=1}^4 v_h(A_j) \phi_{j,T} + \sum_{j=1,2,3} \jump{v_h}_{D_j} \xi_{j,T}+  \jump{\beta \nabla v_h \cdot\bar{\bfn}} ~ \xi_{4,T}.
\end{equation}
\end{remark}

Now, the global IFE spaces are defined as
\begin{subequations}
\begin{align}
    &   \mathcal{S}_{h} = \{ v_h\in L^2(\Omega)~:~ v_h|_T \in \mathcal{P}_1(T), ~ \text{if} ~T\in\mathcal{T}^n_h, ~~ \text{and} ~~  v_h|_T \in \mathcal{S}_{h}(T), ~ \text{if} ~T\in\mathcal{T}^i_h, ~ v_h|_{\partial\Omega} = 0 \},  \label{global_IFE_space_eq1} \\
     &   \mathcal{S}^0_{h} = \{ v_h\in L^2(\Omega)~:~ v_h|_T \in \mathcal{P}_1(T), ~ \text{if} ~T\in\mathcal{T}^n_h, ~~ \text{and} ~~  v_h|_T \in \mathcal{S}^0_{h}(T), ~ \text{if} ~T\in\mathcal{T}^i_h, ~ v_h|_{\partial\Omega} = 0 \},  \label{global_IFE_space_eq2} 
\end{align}
\end{subequations}
where the super script ``$0$" also emphasizes that this space corresponds to the homogeneous jumps.

Next, we present the following estimates for the shape functions.
\begin{theorem}[Bounds of IFE shape functions]
\label{bounds_IFEshapeFun}
For $k=0,1$, there hold
\begin{subequations}
\label{bounds_IFEshapeFun_eq0}
\begin{align}
   &   | \phi_{j,T} |_{W^{k,\infty}(\omega_T)} \lesssim h^{-k}_T,  ~~~~j=1,2,3,4,\\
   &   | \xi_{j,T} |_{W^{k,\infty}(\omega_T)} \lesssim h^{-k}_T, ~~~~~j=1,2,3, ~~~~ | \xi_{4,T} |_{W^{k,\infty}(\omega_T)} \lesssim h^{-k+1}_T. 
\end{align}
\end{subequations}
\end{theorem}
\begin{proof}
We consider a reference cuboid that has the sides of length $1$ from which the tetrahedral elements are cut. 
Let $\hat{\phi}_{j,T}$ and $\hat{\xi}_{j,T}$, $j=1,...,4$ be the IFE functions on the reference elements obtained from the affine mapping. 
Note that the nodal values of $\hat{\phi}_{j,T}$, $j=1,...,4$, and $\hat{\xi}_{j,T}$, $j=1,...,3$, stay unchanged.  
As for $\hat{\xi}_{4,T}$, it is not hard to see that $\jump{\beta\nabla\hat{\xi}_{4,T}\cdot\bar{\bfn}} = h_e$ with $h_e$ being the length of edges of the original cuboid. 
By the analysis in \cite{2005KafafyLinLinWang}, we know that $\text{det}(\hat{\bfB})\gtrsim 1$. 
Using the affine mapping again, we have the desired estimates.
\end{proof}


\subsection{The enriched IFE scheme}
Now, we present the IFE scheme that employs the enrichment functions to handle the non-homogeneous jump conditions. 
Noticing that $\jump{u}_{\Gamma}=q_1$ and $\jump{\beta \nabla u\cdot\bfn}_{\Gamma}=q_2$, we construct 
\begin{equation}
\label{qh1}
q_{T,1}(u) = \sum_{j=1,2,3} q_1(D_j) \xi_{j,T} = \sum_{j=1,2,3} \jump{u}_{D_j} \xi_{j,T} 
\end{equation}
which is well-defined since $q_1\in H^{3/2}(\Gamma)\hookrightarrow C^0(\Gamma)$ with $\Gamma$ being a 2D manifold. As for the flux jump, we construct
\begin{equation}
\label{qh2}
q_{T,2}(u) = \frac{1}{|\Gamma^{\omega_T}|}\int_{\Gamma^{\omega_T}} q_2 \dd s ~ \xi_{4,T} = \frac{1}{|\Gamma^{\omega_T}|}\int_{\Gamma^{\omega_T}}  \jump{\beta \nabla u\cdot\bfn }  \dd s ~ \xi_{4,T},
\end{equation}
where we use the integral here because $q_2$ does not possess sufficient regularity for pointwise evaluation.
Then, the global enrichment function is defined as
\begin{equation}
\label{qh}
q_{h,1} = 
\begin{cases}
      & q_{T,1}(u) ~~~~ \text{if} ~ T\in \mathcal{T}^i_h, \\
      & 0,~~~~~~~~~~~  \text{otherwise},
\end{cases}
~~~~
\text{and}
~~~~
q_{h,2} = 
\begin{cases}
      & q_{T,2}(u) ~~~~ \text{if} ~ T\in \mathcal{T}^i_h, \\
      & 0,~~~~~~~~~~~  \text{otherwise}.
\end{cases}
\end{equation}

\begin{remark}
\label{rem_q2}
If $q_2\in H^s(\Gamma)$, $s>1$, one can then simply take $q_{T,2}(u) = q_2(X)$ for any $X\in \Gamma^T$ in computation. 
Our analysis is also applicable to this case of smooth data. 
In addition, the definition of \eqref{qh2} involves $\Gamma^{\omega_T}$, which is also the key to treat the low regularity of $q_2$. 
\end{remark}

Now, we introduce a mesh-dependent broken space:
\begin{equation}
\label{broken_space}
\mathcal{V}_h = \{ v_h\in L^2(\Omega): v_h|_T \in H^1(T) ~ \text{if} ~ T\in \mathcal{T}^n_h, ~~ v_h|_{T^{\pm}} \in H^1(T^{\pm}) ~ \text{if} ~ T\in \mathcal{T}^i_h, ~ v_h |_{\partial\Omega} = 0 \}.
\end{equation}
Define the bilinear and linear forms:
\begin{subequations}
 \begin{align}
 a_h(u_h,v_h) & = \sum_{T\in\mathcal{T}_h} \int_T \beta \nabla u_h\cdot \nabla v_h \dd X \nonumber   - \sum_{F\in\mathcal{F}^i_h} \int_F \{ \beta \nabla u_h\cdot \mathbf{ n} \} \jump{v_h} \dd s \\
  & - \sum_{F\in\mathcal{F}^i_h} \int_F \{ \beta \nabla v_h \cdot \mathbf{ n} \} \jump{u_h} \dd s + \sum_{F\in\mathcal{F}^i_h} \frac{\sigma }{|e|} \int_F  \jump{u_h} \, \jump{v_h} \dd s, ~~~~~~ \forall u_h,v_h \in \mathcal{V}_h  \label{weak_form_2}
\end{align}
with a constant $\sigma>0$ large enough, and
\begin{align}
   L_{f,q}(v_h)&=\int_{\Omega}fv_h \dd X + \int_{\Gamma}q_2\{v_h\}_{\Gamma} \dd s , ~~~ \forall v_h \in \mathcal{V}_h.
\end{align}
\end{subequations}
Then, the proposed enriched IFE scheme is find $u_h = u^0_h + q_{h,1} + q_{h,2}\in \mathcal{S}_h$ with $u^0_h\in  \mathcal{S}^0_h$ such that
\begin{equation}
\label{IFE_scheme}
a_h(u^0_h,v_h) = L_{f,q}(v_h) - a_h(q_{h,1}+q_{h,2},v_h), ~~ \forall v_h \in \mathcal{S}^0_{h}.
\end{equation}

We end this section with introducing a norm based on the bilinear form:
\begin{equation}
\label{energy_norm}
\vertiii{v_h}^2 := \sum_{T\in\mathcal{T}_h} \| \sqrt{\beta} \nabla v_h\|^2_{L^2(T)} + \sum_{F\in\mathcal{F}^i_h } \sigma \| h^{-1/2} \jump{v_h} \|^2_{L^2(F)} + \sum_{F\in \mathcal{F}^i_h } \frac{1}{\sigma} \| h^{1/2} \aver{\beta\nabla v_h\cdot\mathbf{n}} \|^2_{L^2(F)}. 
\end{equation}
\begin{lemma}
\label{lem_full_norm}
$\vertiii{\cdot}$ is a norm on $\mathcal{V}_h$.
\end{lemma}
\begin{proof}
Assuming $\vertiii{v}=0$, we have $\|\nabla v\|=0$, and thus $v$ is a constant on each element. 
Then, the continuity at mesh nodes implies that $v_h$ must vanish on the whole domain.
\end{proof}
Note that the scheme \eqref{IFE_scheme} essentially follows from a homogenization idea used in PDE analysis. 
However, the error estimation for such a discretization scheme rarely appears in the literature and thus brings us new challenges. 


\section{Interpolation with non-homogeneous data}
\label{sec:analysis}
In this section, we show the approximation capability of $\mathcal{S}_h$ to $H^2(\Omega^-\cup\Omega^+)$. 
For each $u\in H^2(\Omega^-\cup\Omega^+)$, let $u^+_E\in H^2_0(\Omega)$ and $u^-_E\in H^2_0(\Omega)$ 
be the Sobolev extensions of the components $u^{\pm}\in H^2_0(\Omega^{\pm})$ such that
\begin{equation}
\label{sobolev_ext}
\| u^s_E \|_{H^2(\Omega)} \leqslant C_E \| u^s \|_{H^2(\Omega^s)}, ~~~ s=\pm,
\end{equation}
for some constant $C_E$ only depending on the geometry of $\Omega^{\pm}$ \cite{2001GilbargTrudinger}. 
Given such extensions, a straightforward interpolation operator is $\tilde{I}_Tu: = (I_Tu^+_E|_{T^+},I_Tu^-_E|_{T^-})$ with $I_T$ being the standard Lagrange interpolation on an element $T$. 
However, this operator is not suitable for estimation, as the interpolant in $\mathcal{S}^J_h(T)$ must be computable using jump data and $u^{\pm}$. Therefore, a more delicate interpolation operator is required for the analysis.

\subsection{Some fundamental estimates}

We begin with recalling the strip argument from \cite{2010LiMelenkWohlmuthZou} 
that is used to handle the mismatching portion around the interface denoted by $\omega^{\text{int}}_T$, 
i.e., the region sandwiched by the $\Gamma^{\omega_T}$ and $\Gamma^{\omega_T}_h$. 
First define the $\delta$-strip $S_{\delta} = \{ \bfx: \text{dist}(\bfx,\Gamma) \le \delta \}$.
Note that Lemma~\ref{lem_geo_gamma} immediately implies that
\begin{equation}
\label{delta_strip_1}
\cup_{T\in\mathcal{T}^i_h} \omega^{\text{int}}_T \subset S_{\delta}, ~~~~ \text{with} ~~ 0<\delta \le c_0 h^2.
\end{equation}
We can control the $L^2$-norm on the $\delta$-strip by its width through the following lemma.
\begin{lemma}[Lemma 3.1, \cite{2012HiptmairLiZou}]
\label{lem_delta}
It holds for any $z\in H^1_0(\Omega^{\pm})$ that
\begin{equation}
\label{lem_delta_eq0}
\| z \|_{L^2(S_{\delta})} \le C_{\Omega} \sqrt{\delta} \| z \|_{H^1(\Omega)}.
\end{equation}
\end{lemma}

In the following discussion, we let $\Pi^k_{\omega}$ be the $L^2$ projection onto the polynomial space $\mathcal{P}_k(\omega)$ on a region $\omega$.
Note that $q_j$, $j=1,2$, are only defined on $\Gamma^{\omega_T}$ for each element $T$, which are not available on $\Gamma^{\omega_T}_h$. 
So we introduce $\tilde{q}_1:= u^+_E -  u^-_E$ and $\tilde{q}_2:= \beta^+ \nabla u^+_E - \beta^- \nabla u^-_E$ which are well-defined on the whole domain, 
and particularly, $\tilde{q}_1|_{\Gamma} = q_1$ and $\tilde{q}_2\cdot\bfn|_{\Gamma} = q_2$. 
The projections $\Pi^k_{\Gamma^{\omega_T}_h}\tilde{q}_j$, $j=1,2$, are then well-defined. 
In addition, as they are just polynomials,
we shall consider the trivial extension of $\Pi^k_{\Gamma^{\omega_T}_h} \tilde{q}_j$ to the whole space that has zero derivative along the direction perpendicular to $\Gamma^{\omega_T}_h$.
With a little abuse of the notation, we still employ the same notation $\Pi^k_{\Gamma^{\omega_T}_h} \tilde{q}$ for this extension, 
and it now satisfies $\nabla \Pi^k_{\Gamma^{\omega_T}_h} \tilde{q}_j \cdot \bar{\bfn} = 0$.
With the preparation, we present the following estimate:
\begin{lemma}
\label{q_est}
For each $T\in\mathcal{T}^i_h$ and $\tilde{q}\in H^{s+1/2}(\omega_T)$, 
there holds for $t=0$ with $s\ge 1/2$ and $t=1$ with $s\ge 3/2$ that
\begin{equation}
\label{q_est_eq0}
h^t_T \| \tilde{q} - \Pi^k_{\Gamma^{\omega_T}_h} \tilde{q} \|_{H^t(\Gamma^{\omega_T}_h)} \lesssim h^{\min\{k+1,s\}}_T \| \tilde{q} \|_{H^s(\Gamma^{\omega_T})} + h_T \| \tilde{q} \|_{H^1(\omega^{\text{int}}_T)} + h^{1+t}_T \|\tilde{q} \|_{H^{1+t}(\omega^{\text{int}}_T)},~~~ k=0,1.
\end{equation}
\end{lemma}
\begin{proof}
We need to be careful as the norm in the error bound is defined on $\Gamma^{\omega_T}$. The case of $t=0$ for the $L^2$ estimate is technical which is left to \ref{append1}. Here, we focus on $t=1$ which relies on the $L^2$ estimate. We only consider the case of $k \ge 1$, and the argument for $k=0$ is similar and simpler.

Let $\nabla_F$ be the surface gradient parallel to $\Gamma^{\omega_T}_h$. We introduce an auxiliary function $p_h \in \mathcal{P}_1(\Gamma^{\omega_T}_h)$ such that $\nabla_F p_h =  \Pi^{k-1}_{\Gamma^{\omega_T}_h} \nabla_F \tilde{q}$ and $\int_{\partial \Gamma^{\omega_T}_h} p_h \dd s = \int_{\partial \Gamma^{\omega_T}_h} \tilde{q} \dd s$. As $\Gamma^{\omega_T}_h$ is shape regular, i.e., it has an inscribed circle of radius $\mathcal{O}(h_T)$, the Poincar\'e inequality given by (2.15) in \cite{Brenner;Sung:2018Virtual} implies
\begin{equation}
\label{Ih_ext_eq-13_1}
\|  \tilde{q}  -  p_h \|_{L^2(\Gamma^{\omega_T}_h)} \lesssim h_T \| \nabla_F \tilde{q} - \Pi^{k-1}_{\Gamma^{\omega_T}_h} \nabla_F \tilde{q}\|_{L^2(\Gamma^{\omega_T}_h)} \lesssim h^{\min\{k,s-1\}+1}_T \| \nabla_F \tilde{q} \|_{H^{s-1}(\Gamma^{\omega_T})} + h^2_T \| \tilde{q} \|_{H^2(\omega^{\text{int}}_T)},
\end{equation}
where in the last inequality we have used the case of $t=0$.
Next, by the shape regularity of $\Gamma^{\omega_T}_h$ again, the inverse inequality, \eqref{Ih_ext_eq-13_1} and the case of $t=0$ again, we have
\begin{equation*}
\begin{split}
\label{Ih_ext_eq-14}
\| \nabla_F \Pi^k_{\Gamma^{\omega_T}_h}\tilde{q}  - \nabla_F p_h \|_{L^2(\Gamma^{\omega_T}_h)} & \lesssim h^{-1}_T \|  \Pi^k_{\Gamma^{\omega_T}_h}\tilde{q}  -  p_h \|_{L^2(\Gamma^{\omega_T}_h)}  \lesssim h^{-1}_T \|  \Pi^k_{\Gamma^{\omega_T}_h}\tilde{q}  -  \tilde{q} \|_{L^2(\Gamma^{\omega_T}_h)} + h^{-1}_T \|  \tilde{q}  -  p_h \|_{L^2(\Gamma^{\omega_T}_h)} \\
& \lesssim h^{\min\{k+1,s\}-1}_T \| \tilde{q} \|_{H^{s}(\Gamma^{\omega_T})} + \| \tilde{q} \|_{H^1(\omega^{\text{int}}_T)} + h^{\min\{k,s-1\}}_T \| \tilde{q} \|_{H^{s}(\Gamma^{\omega_T})} + h_T \| \tilde{q} \|_{H^2(\omega^{\text{int}}_T)}.
\end{split}
\end{equation*}
Combing the estimates above, we obtain
\begin{equation*}
\begin{split}
\label{Ih_ext_eq-15}
|  \Pi^k_{\Gamma^{\omega_T}_h}\tilde{q}  -  \tilde{q} |_{H^1(\Gamma^{\omega_T}_h)}  & \le \| \nabla_F \Pi^k_{\Gamma^{\omega_T}_h}\tilde{q}  - \Pi^{k-1}_{\Gamma^{\omega_T}_h} \nabla_F \tilde{q} \|_{L^2(\Gamma^{\omega_T}_h)} + \| \Pi^{k-1}_{\Gamma^{\omega_T}_h} \nabla_F \tilde{q}  - \nabla_F \tilde{q} \|_{L^2(\Gamma^{\omega_T}_h)} \\
& \lesssim h^{\min\{k,s-1\}}_T \| \tilde{q} \|_{H^{s}(\Gamma^{\omega_T})} + \| \tilde{q} \|_{H^1(\omega^{\text{int}}_T)} + h_T \| \tilde{q} \|_{H^2(\omega^{\text{int}}_T)},
\end{split}
\end{equation*}
which finishes the proof.
\end{proof}

\begin{remark}
The key of Lemma \ref{q_est} is to obtain the error bound in terms of the norm $\|\cdot\|_{H^s(\Gamma^{\omega_T})}$, 
which highly relies on the fact that $\Gamma^{\omega_T}$ is flat enough and close to $\Gamma^{\omega_T}_h$ locally. 
In fact, using the techniques above, we can also show
\begin{equation}
\label{remark_qest_eq0}
\| v \|_{L^2(\Gamma^{\omega_T}_h)} \lesssim \| v \|_{L^2(\Gamma^{\omega_T})} + h_T \| v \|_{H^1(\omega^{\text{int}}_T)}, ~~~ v\in H^1(\omega_T),
\end{equation}
meaning that the difference between $\| v \|_{L^2(\Gamma^{\omega_T}_h)}$ and $\| v \|_{L^2(\Gamma^{\omega_T})}$ is a high-order term.
\end{remark}

\subsection{A data-aware quasi-interpolation}

For each $\phi_T\in\mathcal{S}^0_h(T)$, we note that it admits the following representation:
\begin{equation}
\label{eq_phi0}
\phi_T|_{T^-} = p_h, ~~~~~~\phi_T|_{T^+} = p_h + \frac{\beta^- - \beta^+}{\beta^+} \nabla p_h\cdot \bar{\bfn} (X-F) \cdot \bar{\bfn},  ~~~~~ \forall p_h\in \mathcal{P}_1.
\end{equation}
One can verify that such a function satisfies the homogeneous jump conditions on $\Gamma^T_h$. This observation motivates the following operator:
\begin{equation}
\label{extension0}
\mathfrak{C}_T: \mathcal{P}_1 \rightarrow  \mathcal{P}_1 , ~~~ \text{with} ~\mathfrak{C}_T(p_h) = p_h + \frac{\beta^- - \beta^+}{\beta^+} \nabla p_h\cdot \bar{\bfn} (X-F) \cdot \bar{\bfn}.
\end{equation}
Then, we introduce a special interpolation operator on each patch $\omega_T$ denoted by $J_T:H^2(\omega^{+}_T\cup\omega^{-}_T)\rightarrow \mathcal{S}_h(T)$ satisfying
\begin{equation}
\label{Jh_def}
J_{T} u =
\left\{\begin{array}{cc}
J^-_{T} u := \Pi^1_{\omega_T} u^-_E  ~~~~~~~~~~~~~~~~~~~~~~~~~~~~~~~~~~~~~~~~~~~~~~~~~~~~~~~~~~~~~~~~~~~~~~ &\text{in} ~ \omega_T^- ,\\
J^+_{T} u := \mathfrak{C}_T (\Pi^1_{\omega_T} u^-_E) + \Pi^1_{\Gamma^{\omega_T}_h} \tilde{q}_1 + (\beta^+)^{-1} \Pi^0_{\Gamma^{\omega_T}_h}\tilde{q}_2 (X-F)\cdot\bar{\bfn} ~~~ &\text{in} ~ \omega_T^+.
\end{array}\right.
\end{equation}
Note that each polynomial component of an IFE function can be trivially and naturally extended to the whole space. 
In the following discussion, we shall use both of the two components $J^{\pm}_T u$ on the entire patch $\omega_T$. 
\begin{remark}
The quasi-interpolation does not rely on Lagrange-type shape functions; 
instead, it is similar to Hermite-type interpolation. 
Although we require an estimate for the Lagrange interpolation, the quasi-interpolation can serve as an intermediate step, facilitating the derivation of the desired estimate.
\end{remark}

Next, we introduce a special norm to gauge the error across the interface:
\begin{equation}
\label{triNorm}
 \vertiii{p_h}_{h,T} =  \| p_h \|_{L^2(\Gamma^{\omega_T}_h)} +  h_T \| \nabla p_h\cdot\bar{\mathbf{ n}} \|_{L^2(\Gamma^{\omega_T}_h)}, ~~~ v_h \in \mathcal{P}_1(\omega_T).
\end{equation}
\begin{lemma}
\label{lem_norm_equiv}
The following norm equivalence holds regardless of interface location.
\begin{equation}
\label{lem_norm_equiv_eq}
 \vertiii{\cdot}_{h,T} \simeq h^{-1/2}_T \| \cdot \|_{L^2(T)}, ~~~ \text{in} ~ \mathcal{P}_1(\omega_T).
\end{equation}
\end{lemma}
\begin{proof}
The argument is similar to Lemma 4.4 in \cite{2020GuoLin2}.
\end{proof}

\begin{theorem}
\label{thm_J_error}
Let $u\in H^2(\Omega^-\cup\Omega^+)$. Then, on each interface element $T$ with each patch $\omega_T$, there holds
\begin{equation}
\label{thm_J_error_eq0}
h^k_T| J^{\pm}_Tu - u^{\pm}_E |_{H^k(\omega_T)} \lesssim h^{2}_T ( \| u^-_E \|_{H^2(\omega_T)} + \| u^+_E \|_{H^2(\omega_T)} ) + h^{3/2}_T ( \| u^-_E \|_{H^1(\omega^{\text{int}}_T)} + \| u^+_E \|_{H^1(\omega^{\text{int}}_T)} ), ~~~ k = 0,1.
\end{equation}
\end{theorem}
\begin{proof}
The estimate for the ``-" piece is trivial. Let us only consider the ``$+$" component, and write down
\begin{equation}
\label{thm_J_error_eq1-1}
|J^{+}_Tu - u^{+}_E |_{H^k(\omega_K)} \le |J^{+}_Tu - \Pi^1_{\omega_T} u^{+}_E |_{H^k(\omega_K)} + | \Pi^1_{\omega_T} u^+_E  - u^{+}_E |_{H^k(\omega_K)}.
\end{equation}
The estimate of the second term in \eqref{thm_J_error_eq1-1} is trivial, and we focus on the first term. Denote $v_h = J^+_T u - \Pi^1_{\omega_T}u^+_E \in \mathcal{P}_1(\omega_T)$, and consider its norm in \eqref{lem_norm_equiv_eq}:
\begin{equation}
\label{thm_J_error_eq1}
\vertiii{v_h}_{h,T} = \| v_h \|_{L^2(\Gamma^{\omega_T}_h)} +  h_T \| \nabla v_h\cdot\bar{\mathbf{ n}} \|_{L^2(\Gamma^{\omega_T}_h)} =: (I) + (II).
\end{equation}
For $(I)$, since $(\bfx-F)\cdot\bar{\bfn}=0$ for $\bfx\in \Gamma^{\omega_T}_h$, we have the following decomposition:
\begin{equation}
\label{thm_J_error_eq2}
(I) \le \underbrace{ \| \mathfrak{C}_T (\Pi^1_{\omega_T} u^-_E)  - u^-_E \|_{L^2(\Gamma^{\omega_T}_h )} }_{(Ia)} + \underbrace{ \| u^-_E - \Pi^1_{\omega_T}u^+_E + \Pi^1_{\Gamma^{\omega_T}_h} \tilde{q}_1   \|_{L^2(\Gamma^{\omega_T}_h)} }_{(Ib)} .
\end{equation}
Then, the trace inequality of Lemma 3.4 in \cite{2016WangXiaoXu} implies
\begin{equation}
\label{thm_J_error_eq3}
(Ia) =  \| \Pi^1_{\omega_T} u^-_E  - u^-_E \|_{L^2(\Gamma^{\omega_T}_h)} \lesssim h^{-1/2}_T \| \Pi^1_{\omega_T} u^-_E  - u^-_E \|_{L^2(\omega_T)} + h^{1/2}_T | \Pi^1_{\omega_T} u^-_E  - u^-_E |_{H^1(\omega_T)} \lesssim h^{3/2} \| u^-_E \|_{H^2(\omega_T)}.
\end{equation}
As for $(Ib)$, by Lemma \ref{q_est} and the trace inequality of Lemma 3.4 in \cite{2016WangXiaoXu} again, we obtain
\begin{equation}
\begin{split}
\label{thm_J_error_eq4}
(Ib) & = \| u^+_E - \Pi^1_{\omega_T}u^+_E + \Pi^1_{\Gamma^{\omega_T}_h}\tilde{q}_1 - \tilde{q}_1  \|_{L^2(\Gamma^{\omega_T}_h)}  \le \| \Pi^1_{\omega_T}u^+_E -u^+_E   \|_{L^2(\Gamma^{\omega_T}_h)}
+ \|  \Pi^1_{\Gamma^{\omega_T}_h} \tilde{q}_1 - \tilde{q}_1  \|_{L^2(\Gamma^{\omega_T}_h)} \\
& \lesssim h^{3/2}_T \| u^+_E \|_{H^2(\omega_T)} + h^{3/2}_T |q_1|_{H^{3/2}(\Gamma^{\omega_T})} + h_T (  \| u^+_E \|_{H^1(\omega^{\text{int}}_T)} +  \| u^-_E \|_{H^1(\omega^{\text{int}}_T)}).
\end{split}
\end{equation}
Putting \eqref{thm_J_error_eq3} and \eqref{thm_J_error_eq4} into \eqref{thm_J_error_eq2}, we get the estimate for $(I)$

As for $(II)$, denoting $\rho=\beta^-/\beta^+$ for simplicity, we use the triangular inequality to write down
\begin{equation}
\begin{split}
\label{thm_J_error_eq5}
(II) \leqslant  h_T \underbrace{ \| \nabla (\mathfrak{C}_T (\Pi^1_{\omega_T} u^-_E)  - \rho u^-_E)\cdot\bar{\bfn} \|_{L^2(\Gamma^{\omega_T}_h)} }_{(IIa)}   
+ h_T \underbrace{ \| (\beta^+)^{-1} \Pi^0_{\Gamma^{\omega_T}_h}\tilde{q}_2 + \nabla ( \rho u^-_E - \Pi^1_{\omega_T}u^+_E)\cdot\bar{\bfn}  \|_{L^2(\Gamma^{\omega_T}_h )} }_{(IIb)},
\end{split}
\end{equation}
where we have used $\nabla \Pi^1_{\Gamma^{\omega_T}_h}\tilde{q}_1 \cdot\bar\bfn =0$.
The estimate of $(IIa)$ is similar to \eqref{thm_J_error_eq3}:
\begin{equation}
\label{thm_J_error_eq6}
(IIa) =  \rho \| \nabla ( \Pi^1_{\omega_T} u^-_E  -  u^-_E)\cdot\bar{\bfn} \|_{L^2(\Gamma^{\omega_T}_h)}\lesssim h^{1/2} \| u^-_E \|_{H^2(\omega_T)}.
\end{equation}
The estimate of $(IIb)$ is similar to \eqref{thm_J_error_eq4}:
\begin{equation}
\begin{split}
\label{thm_J_error_eq7}
(IIb) & \le \| (\beta^+)^{-1} \Pi^0_{\Gamma^{\omega_T}_h}\tilde{q}_2 - (\beta^+)^{-1} \tilde{q}_2  \|_{L^2(\Gamma^{\omega_T}_h)} +\|  \nabla (  u^+_E - \Pi^1_{\omega_T}u^+_E)\cdot\bar{\bfn}  \|_{L^2(\Gamma^{\omega_T}_h)} \\
& \lesssim  h^{1/2}_T ( \| u^-_E \|_{H^2(\omega_T)} + \| u^+_E \|_{H^2(\omega_T)} ) .
\end{split}
\end{equation}
Now, putting \eqref{thm_J_error_eq6} and \eqref{thm_J_error_eq7} into \eqref{thm_J_error_eq5}, we have the estimate for $(II)$. 
Combining the estimates for $(I)$ and $(II)$, we have the estimate for $\vertiii{v_h}_{h,T}$. 
Then, the desired result in terms of the $L^2$ norm comes from the norm equivalence by Lemma \ref{lem_norm_equiv}. 
The estimate in terms of the $H^1$ norm follows from the inverse inequality for $|J^{+}_Tu - \Pi^1_{\omega_T} u^{+}_E |_{H^1(\omega_K)}$.
\end{proof}



\subsection{Lagrange-type interpolation}
On a tetrahedron $T=A_1A_2A_3A_4$, we define the IFE version of the nodal-value interpolation incorporating non-homogeneous jumps:
\begin{equation}
\label{interp}
I_T: H^2(T^-\cup T^+)  \rightarrow \mathcal{S}_h(T), ~~~ \text{with}~ I_T u = \sum_{j=1}^4 u(A_j) \phi_{j,T} + q_{T,1}(u) + q_{T,2}(u)
\end{equation}
where $q_{T,1}(u)$ and $q_{T,2}(u)$ are defined in \eqref{qh1} and \eqref{qh2}. Then, the global interpolation is defined as
\begin{equation}
\label{interp_glob}
I_h: H^2(\Omega^-\cup\Omega^+)  \rightarrow \mathcal{S}_h, ~~~ \text{with}~ I_h u|_T = I_T u, ~ \forall T\in \mathcal{T}_h,
\end{equation}
where $I_T u$ is given as \eqref{interp} on interface elements and the standard Lagrange interpolation on non-interface elements.
\begin{remark}
\label{rem_IT}
Different from the case of homogeneous jumps, $I_T$ does not reduce to an identity on $\mathcal{S}_h(T)$, i.e., 
$I_T$ is not reproductive, since $q_{T,2}$ is defined with the actual interface. Delicate analysis techniques are demanded for this issue. 
\end{remark}


The error estimation relies on the following decomposition:
\begin{equation}
\label{interp_decomp}
u - I_T u  =  (u - J_T u ) + (J_Tu - I_TJ_Tu ) + ( I_TJ_Tu - I_Tu ) =: \xi_1 + \xi_2 + \xi_3. 
\end{equation}
Note that the estimate of $\xi_1$ directly follows from Lemma \ref{thm_J_error}. So we proceed to estimate the remaining two terms.
For $\xi_2$, the main difficulty comes from the issue in Remark \ref{rem_IT}, but we can show that $I_T v_h - v_h$ produces optimal errors.

\begin{lemma}
\label{lem_ITv_est}
Let $u\in H^2(\Omega^-\cup\Omega^+)$. Then, on each interface element $T$ with each patch $\omega_T$, there holds for each $v_h\in \mathcal{S}_h(T)$ that
\label{lem_ITv}
\begin{equation}
\label{lem_ITv_est_eq0}
|\xi_2|_{H^k(\omega_T)} \lesssim h^{-k+2}_T( \| u^{+}_E\|_{H^2(\omega_T)} + \| u^-_E\|_{H^2(\omega_T)} ), ~~ k=0,1.
\end{equation}
\end{lemma}
\begin{proof}
Given each $v_h\in V_h(T)$, we first show
\begin{equation}
\label{lem_ITv_est_eq1}
| I_T v_h  - v_h|_{H^k(\omega_T)} \lesssim h^{-k+2}_T( \|\nabla v^{+}_h\|_{L^2(\omega_T)} + \|\nabla v^-_h\|_{L^2(\omega_T)} ).
\end{equation}
By the identity \eqref{interp_1} and the Mean Value Theorem, there exists $\tilde{\bfx}\in\Gamma^T$ such that
\begin{equation}
\begin{split}
I_{T}v_h - v_h = ( \beta^+\nabla v^+_h -  \beta^-\nabla v^-_h) \cdot (\bfn(\tilde{\bfx}) - \bar{\bfn})  \xi_4, 
\end{split}
\end{equation}
where we have used $\nabla v^{\pm}_h$ being constant vectors. 
Then, \eqref{lem_ITv_est_eq1} follows from \eqref{lem_geo_gamma_eq02} together with Theorem \ref{bounds_IFEshapeFun}:
\begin{equation}
\begin{split}
\label{IhT_diff}
| I_{T}v_h - v_h |_{H^k(\omega_T)} &\lesssim h^{7/2-k}_T  ( \|\nabla v^{+}_h\| + \|\nabla v^-_h\| ) \lesssim 
h^{-k+2}_T( \|\nabla v^{+}_h\|_{L^2(\omega_T)} + \|\nabla v^-_h\|_{L^2(\omega_T)} ).
\end{split}
\end{equation}
Now, taking $v_h = J_Tu$ and applying the boundedness property of $J_Tu$ given by Theorem \ref{thm_J_error}, we arrive at \eqref{lem_ITv_est_eq0}.
\end{proof}

Next, we estimate $\xi_3 = I_T(J_Tu - u)$. 
Define $\eta_T = J_T u - u$ and notice that each $\eta^{\pm}_T = J^{\pm}_T u - u^{\pm}_E$ can be naturally used on the whole $\omega_T$. 
By the definition \eqref{interp}, one key is to estimate $q_{T,1}(\eta_T)$ and $q_{T,2}(\eta_T)$ given by the following lemma.
\begin{lemma}
\label{lem_IqT_est}
Let $q_1\in H^{3/2}(\Gamma)$ and $q_2\in H^{1/2}(\Gamma)$. Then, on each interface element $T$ with each patch $\omega_T$, there hold
\begin{subequations}
\label{lem_IqT_est_eq0}
\begin{align}
    & | {q}_{T,1}(\eta_T)|_{H^k(\omega_T)} \lesssim h^{2-k}_T ( \| q_1 \|_{H^{3/2}(\Gamma^{\omega_T})} +  \| u^-_E \|_{H^2(\omega_T)} + \| u^+_E \|_{H^2(\omega_T)}  )+ h^{3/2-k}_T( \| u^-_E \|_{H^1(\omega^{\text{int}}_T)} + \| u^+_E \|_{H^1(\omega^{\text{int}}_T)} ), \label{lem_IqT_est_eq01} \\
    &  |q_{T,2}(\eta_T)|_{H^k(\omega_T)} \lesssim h^{2-k}_T ( \| q_2 \|_{H^{1/2}(\Gamma^{\omega_T})} +  \| u^-_E \|_{H^2(\omega_T)} + \| u^+_E \|_{H^2(\omega_T)} ), ~~~~ k=0,1. \label{lem_IqT_est_eq02} 
\end{align}
\end{subequations}
\end{lemma}
\begin{proof}
The jump conditions of $u$ and $J_Tu$ together with $(D_j-F)\cdot\bar{\bfn}=0$, $j=1,2,3$, yield  
\begin{subequations}
\label{Ih_ext_eq-1}
\begin{align}
 &  \jump{\eta_T}_{D_j} = \Pi^1_{\Gamma^{\omega_T}_h}\tilde{q}_1(D_j)  - {q}_1(D_j), ~ j=1,2,3, \label{Ih_ext_eq-11} \\   
 &  \jump{\beta \nabla \eta_T\cdot \bfn}_{\Gamma^{\omega_T}} = (\beta^+ - \beta^-) ( \nabla \Pi^1_{\omega_T} u^-_E\cdot \bfn - \bar{\bfn}\cdot\bfn \nabla \Pi^1_{\omega_T} u^-_E\cdot \bar{\bfn}  )+ \beta^+ \nabla  \Pi^1_{\Gamma^{\omega_T}_h} \tilde{q}_1\cdot\bfn + \Pi^0_{\Gamma^{\omega_T}_h} \tilde{q}_2 \bfn\cdot\bar{\bfn} - q_2. \label{Ih_ext_eq-12}
\end{align}
\end{subequations}

Let us first estimate $q_{T,1}(\eta_T)$ from \eqref{Ih_ext_eq-11}. 
As $\Gamma^{\omega_T}_h$ and $\omega_T$ are both shape-regular, we can find a shape regular tetrahedron $K\subseteq\omega_T$ 
such that one of its faces denoted by $F$ also locates inside $\Gamma^{\omega_T}_h$, i.e., $F\subseteq \Gamma^{\omega_T}_h$. 
With $\tilde{q}(D_j)=q(D_j)$, $j=1,2,3$, by Sobolev embedding theorem, we have
\begin{equation}
\label{Ih_ext_eq-13}
|\Pi^1_{\Gamma^{\omega_T}_h}\tilde{q}_1(D_j)  - \tilde{q}_1(D_j)| \lesssim h^{-1}_T \| \Pi^1_{\Gamma^{\omega_T}_h}\tilde{q}_1  - \tilde{q}_1 \|_{L^2(F)}
+  | \Pi^1_{\Gamma^{\omega_T}_h}\tilde{q}_1  - \tilde{q}_1 |_{H^1(F)} + h^{1/2}_T | \Pi^1_{\Gamma^{\omega_T}_h}\tilde{q}_1  - \tilde{q}_1 |_{H^{3/2}(F)}.
\end{equation}
The estimates of the first and second terms above immediately follow from Lemma \ref{q_est}. As $K$ is shape regular, the estimate for the third term can be derived from the trace inequality:
\begin{equation}
\label{Ih_ext_eq-13_2}
h^{1/2}_T | \Pi^1_{\Gamma^{\omega_T}_h}\tilde{q}_1  - \tilde{q}_1 |_{H^{3/2}(F)} = h^{1/2}_T | \tilde{q}_1 |_{H^{3/2}(F)} \lesssim h^{1/2}_T \| \tilde{q}_1 \|_{H^2(K)} \lesssim h^{1/2}_T \| \tilde{q}_1 \|_{H^2(\omega_T)}.
\end{equation}
Then, using \eqref{qh1} with Theorem \ref{bounds_IFEshapeFun}, we have \eqref{lem_IqT_est_eq01} by the definition of $\tilde{q}_1$.

Next, we estimate $q_{T,2}(\eta_T)$.
By H\"older's inequality, we have
\begin{equation}
\begin{split}
\label{Ih_ext_eq-2}
&|\Gamma^{\omega_T}|^{-1} \int_{\Gamma^{\omega_T}} \jump{\beta \nabla \eta_T\cdot\bfn}_{\Gamma^{\omega_T}} \dd s \\
 \lesssim & \underbrace{  |\Gamma^{\omega_T}|^{-1/2} \| \nabla \Pi^1_{\omega_T} u^-_E\cdot \bfn - \bar{\bfn}\cdot\bfn \nabla \Pi^1_{\omega_T} u^-_E\cdot \bar{\bfn}  \|_{L^2(\Gamma^{\omega_T})} }_{(I)} + \underbrace{ |\Gamma^{\omega_T}|^{-1/2} \| \beta^+ \nabla  \Pi^1_{\Gamma^{\omega_T}_h} \tilde{q}_1\cdot\bfn \|_{L^2(\Gamma^{\omega_T})} }_{(II)} \\
+& \underbrace{ |\Gamma^{\omega_T}|^{-1/2} \| \Pi^0_{\Gamma^{\omega_T}_h} \tilde{q}_2 - q_2 \|_{L^2(\Gamma^{\omega_T})} }_{(III)} + \underbrace{ |\Gamma^{\omega_T}|^{-1/2} \| \Pi^0_{\Gamma^{\omega_T}_h} \tilde{q}_2 (1-\bfn\cdot\bar{\bfn}) \|_{L^2(\Gamma^{\omega_T})} }_{(IV)}.
\end{split}
\end{equation} 
For $(I)$, by \eqref{lem_geo_gamma_eq02}, the trace inequality and the boundedness of $\Pi^1_{\omega_T}$, we have
\begin{equation}
\label{Ih_ext_eq-2_1}
(I) =  |\Gamma^{\omega_T}|^{-1/2} \| \nabla \Pi^1_{\omega_T} u^-_E\cdot (\bfn - \bar{\bfn}) - (\bar{\bfn}\cdot\bfn-1) \nabla \Pi^1_{\omega_T} u^-_E\cdot \bar{\bfn}  \|_{L^2(\Gamma^{\omega_T})}  \lesssim \|  \nabla \Pi^1_{\omega_T} u^-_E \|_{L^2(\Gamma^{\omega_T})}
\lesssim h^{-1/2}_T \|  \nabla u^-_E \|_{L^2(\omega_T)}.
\end{equation}
For $(II)$, by applying \eqref{lem_geo_gamma_eq02} and \eqref{remark_qest_eq0}, the boundedness of $\Pi^1_{\Gamma^{\omega_T}_h}$ and noticing that $\nabla  \Pi^1_{\Gamma^{\omega_T}_h}\cdot$ is a constant, we have
\begin{equation}
\begin{split}
\label{Ih_ext_eq-3}
(II) & = |\Gamma^{\omega_T}|^{-1/2} \| \nabla  \Pi^1_{\Gamma^{\omega_T}_h} \tilde{q}_1\cdot (\bfn - \bar{\bfn}) \|_{L^2(\Gamma^{\omega_T})}   \lesssim \| \nabla  \Pi^1_{\Gamma^{\omega_T}_h} \tilde{q}_1 \|_{L^2(\Gamma^{\omega_T}_h)} \lesssim \| \nabla \tilde{q}_1 \|_{L^2(\Gamma^{\omega_T}_h)}  \lesssim \| q_1 \|_{H^1(\Gamma^{\omega_T})} + h_T \| \tilde{q}_1 \|_{H^2(\omega_T)}. 
\end{split}
\end{equation}
Next, for $(III)$, we apply Lemma \ref{q_est} to obtain
\begin{equation}
\label{Ih_ext_eq-4}
(III) \le |\Gamma^{\omega_T}|^{-1/2}  \| \tilde{q}_2 - \Pi^0_{\Gamma^{\omega_T}_h} \tilde{q}_2\|_{L^2(\Gamma^{\omega_T}_h)}\lesssim h^{-1/2}_T \| {q}_2 \|_{H^{1/2}(\Gamma^{\omega_T})} +  \| \tilde{q}_2 \|_{H^1(\omega^{\text{int}}_T)} .
\end{equation}
As for $(IV)$, using the similar argument to \eqref{Ih_ext_eq-3}, we have
\begin{equation}
\label{Ih_ext_eq-5}
(IV) \lesssim h_T \| \Pi^0_{\Gamma^{\omega_T}_h} \tilde{q}_2  \|_{L^2(\Gamma^{\omega_T}_h)} \lesssim h_T \| \tilde{q}_2 \|_{L^2(\Gamma^{\omega_T}_h)}
\lesssim h_T \| q_2 \|_{L^2(\Gamma^{\omega_T})} + h^2_T \| \tilde{q}_2 \|_{H^1(\omega^{\text{int}}_T)}.
\end{equation}
Noticing $\|\tilde{q}_1\|_{H^2(\omega_T)} \simeq \|\tilde{q}_2\|_{H^1(\omega_T)}$, 
using \eqref{qh2} with Theorem \ref{bounds_IFEshapeFun} and putting \eqref{Ih_ext_eq-2_1}-\eqref{Ih_ext_eq-5} to \eqref{Ih_ext_eq-2}, 
we obtain \eqref{lem_IqT_est_eq02} by the definition of $\tilde{q}_1$ and $\tilde{q}_2$.
\end{proof}

\begin{remark}
\label{rem_IqT_est}
Note that the polynomial order used in the approximation of $q_1$ and $q_2$ are $1$ and $0$, respectively. 
Thus, if $q_1$, $q_2$ and their associated $\tilde{q}_1$ and $\tilde{q}_2$ have higher regularity, 
the convergence order will also be enhanced correspondingly. 
In fact, using the arguments above, we can show the following result:
if $q_1\in H^{s}(\Gamma)$, $\tilde{q}_1\in H^{s+1/2}(\omega_T)$, $s\ge 3/2,$ and $q_2\in H^{t}(\Gamma)$, $\tilde{q}_2\in H^{t+1/2}(\omega_T)$, $t \ge 1/2$, 
then there holds for $k=0,1$
\begin{subequations}
\label{rem_IqT_est_eq0}
\begin{align*}
    & | {q}_{T,1}(\eta_T)|_{H^k(\omega_T)} \lesssim h^{1/2+\min\{s,2\}-k}_T \| q_1 \|_{H^{s}(\Gamma^{\omega_T})} + h^{2-k}_T (   \| u^-_E \|_{H^2(\omega_T)} + \| u^+_E \|_{H^2(\omega_T)}  )+ h^{3/2-k}_T( \| u^-_E \|_{H^1(\omega^{\text{int}}_T)} + \| u^+_E \|_{H^1(\omega^{\text{int}}_T)} ), \\
    &  |q_{T,2}(\eta_T)|_{H^k(\omega_T)} \lesssim h^{3/2+\min\{t,1\}-k}_T  \| q_2 \|_{H^{t}(\Gamma^{\omega_T})} + h^{2-k}_T (  \| u^-_E \|_{H^2(\omega_T)} + \| u^+_E \|_{H^2(\omega_T)} ) . 
\end{align*}
\end{subequations}
\end{remark}

\begin{lemma}
\label{thm_Ih_ext}
Let $u\in H^2(\Omega^-\cup\Omega^+)$. Then, on each interface element $T$ with each patch $\omega_T$, there holds 
\begin{equation}
\label{Ih_est_eq0}
|\xi_3 |_{H^k(\omega_T)} \lesssim h^{2-k} ( \| u^-_E \|_{H^2(\omega_T)} + \| u^+_E \|_{H^2(\omega_T)} + \| q_1 \|_{H^{3/2}(\Gamma^{\omega_T})} + \| q_2 \|_{H^{1/2}(\Gamma^{\omega_T})} ) + h^{3/2-k}_T ( \| u^-_E \|_{H^1(\omega^{\text{int}}_T)} + \| u^+_E \|_{H^1(\omega^{\text{int}}_T)} ), ~ k=0,1.
\end{equation}
\end{lemma}
\begin{proof}
Assume the vertex $A_j$ is at the piece $T^{s_j}$, with $s_j=\pm$, $j=1,2,3,4$. 
 By Sobolev embedding Theorem, the scaling argument, Theorems \ref{bounds_IFEshapeFun} and \ref{thm_J_error}, we have
\begin{equation}
\begin{split}
\label{Ih_ext_eq1}
|\phi_{j,T}|_{H^{k}(\omega_T)} |\eta_T(A_j)| & \lesssim h^{3/2-k}_T \left( h^{-3/2}_T \|\eta_T^{s_j} \|_{L^2(\omega_T)} +  h^{-1/2}_T |\eta_T^{s_j} |_{H^1(\omega_T)}  +  h^{1/2}_T |\eta_T^{s_j} |_{H^2(\omega_T)} \right) \\
& \lesssim h^{2-k} ( \| u^-_E \|_{H^2(\omega_T)} + \| u^+_E \|_{H^2(\omega_T)} ) + h^{3/2-k}_T ( \| u^-_E \|_{H^1(\omega^{\text{int}}_T)} + \| u^+_E \|_{H^1(\omega^{\text{int}}_T)} ).
\end{split}
\end{equation}
Then, by Lemma \ref{lem_IqT_est} we have
\begin{equation*}
\begin{split}
\label{Ih_ext_eq2}
| {I}_{T} \eta_T|_{H^k(\omega_T)} &\leqslant \sum_{j=1,2,3,4} |\phi_{j,T}|_{H^{k}(\omega_T)} |\eta_T(A_j)| + |{q}_{T,1}(\eta_T)|_{H^k(\omega_T)} + |{q}_{T,2}(\eta_T)|_{H^k(\omega_T)} \\
& \lesssim h^{2-k} ( \| u^-_E \|_{H^2(\omega_T)} + \| u^+_E \|_{H^2(\omega_T)} + \| q_1 \|_{H^{3/2}(\Gamma^{\omega_T})} + \| q_2 \|_{H^{1/2}(\Gamma^{\omega_T})} ) \\ 
&+ h^{3/2-k}_T ( \| u^-_E \|_{H^1(\omega^{\text{int}}_T)} + \| u^+_E \|_{H^1(\omega^{\text{int}}_T)} ),
\end{split}
\end{equation*}
which finishes the proof.
\end{proof}

Combining the estimates above, we finally have the following result.
\begin{lemma}
\label{thm_interp_loc_est}
Let $u\in H^2(\Omega^-\cup\Omega^+)$. Then, on each interface element $T$ with each patch $\omega_T$, there holds
\begin{equation}
\label{thm_interp_loc_est_eq0}
|I_T u - u|_{H^k(\omega_T)} \lesssim h^{2-k}_T ( \| u^-_E \|_{H^2(\omega_T)} + \| u^+_E \|_{H^2(\omega_T)} + \| q_1 \|_{H^{3/2}(\Gamma^{\omega_T})} + \| q_2 \|_{H^{1/2}(\Gamma^{\omega_T})} ) + h^{3/2-k}_T ( \| u^-_E \|_{H^1(\omega^{\text{int}}_T)} + \| u^+_E \|_{H^1(\omega^{\text{int}}_T)} ) .
\end{equation}
\end{lemma}
\begin{proof}
The result directly follows from applying Theorem \ref{thm_J_error} and Lemma \ref{thm_Ih_ext} to \eqref{interp_decomp}.
\end{proof}


\begin{lemma}
\label{lem_interp_enrg_error}
Let $u\in H^2(\Omega^-\cup\Omega^+)$. Then, there holds that
\begin{equation}
\label{lem_interp_enrg_error_eq0}
\vertiii{u - I_hu} \lesssim h ( \|u\|_{H^2(\Omega^+\cup\Omega^-)} +  \| q_1 \|_{H^{3/2}(\Gamma)} + \| q_2 \|_{H^{1/2}(\Gamma)} ).
\end{equation}
\end{lemma}
\begin{proof}
The standard estimate of the Lagrange interpolation and Theorem \ref{thm_interp_loc_est} give
\begin{equation*}
\label{lem_interp_enrg_error_eq1}
\sum_{T\in\mathcal{T}_h} \| \nabla (u - I_hu)\|^2_{L^2(T)}  \lesssim h^2_T ( \| u^+_E \|^2_{H^2(\Omega)} +  \|u^-_E\|^2_{H^2(\Omega)} + \| q_1 \|^2_{H^{3/2}(\Gamma)} + \| q_2 \|^2_{H^{1/2}(\Gamma)} ) + h_T ( \| u^-_E \|^2_{H^1(S_{\delta})} + \| u^+_E \|^2_{H^1(S_{\delta} )})  ,
\end{equation*}
where we have used the finite overlapping property of the patches $\omega_T$, $T\in\mathcal{T}^i_h$, and $S_{\delta}$ is the $\delta$-strip given by \eqref{delta_strip_1}. For the second term in \eqref{energy_norm}, for each $F\in\mathcal{F}^i_h$ and one of its associated element $T_F$, there exist pyramids $P^{\pm}\subseteq \omega^{\pm}_{T_F}$ which have the base $F^{\pm}$ such that their heights are in the order of $\mathcal{O}(h_T)$. Then, the trace inequality given by \cite[Lemma 3.2]{2018ChenCao} and Theorem \ref{thm_interp_loc_est} lead to 
\begin{equation}
\begin{split}
\label{lem_interp_enrg_error_eq2}
& \sqrt{\sigma} h^{-1/2}_T \|  u - I_hu \|_{L^2(F^{\pm})}  \lesssim  h^{-1}_T \| u - I_hu \|_{L^2(P^{\pm})} + | u - I_hu |_{H^1(P^{\pm})}  \\
 \lesssim & h_T ( \| u^+_E  \|_{H^2(\omega_{T_F})} + \| u^-_E  \|_{H^2(\omega_{T_F})} + \| q_1 \|_{H^{3/2}(\Gamma^{\omega_T})} + \| q_2 \|_{H^{1/2}(\Gamma^{\omega_T})}  )  + h^{1/2}_T ( \| u^-_E \|_{H^1(\omega^{\text{int}}_T)} + \| u^+_E \|_{H^1(\omega^{\text{int}}_T)} ).
\end{split}
\end{equation}
The analysis of the third term in \eqref{energy_norm} is similar to \eqref{lem_interp_enrg_error_eq2}.
Summing these estimates above over all the interface faces, we have the global estimates. 
It remains to estimate the norm on $S_{\delta}$. We apply Lemma \ref{lem_delta} to obtain
\begin{equation}
\label{lem_interp_enrg_error_eq4}
\| u^{\pm}_E \|_{H^1(S_{\delta})} \lesssim \sqrt{\delta} \| u^{\pm}_E \|_{H^2(\Omega)} \lesssim  h_T \| u^{\pm}_E \|_{H^2(\Omega^{\pm})},
\end{equation}
where the last inequality follows from \eqref{sobolev_ext}.
\end{proof}


\section{Scheme analysis}
\label{sec:scheme}

In this section, we analyze the enriched IFE scheme \eqref{IFE_scheme} including the convergence order and conditioning.

\subsection{Error analysis}

Let us first recall the trace and inverse inequalities in the following lemma.

\begin{theorem}[Trace and Inverse Inequalities]
\label{thm_trace_inequa}
On each interface element $T$ and its face $F$, the following trace and inverse inequalities hold for $v_h\in \mathcal{S}^0_h(T)$
\begin{subequations}
\label{thm_trace_inequa_eq0}
\begin{align}
    &  \| \nabla v_h\cdot\mathbf{ n} \|_{L^2(F)} \lesssim  h^{-1/2} \| \nabla v_h \|_{L^2(T)}, \label{thm_trace_inequa_eq01}  \\
    & \| \nabla v_h \|_{L^2(T)} \lesssim  h^{-1} \| v_h \|_{L^2(T)}.  \label{thm_trace_inequa_eq02}
\end{align}
\end{subequations}
\end{theorem}
\begin{proof}
As the linear homogeneous IFE space is a subspace of the bilinear IFE space in \cite{2020GuoZhang}, 
the results are immediately given by Theorems 4.5 and 4.11 there.
\end{proof}

The following theorem gives the continuity and coercivity.
\begin{lemma}
\label{thm_bound}
Assume $\sigma$ is large enough. Then, there hold
\begin{subequations}
\begin{align}
    &  a_h(v,w) \lesssim \vertiii{v}  \vertiii{w}, ~~~~ \forall v,w \in \mathcal{V}_h(\Omega), \label{thm_bound_0}  \\
    &  a_h(v,v) \gtrsim \vertiii{v}^2, ~~~~ \forall v\in \mathcal{S}_h(\Omega).  \label{thm_coer_eq0}
\end{align}
\end{subequations}
\end{lemma}
\begin{proof}
\eqref{thm_bound_0} follows from the H\"older's inequality, while \eqref{thm_coer_eq0} directly follows from the inverse inequality \eqref{thm_trace_inequa_eq01} through the standard argument.
\end{proof}



Next, we estimate the inconsistent error. 
Define $\Omega^{\Gamma}_h = \cup_{T\in \mathcal{T}^i_h}T$ as the subregion formed by all the interface elements.
\begin{lemma}
\label{lem_ErrEqn}
Let $u\in H^2(\Omega^-\cup\Omega^+)$ be the solution to \eqref{model}. Then, $\forall v_h\in \mathcal{S}^0_h(\Omega)$ there holds
\begin{equation}
\label{lem_ErrEqn_eq0}
\verti{ a_h(u,v_h) - a_h(u_h,v_h) } \lesssim h^{3/2} \| u \|_{H^2(\Omega^-\cup\Omega^+)}  \| \nabla v_h \|_{L^2(\Omega^{\Gamma}_h)}  .
\end{equation}
\end{lemma}
\begin{proof}
We first prove the following identity
\begin{equation}
\begin{split}
\label{lem_ErrEqn_eq01}
a_h(u,v_h) - a_h(u_h,v_h) = b_h(u,v_h):= \int_{\Gamma}  \{\beta \nabla u\cdot \mathbf{ n}\} \jump{v_h} ds.
 \end{split}
\end{equation}
We test \eqref{inter_PDE} by $v_h\in \mathcal{S}^0_h$. On an interface element, the integration by parts yields
\begin{equation*}
\begin{split}
\label{lem_ErrEqn_eq1}
 \int_T -\nabla\cdot(\beta \nabla u) v_h dx& = \int_T \beta \nabla u\cdot \nabla v_h dx - \int_{\partial T^-} (\beta\nabla u \cdot \bfn) v_h ds - \int_{\partial T^+} (\beta\nabla u \cdot \bfn) v_h ds \\
 & = \int_T \beta \nabla u\cdot \nabla v_h dx - \int_{\partial T} (\beta\nabla u \cdot \bfn) v_h ds - \int_{\Gamma^T} \{ \beta \nabla u \cdot \bfn \} \jump{v_h} ds
 - \int_{\Gamma^T} \jump{ \beta \nabla u \cdot \bfn } \{ v_h\} ds.
 \end{split}
\end{equation*}
The identities on non-interface elements are trivial. Now, adding these identities for all elements, we achieve
\begin{equation*}
\begin{split}
\label{lem_ErrEqn_eq2}
 \sum_{T\in\mathcal{T}_h} \int_T \beta \nabla u \cdot \nabla v_h dX 
 - \sum_{F\in \mathcal{F}^i_h} \int_F \aver{ \beta \nabla u\cdot \mathbf{ n} } \jump{v_h} ds  
 -  \sum_{T\in\mathcal{T}^i_h} \int_{\Gamma^T} \{  \beta \nabla u\cdot \mathbf{ n} \} \jump{v_h} ds = \int_{\Omega} fv_h dX + \int_{\Gamma^T} q_2 \{ v_h\} ds.
 \end{split}
\end{equation*}
Then, the identity \eqref{lem_ErrEqn_eq01} follows from $\jump{u}=0$ on each non-interface face $F$.
Next, let us estimate $b_h(u,v_h)$. Using \eqref{thm_trace_interface_eq0} and H\"older's inequality, we have 
\begin{equation*}
\begin{split}
\label{thm_error_bound_eq4}
|b_h(u,v_h)| 
& \lesssim \sum_{T\in\mathcal{T}^i_h} ( \| \nabla u^-\cdot \mathbf{ n} \|_{L^2(\Gamma^T)}  +  \| \nabla u^+\cdot \mathbf{ n} \|_{L^2(\Gamma^T)} ) h^{3/2}_T \|  \nabla v_h \|_{L^2(T)} \\
& \lesssim h^{3/2} ( \| \nabla u^-\cdot \mathbf{ n} \|_{L^2(\Gamma)}  +  \| \nabla u^+\cdot \mathbf{ n} \|_{L^2(\Gamma)} )  \| \nabla v_h \|_{L^2(\cup\mathcal{T}^i_h)} \\
& \lesssim h^{3/2}  \| u \|_{H^{2}(\Omega^-\cup\Omega^+)}   \| \nabla v_h \|_{L^2(\Omega^{\Gamma}_h)},
\end{split}
\end{equation*}
where we have used the discrete H\"older's inequality and the trace inequality on the whole domain.
\end{proof}

With the preparation above, we are able to establish the error estimates for the IFE solutions.

\begin{theorem}
\label{thm_error_bound}
Let $u\in H^2(\Omega^-\cup\Omega^+)$ be the solution to \eqref{model}. Then there holds
\begin{equation}
\label{thm_error_bound_eq0}
\vertiii{u - u_h} \lesssim h  ( \|u\|_{H^2(\Omega^+\cup\Omega^-)} +  \| q_1 \|_{H^{3/2}(\Gamma)} + \| q_2 \|_{H^{1/2}(\Gamma)} ).
\end{equation}
\end{theorem}
\begin{proof}
Note that $u_h - I_h u \in \mathcal{S}^0_h$ as the enrichment are cancelled. Then, by Lemmas \ref{thm_bound}, \ref{lem_ErrEqn} and \ref{lem_interp_enrg_error}, we have
\begin{equation*}
\begin{split}
\label{thm_error_bound_eq1}
\vertiii{ u_h - I_h u }^2 & \lesssim a_h( u_h - I_h u, u_h - I_h u ) = a_h( u - I_h u, u_h - I_h u ) - b_h(u, u_h - I_h u) \\
& \lesssim h ( \|u\|_{H^2(\Omega^+\cup\Omega^-)} +  \| q_1 \|_{H^{3/2}(\Gamma)} + \| q_2 \|_{H^{1/2}(\Gamma)} ) \vertiii{ u_h - I_h u } + h^{3/2} \| u \|_{H^2(\Omega^-\cup\Omega^+ )} \vertiii{  u_h - I_h u },
\end{split}
\end{equation*}
which finishes the proof by cancelling one term of $\vertiii{ u_h - I_h u }$.
\end{proof}
%

The error estimate for the $L^2$ is more difficult, which relies on more regular jump data.

\begin{theorem}
\label{thm_l2_error_bound}
Let $u\in H^2(\Omega^-\cup\Omega^+)$ be the solution to \eqref{model}, and assume $q_1\in H^{5/2}(\Gamma)$. Then, there holds
\begin{equation}
\label{thm_l2_error_bound_eq0}
\| u - u_h \|_{L^2(\Omega)} \lesssim h^2( \| f \|_{L^2(\Omega)} + \|q_1 \|_{H^{5/2}(\Gamma)} + \| q_2 \|_{H^{1/2}(\Gamma)} ) .
\end{equation}
\end{theorem}
\begin{proof}
We follow the duality argument. Define an auxiliary function $z\in H^2(\Omega^-\cup\Omega^+)$ to the interface problem \eqref{model} with the right-hand side $f$ replaced by $u-u_h\in L^2(\Omega)$. 
As $u-u_h$ is not continuous across faces in $\mathcal{F}^i_h$ and $\Gamma$, we apply integration by parts to obtain
\begin{equation}
\begin{split}
\label{thm_l2_error_bound_eq1}
\| u - u_h \|_{L^2(\Omega)}^2 =  a_h(z , u -u_h ) - b_h(z,u - u_h) .
 \end{split}
\end{equation}
Note that $z$ has homogeneous jumps, and thus $I_hz$ involves no non-homogeneous enrichment. 
Lemmas \ref{lem_interp_enrg_error}, \ref{lem_ErrEqn} and Theorem \ref{thm_error_bound} for $z$ show that
\begin{equation*}
\begin{split}
\label{thm_l2_error_bound_eq2}
 & a_h(z, u -u_h ) =  a_h(z - I_h z, u -u_h ) + a_h(I_h z, u -u_h )   \\
\lesssim & h^2 \| z \|_{H^2(\Omega^-\cup\Omega^+)} \| u \|_{H^2(\Omega^-\cup\Omega^+)} + h^{3/2}  \| u \|_{H^{2}(\Omega^-\cup\Omega^+)}  \| \nabla I_hz \|_{L^2(\Omega^{\Gamma}_h)} .
\end{split}
\end{equation*}
We then recover the additional $h^{1/2}$ order for the second term in \eqref{thm_l2_error_bound_eq2}. 
We note that $\| z \|_{H^2(\Omega^{\Gamma}_h)}$ only appears on the interface elements and $\Omega^{\Gamma}_h \subseteq S_{\hat{\delta}}:= \{ X\in\Omega: \text{dist}(X,\Gamma) \le \hat{\delta} \}$ with $\hat{\delta}= \mathcal{O}(h)$. Then, applying Lemma \ref{lem_delta}, we have
\begin{equation*}
\label{thm_l2_error_bound_eq9}
\| \nabla z \|_{L^2(\Omega^{\Gamma}_h)} \lesssim \sqrt{\hat{\delta}} \|  z \|_{H^2(\Omega^-\cup\Omega^+)} \lesssim h^{1/2} \|  z \|_{H^2(\Omega^-\cup\Omega^+)}.
\end{equation*}
Then, we obtain
\begin{equation*}
\label{thm_l2_error_bound_eq3}
\| \nabla I_hz \|_{L^2(\Omega^{\Gamma}_h)}  \le \| \nabla (I_hz - z) \|_{L^2(\Omega^{\Gamma}_h)} + \| \nabla z \|_{L^2(\Omega^{\Gamma}_h)} \lesssim h \| z \|_{H^2(\Omega^-\cup\Omega^+)} + h^{1/2} \| z \|_{H^2(\Omega^-\cup\Omega^+)}.
\end{equation*}
Putting \eqref{thm_l2_error_bound_eq3} into \eqref{thm_l2_error_bound_eq2} and using the regularity estimate $\| z \|_{H^2(\Omega^-\cup\Omega^+)}\lesssim \| u-u_h \| _{L^2(\Omega)}$, we have
\begin{equation}
\label{thm_l2_error_bound_eq3_1}
a_h(z, u -u_h )  \lesssim h^2 \| u - u_h \|_{L^2(\Omega)} \| u \|_{H^2(\Omega^-\cup\Omega^+)} .
\end{equation}

The main difficulty is on $b_h(z,u - u_h) $. By  H\"older's inequality, we have
\begin{equation*}
\label{thm_l2_error_bound_eq4}
b_h(z,u - u_h)  =  \int_{\Gamma}  \{\beta \nabla z\cdot \mathbf{ n}\} \jump{u - u_h} \dd s
\le \|  \{\beta \nabla z\cdot \mathbf{ n}\} \|_{L^2(\Gamma)}  \|  \jump{u - u_h}  \|^2_{L^2(\Gamma)} .
\end{equation*}
By the assumption that $q_1\in H^{5/2}(\Gamma)$ and $\partial\Omega$ and $\Gamma$ are smooth, 
we can construct auxiliary functions $\hat{u}^{\pm}|_{\Omega^{\pm}}\in H^{3}(\Omega^{\pm})$ such that $\hat{u}^+=0$ in $\Omega^+$ and $\hat{u}^-|_{\Gamma}=-q_1$. 
The existence of $\hat{u}^-$ is based on Theorem 3.37 in \cite{2000Mclean} which also gives the following stability estimate
\begin{equation}
\label{thm_l2_error_bound_eq4_2}
\| \hat{u}^- \|_{H^{3}(\Omega^-)}\lesssim \| q_1 \|_{H^{5/2}(\Gamma)}.
\end{equation}
Their Sobolev extensions are also denoted as $\hat{u}^{\pm}_E$ where $\hat{u}^-=0$ in $\Omega$ but $\hat{u}^-_E$ is non-trivial. 
Here, we do not use the original functions $u^{\pm}$ as they may not have such higher regularity. 
Then, $\jump{u - u_h}_{\Gamma} = \jump{\hat{u} - u_h}_{\Gamma} =  \jump{\hat{u} - I_h \hat{u}}_{\Gamma}$. 
Using Remark \ref{rem_IqT_est} and the same argument to Lemma \ref{thm_Ih_ext}, we can show that for $k=0,1$
\begin{equation}
\begin{split}
\label{thm_l2_error_bound_eq4_1}
|I_T \hat{u} - \hat{u} |_{H^k(\omega_T)} & \lesssim h^{2-k}_T ( \| \hat{u}^-_E \|_{H^2(\omega_T)} + \| \hat{u}^+_E \|_{H^2(\omega_T)} ) \\
& +  h^{5/2-k}_T ( \| q_1 \|_{H^{5/2}(\Gamma^{\omega_T})} + \| \jump{ \beta \nabla \hat{u}\cdot\bfn } \|_{H^{1}(\Gamma^{\omega_T})} )  + h^{3/2-k}_T ( \| \hat{u}^-_E \|_{H^1(\omega^{\text{int}}_T)} + \| \hat{u}^+_E \|_{H^1(\omega^{\text{int}}_T)} ).
\end{split}
\end{equation}
Then, applying the trace inequality into the inequality above, we obtain
\begin{equation*}
\begin{split}
\label{thm_l2_error_bound_eq5}
\| \jump{\hat{u} - I_h \hat{u}} \|_{\Gamma^T} & \lesssim \sum_{s = \pm} h^{-1/2}_T \| \hat{u}^s - I^s_h \hat{u} \|_{L^2(\omega^s_T)} +  h^{1/2}_T| \hat{u}^s - I^s_h \hat{u} |_{H^1(\omega^s_T)}\\
& \lesssim  h^{3/2}_T ( \| \hat{u}^-_E \|_{H^2(\omega_T)} + \| \hat{u}^+_E \|_{H^2(\omega_T)} ) +  h^{2}_T ( \| q_1 \|_{H^{5/2}(\Gamma^{\omega_T})} + \| \jump{ \beta \nabla \hat{u} \cdot\bfn} \|_{H^{1}(\Gamma^{\omega_T})} ) \\
& + h_T ( \| \hat{u}^-_E \|_{H^1(\omega^{\text{int}}_T)} + \| \hat{u}^+_E \|_{H^1(\omega^{\text{int}}_T)} ).
\end{split}
\end{equation*}
Recall the $\delta$ strips $S_{\delta}\supseteq \cup_{T\in\mathcal{T}^i_h}\omega^{\text{int}}_T$ and $S_{\hat{\delta}} \supseteq \cup_{T\in\mathcal{T}^i_h}\omega_T$ with $\delta \simeq h^2$ and $\hat{\delta} = h$.
Then, by the finite overlapping property of $\omega_T$ and Lemma \ref{lem_delta}, we obtain
\begin{equation}
\begin{split}
\label{thm_l2_error_bound_eq6}
\| \jump{\hat{u} - I_h \hat{u}} \|_{\Gamma} & \lesssim h^{3/2} ( \| \hat{u}^-_E \|_{H^2(S_{\hat{\delta}})} + \| \hat{u}^+_E \|_{H^2(S_{\hat{\delta}})} ) + h ( \| \hat{u}^-_E \|_{H^1(S_{\delta})} + \| \hat{u}^+_E \|_{H^1(S_{\delta})} )\\
&+ h^2 ( \| q_1 \|_{H^{5/2}(\Gamma)} + \| \jump{ \beta \nabla \hat{u} \cdot\bfn } \|_{H^{1}(\Gamma)} )  + h \| \tilde{q}_1 \|_{H^1(S^2)} \\
& \lesssim h^2 \left( \| \hat{u}^-_E \|_{H^3(S_{\hat{\delta}})} + \| \hat{u}^+_E \|_{H^3(S_{\hat{\delta}})} + \| q_1 \|_{H^{5/2}(\Gamma)} \right) \lesssim h^2 \| q_1 \|_{H^{5/2}(\Gamma)},
\end{split}
\end{equation}
where we have also used $\| \jump{ \beta \nabla \hat{u} \cdot\bfn } \|_{H^{1}(\Gamma)} \lesssim \| \hat{u}^- \|_{H^2(\Omega^-)}\lesssim \| \hat{u}^- \|_{H^{3}(\Omega^-)} \lesssim \| q_1 \|_{H^{5/2}(\Gamma)} $ by the boundedness of Sobolev extensions and \eqref{thm_l2_error_bound_eq4_2}. 
Putting \eqref{thm_l2_error_bound_eq6} and \eqref{thm_l2_error_bound_eq3_1} into \eqref{thm_l2_error_bound_eq1} leads to the desired estimate.
\end{proof}

\subsection{Conditioning analysis}

Now, we proceed to analyze the conditioning of the resulting linear system
\begin{equation}
\label{linearsys}
\bfA \bar{\bfu} = \bar{\bff}.
\end{equation}
Here, we highlight that the size and algebraic structure of $\bfA$ are identical to the matrix of the standard FE spaces on the same mesh regardless of interface. 
This very feature does not only benefit solving moving interface but also is the key for the robustness of the linear system.
We begin with giving the rigorous definition of the isomorphism \eqref{isomap0}:
\begin{equation}
\label{isomap1}
\mathbb{I}_h v_h(X) = v_h(X), ~~~ \forall X\in\mathcal{N}_h.
\end{equation}
It preserves the algebraic structure of the linear system, ensuring that the resulting condition number exhibits the standard $\mathcal{O}(h^2)$ growth, independent of the interface location. In this subsection, we provide a theoretical proof of this statement.
One key is that the isomorphism $\mathbb{I}_h$ is the following equivalence.
\begin{lemma}
\label{lem_iso}
The following equivalence holds regardless of interface location.
\begin{equation}
\label{lem_iso_eq0}
\| \mathbb{I}_h v_h \|_{L^2(\Omega)} \simeq \| v_h \|_{L^2(\Omega)}, ~~~ \forall v_h \in \mathcal{S}^0_h. 
\end{equation}
\end{lemma} 
\begin{proof}
The argument is similar to Lemma 4.2 in \cite{2020Guo}.
\end{proof}

For each $v_h\in \mathcal{S}^0_{h}$, we define $\mathfrak{I}v_h$ as the vector of the coefficients of the global shape functions associated with each node. 
Trivially, there holds that $\mathfrak{I}v_h = \mathfrak{I}\mathbb{I}_hv_h$. Denote $\|\cdot\|_2$ as the $l^2$ norm for vectors. We have the following result.
\begin{lemma}
\label{lemma_l2equiv}
The following equivalence holds regardless of interface location:
\begin{equation}
\label{FE_condNum_est1}
 \| v_h \|_{L^2(\Omega)} \simeq h \| \mathfrak{I}v_h \|_2, ~~~ \forall v_h\in \mathcal{S}^0_{h}.
\end{equation}
\end{lemma}
\begin{proof}
By the results of standard Lagrange elements, we clearly have
\begin{equation*}
\label{FE_condNum_est2}
\| \mathbb{I}_h v_h \|_{L^2(\Omega)} \simeq h \| \mathfrak{I} \mathbb{I}_h v_h \|_2, ~~~ \forall v_h\in \mathcal{S}^0_{h}.
\end{equation*}
Then, the desired result follows from Lemma \ref{lem_iso}.
\end{proof}

Now, we are ready to state and prove an equivalence between $\|\cdot\|_{L^2(\Omega)}$ and $\vertiii{\cdot}_h$ norms on $\mathcal{S}^0_h$.

\begin{lemma}
\label{lem_glob_energy_norm_equiv}
For every  $v_h\in \mathcal{S}^0_h$, there holds
\begin{equation}
\label{rem_L2_energy_equiv_eq}
\| v_h \|_{L^2(\Omega)}\lesssim  \vertiii{v_h} \lesssim h^{-1} \| v_h \|_{L^2(\Omega)}.
\end{equation}
\end{lemma}
\begin{proof}
Let us first prove the left inequality in \eqref{rem_L2_energy_equiv_eq}. 
The argument follows from Lemma 2.1 in \cite{1982Arnold}. Given a $v_h\in \mathcal{S}^0_h$,
we define an auxiliary function $z=(z^+,z^-) \in H^2(\Omega^-\cup\Omega^+)$ as the solution of the
interface problem \eqref{model} satisfying homogeneous jump conditions with $f=v_h$. 
Still let $z^{\pm}_E \in H^2(\Omega)$ be the Sobolev extensions of $z$ from $\Omega^{\pm}$ to $\Omega$. 
Testing $\nabla\cdot(\beta\nabla z) = v_h$ by $v_h$, we write 
\begin{equation}
\begin{split}
\label{lem_poc_eq3}
\| v_h \|^2_{L^2(\Omega)} &= \int_{\Omega} \beta \nabla z \cdot \nabla v_h \dd X - \sum_{F\in\mathcal{F}^i_h} \int_F \{ \beta \nabla z\cdot \mathbf{ n} \} \jump{v_h} \dd s - \sum_{T\in\mathcal{T}^i_h} \int_{\Gamma^T} \{ \beta \nabla z\cdot \mathbf{ n} \}_{\Gamma} \jump{v_h}_{\Gamma} \dd s .
\end{split}
\end{equation}
Applying \eqref{thm_trace_interface_eq0} and trace inequality, we obtain
\begin{equation*}
\label{lem_poc_eq3_1}
\int_{\Gamma^T} \{ \beta \nabla z\cdot \mathbf{ n} \}_{\Gamma} \jump{v_h}_{\Gamma} ds \le \| \{ \beta \nabla z\cdot \mathbf{ n} \} \|_{L^2(\Gamma^T)} \| \jump{v_h} \|_{L^2(\Gamma^T)} \lesssim h_T \sum_{s=\pm}( |z^s_E|_{H^1(T)} + |z^s_E|_{H^2(T)} ) \| \nabla v_h \|_{L^2(T)} .
\end{equation*}
Then, by H\"older's inequality we have
\begin{equation*}
\begin{split}
\label{lem_poc_eq3_2}
\| v_h \|^2_{L^2(\Omega)} &\leqslant \left( \| \ \nabla v_h \|^2_{L^2(\Omega)} +  \sum_{F\in\mathcal{F}^i_h} h^{-1}_T \| [v_h] \|^2_{L^2(F)}  \right)^{1/2}  \cdot  \left( \| \beta \nabla z \|^2_{L^2(\Omega)} + \sum_{F\in\mathcal{F}^i_h}  h_T \| \{ \beta \nabla z\cdot \mathbf{ n} \} \|^2_{L^2(F)} \right)^{1/2}.
\end{split}
\end{equation*}
In addition, for each $F\in\mathcal{F}^i_h$, let $T$ be one element that has $e$ as its edge. 
Then, the trace inequality on $F$ and $T$ yields
\begin{align}
\| h^{1/2}_T \{ \beta \nabla z\cdot \mathbf{ n} \} \|_{L^2(F)}   \lesssim \sum_{s=\pm} \left( | z^{s} |_{H^1(T)} + h_T | z^{s} |_{H^2(T)} \right). \label{lem_poc_eq5}
\end{align}
We sum \eqref{lem_poc_eq5} over all the faces in $\mathcal{F}^i_h$ to obtain an upper bound for the jump terms. 
Then, combining the resulting bounds and applying the elliptic regularity, we obtain
\begin{equation}
\begin{split}
\label{lem_poc_eq7}
& \left( \| \beta \nabla z \|^2_{L^2(\Omega)} +  \sum_{F\in\mathcal{F}^i_h} h_T \| \{ \beta \nabla z\cdot \mathbf{ n} \} \|^2_{L^2(F)} \right)^{1/2} 
 \lesssim   \sum_{s=\pm} \beta^s  (| z^s |_{H^1(\Omega^s)} + | z^s |_{H^2(\Omega^s)} ) \lesssim \| v_h \|_{L^2(\Omega)}.
\end{split}
\end{equation}
Finally, substituting \eqref{lem_poc_eq7} into \eqref{lem_poc_eq3} leads to the desired result. 
The right one follows from the trace and inverse inequality given by Theorem \ref{thm_trace_inequa},
\end{proof}

\begin{lemma}
\label{lem_condNum}
Let $\kappa_2(\cdot)$ be the spectral condition number of a matrix. Then $\kappa_2(\bfA)\lesssim h^{-2}$. 
\end{lemma}
\begin{proof}
The estimate immediately follows from Theorem 3.1 in \cite{2006ErnGuermond} together with Lemmas \ref{lem_glob_energy_norm_equiv} and \ref{lemma_l2equiv}.
\end{proof}


\section{Numerical Experiments}
\label{sec:num_sec}

In this section, we present numerical experiments to demonstrate the performance of the proposed IFE method. All tests are conducted on structured tetrahedral meshes. The solution domain is first partitioned into $N\times N\times N$ cuboids, and each cuboid is further subdivided into six tetrahedra, resulting in a total of $6N^3$ elements. The enriched IFE method \eqref{IFE_scheme} is then employed to solve the interface problems. 
Due to the isomorphism of DoFs between the enriched IFE and the standard FE methods, a algebraic multigrid (AMG) preconditioner is immediately applicable.
The AMG solver is used in all the computation presented below, 
and a comprehensive test regarding its stability with respect to interface location is given in Example 4 below.
All numerical experiments are carried out in MATLAB 2021b on iMac 8-core 4GHz i9 and 64Gb RAM.

\subsubsection*{Example 1 (Spherical Interface)}
In this example, we consider an elliptic equation with nonhomogeneous interface jump conditions across a spherical interface. Let $\Omega = (-1,1)^3$ and the interface $\Gamma =\{(x,y,z): \gamma(x,y,z)=0\}$ be a sphere given by the level set function 
\[\gamma(x,y,z) =x^2+y^2+z^2 -r_0^2,\]
where $r_0=\pi/4$.
The exact solution of the elliptic PDE is given by
\begin{equation}
u(x,y,z) = 
\left\{
\begin{split}
&\sin(x^2+y^2+z^2)~~~& \text{in}~\Omega^- := \{(x,y,z)\in\Omega:\gamma(x,y,z)<0\},\\
&\cos(x^2+y^2+z^2)~~~& \text{in}~\Omega^+ := \{(x,y,z)\in\Omega:\gamma(x,y,z)>0\}.
\end{split}
\right.
\end{equation}
The coefficients are chosen to be  and $\beta^- = 1$, and $\beta^+ = 100$. The jump data $q_1$ and $q_2$ are computed accordingly. 

\begin{figure}[!]
\begin{center}
\includegraphics[width=.49\textwidth]{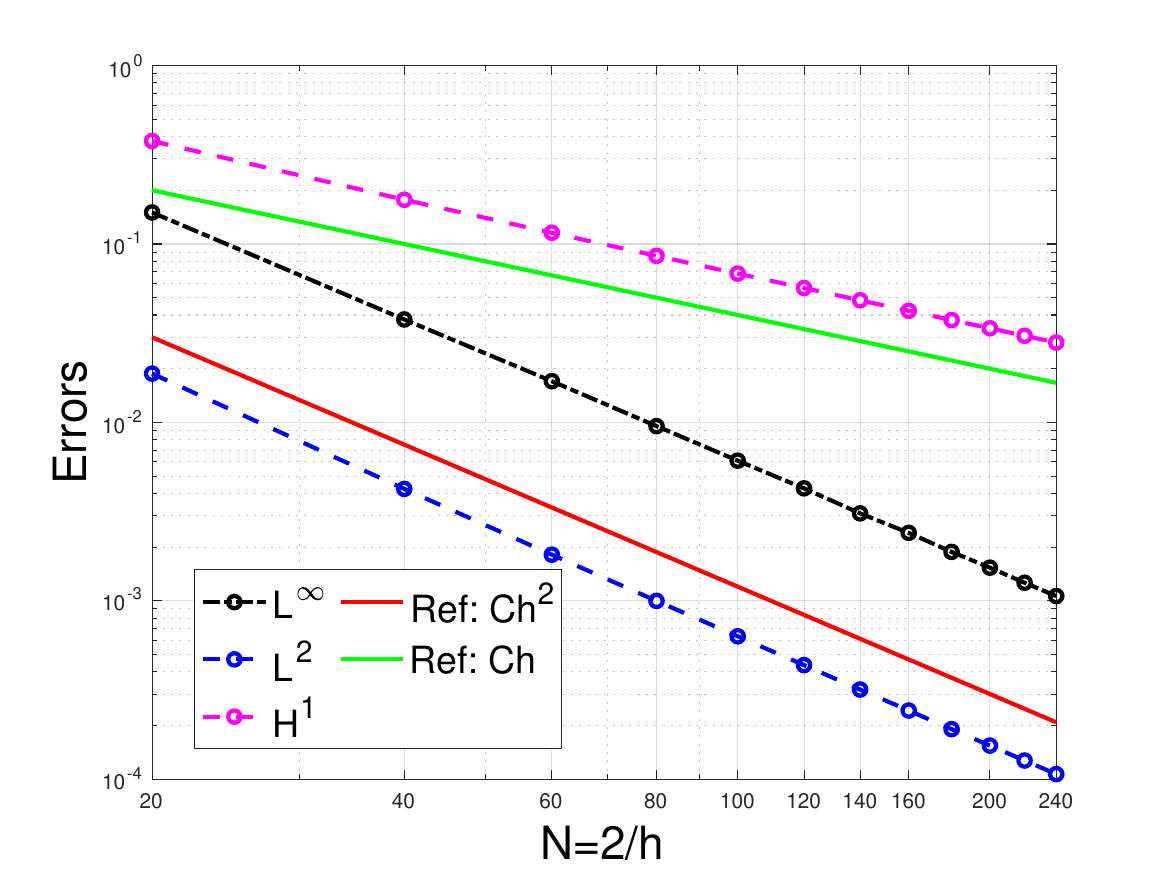}
\includegraphics[width=.49\textwidth]{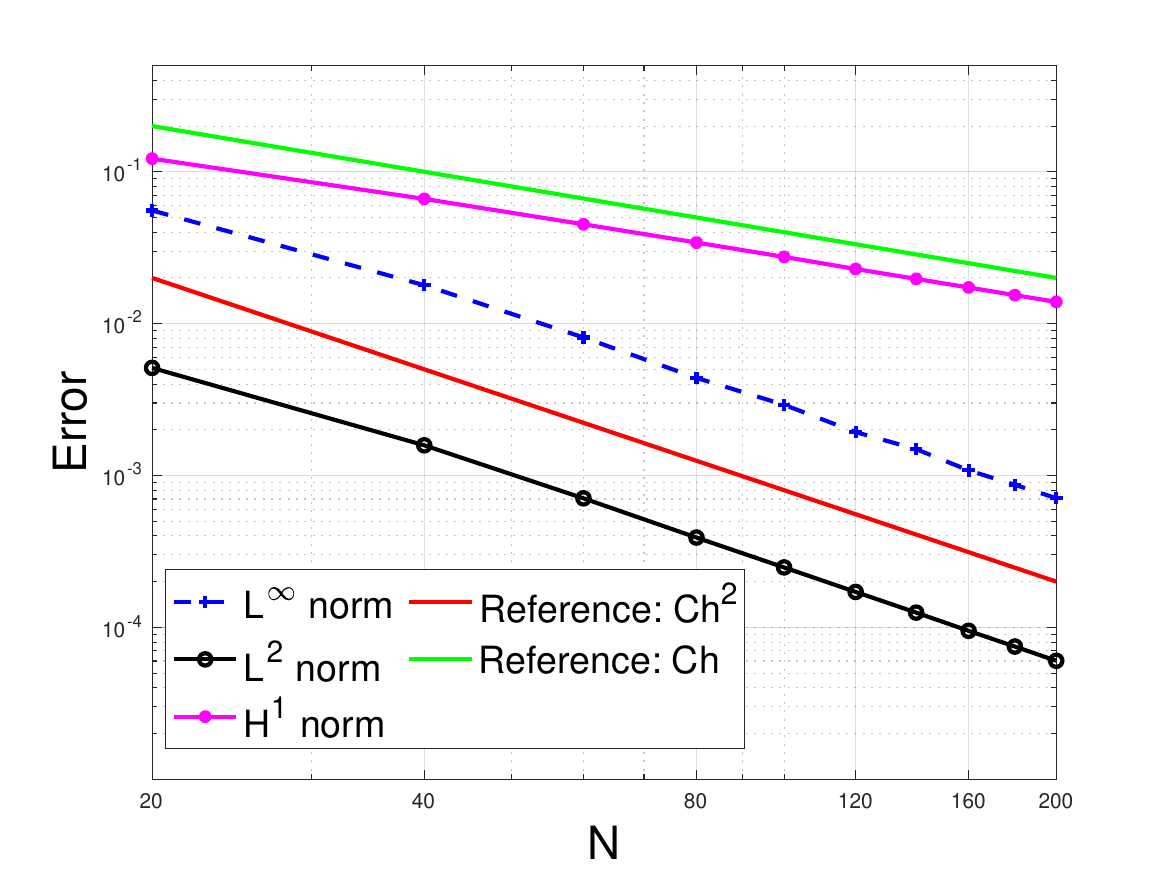}
\end{center}
\caption{Convergence of IFE solutions: Example 1 spherical interface (left) and Example 2 orthocircle interface (right).}
\label{fig: cvg ex1}
\end{figure}

We record IFE solutions on a sequence of structured tetrahedral meshes with $N=20,40,60,\cdots, 240$. The finest mesh contains approximately $83$ million elements, of which $765$ thousands are interface elements, giving the ratio of $|\mathcal{N}_h^i|/|\mathcal{N}_h| \approx 0.92\%$. The total number of DoFs on the finest mesh is about 14 million. Figure \ref{fig: cvg ex1} shows the errors in the $L^\infty$, $L^2$, and $H^1$ norms, corroborating the theoretical estimates \eqref{thm_error_bound_eq0} and \eqref{thm_l2_error_bound_eq0}. 

To accelerate the 3D matrix assembly, we employed the vectorization technique described in \cite{2008Chen}, which dramatically reduces time required to construct global matrices and vectors. Table \ref{tab:examp1} lists the CPU time (in seconds) for Example 1 with the linear enriched IFE method 
The time of assembling the matrix system and solving the linear equations are comparable to those reported for fitted-mesh finite element methods such as \cite{2017ChenWeiWen}, while the mesh-generation cost is negligible. This computational efficiency is particularly advantageous for time-dependent problems with moving interfaces.  
\vspace{2mm}
\begin{table}
\centering
\begin{tabular}{r r r r r r r r c}
\hline
N     & \#Cell    & DoF       & {Mesh} & IFEM  & Matrix  & {Solve} \\\hline 
40    & 384,000   & 68,921    & {0.01} & 4.31 &  2.47   & { 0.12} \\ 
80    & 3,072,000 & 531,441   & {0.05} & 16.20& 15.06   & { 1.79} \\ 
120   & 10,368,000& 1,771,561 & {0.20} & 36.73& 47.20   & { 8.51} \\ 
160   & 24,576,000& 4,173,281 & {0.63} & 64.17&123.19   &{ 27.65} \\ \hline
\end{tabular}
\caption{Computation time of each module in the IFEM for Example 1.}
\label{tab:examp1}
\end{table}


\subsubsection*{Example 2 (Complex-Shape Interface and Large Jump)}
\begin{figure}
\begin{center}
\includegraphics[width=.26\textwidth]{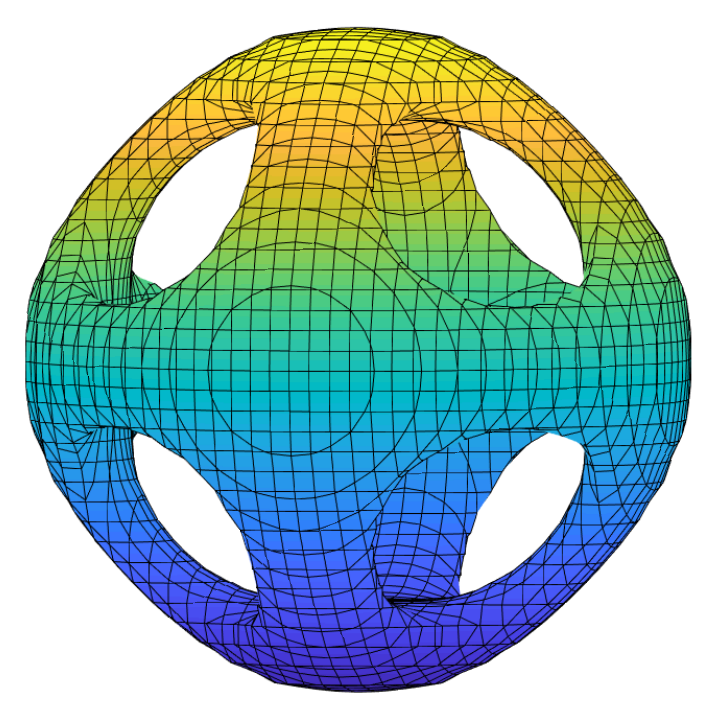}~~~~~
\includegraphics[width=.33\textwidth]{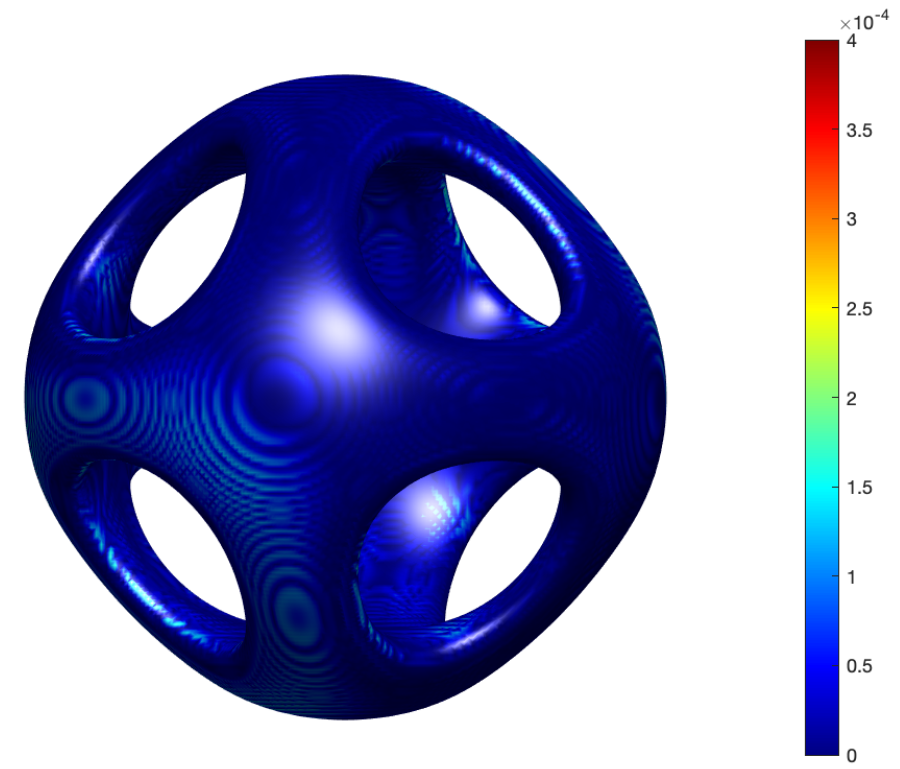}~~~~~
\includegraphics[width=.33\textwidth]{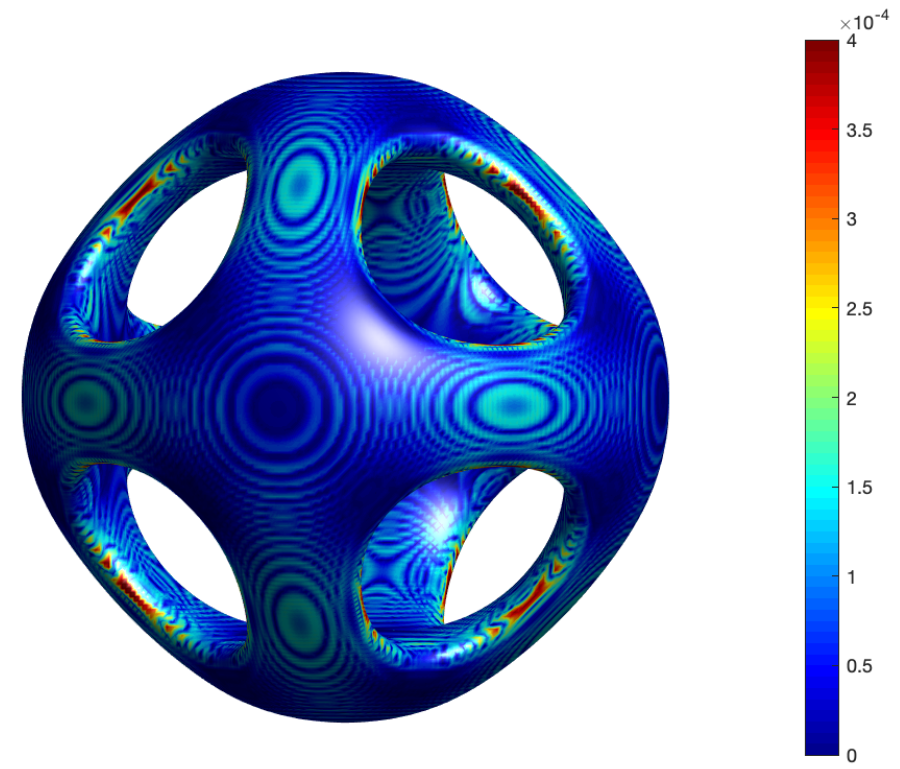}
\caption{From left: the geometry of orthocircle interface; error of linear IFE solution on interface; error of trilinear IFE solution on interface.}
\label{fig: comp Ex2}
\end{center}
\end{figure}

We next consider an interface problem whose geometry is considerably more intricate than Example 1. The interface $\gamma$ is an orthocircle defined by
\[\gamma(x,y,z) = [(x^2+y^2-1)^2+z^2][(x^2+z^2-1)^2+y^2][(y^2+z^2-1)^2+x^2]-0.075^2[1+3(x^2+y^2+z^2)].\]
The computational domain is $\Omega = (-1.2,1.2)^3$, and the interface is given by 
$\Gamma = \{(x,y,z)\in\Omega: \gamma(x,y,z) = 0\}$. 
The diffusion coefficients are piecewise constant, with $\beta^- =1$, $\beta^+=100$. The orthocircle surface is illustrated in the left plot of Figure \ref{fig: comp Ex2}. A similar interface geometry has been studied in  \cite{2017ChenWeiWen, 2021GuoZhang}. 
The exact solution is given by 
\begin{equation}
u(x,y,z) = 
\left\{
\begin{split}
&\sin(x+2y+3z)~~~& \text{in}~\Omega^- := \{(x,y,z)\in\Omega:\gamma(x,y,z)<0\},\\
&x^2+y^3-z~~~& \text{in}~\Omega^+ := \{(x,y,z)\in\Omega:\gamma(x,y,z)>0\}.
\end{split}
\right.
\end{equation}
The jump data $q_1$ and $q_2$ are computed accordingly. 

We solve the interface problem on a sequence of structured tetrahedral meshes with $N=20,40,60,\cdots, 200$. The finest mesh contains approximately $48$ million elements, of which $936$ are interface elements, so that the ratio of $|\mathcal{N}_h^i|/|\mathcal{N}_h| \approx 1.95\%$. This ratio is roughly twice of the  observed
in the spherical‑interface case in Example 1. The corresponding number of DoFs is about $8$ million. 
The right plot of Figure \ref{fig: cvg ex1} shows the errors in the $L^\infty$, $L^2$, and $H^1$ norms, once again validating the theoretical estimates \eqref{thm_error_bound_eq0} and \eqref{thm_l2_error_bound_eq0}. 

We next compare the accuracy of linear IFE basis functions on tetrahedral meshes with that of trilinear IFE basis functions on structured cuboid meshes \cite{2021GuoZhang}. For a background grid of $N^{3}$ cuboids with $N=160$, the two discretizations possess the same total number of degrees of freedom, whereas the tetrahedral mesh contains six times as many elements as its cuboid counterpart. 
On this grid roughly \(2.97\%\) of the cuboids are interface elements, compared with \(1.95\%\) of the tetrahedra.  
Numerical results show that the linear tetrahedral IFE approximation yields slightly lower errors on the interface surface, due to the finer geometric resolution of the interface inside each tetrahedral element. See the right two plots in Figure \ref{fig: comp Ex2}.

\begin{figure}[h!]
\centering
\includegraphics[width=.5\textwidth]{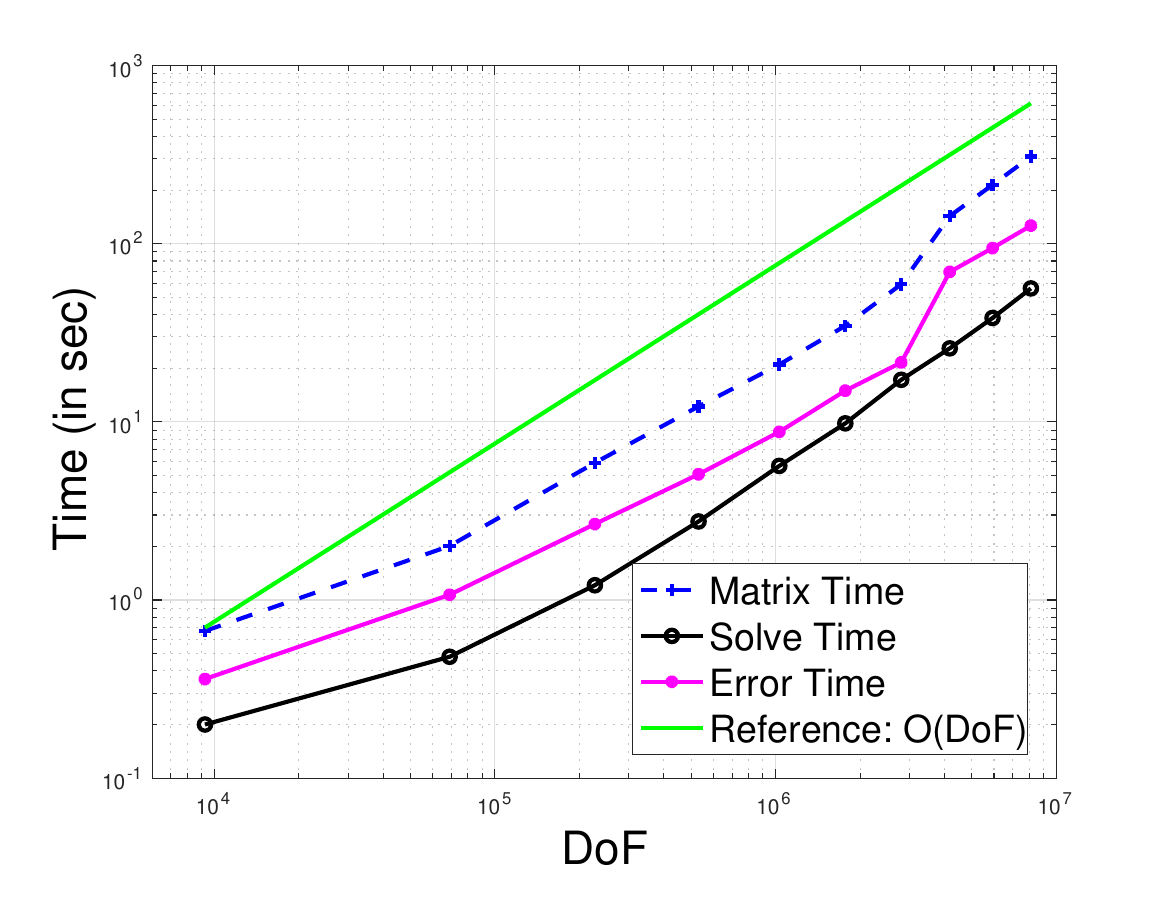}
\caption{CPU Time of IFE solutions of Example 2.}
\label{fig: time ex2}
\end{figure}

Figure \ref{fig: time ex2} records the CPU times for assembling the global system, solving the linear equations, and computing the errors; the recorded times exhibit a linear growth with respect to the number of DoFs, a trend identical to that observed for Example 1. 

\begin{figure}
\begin{center}
\includegraphics[width=.48\textwidth]{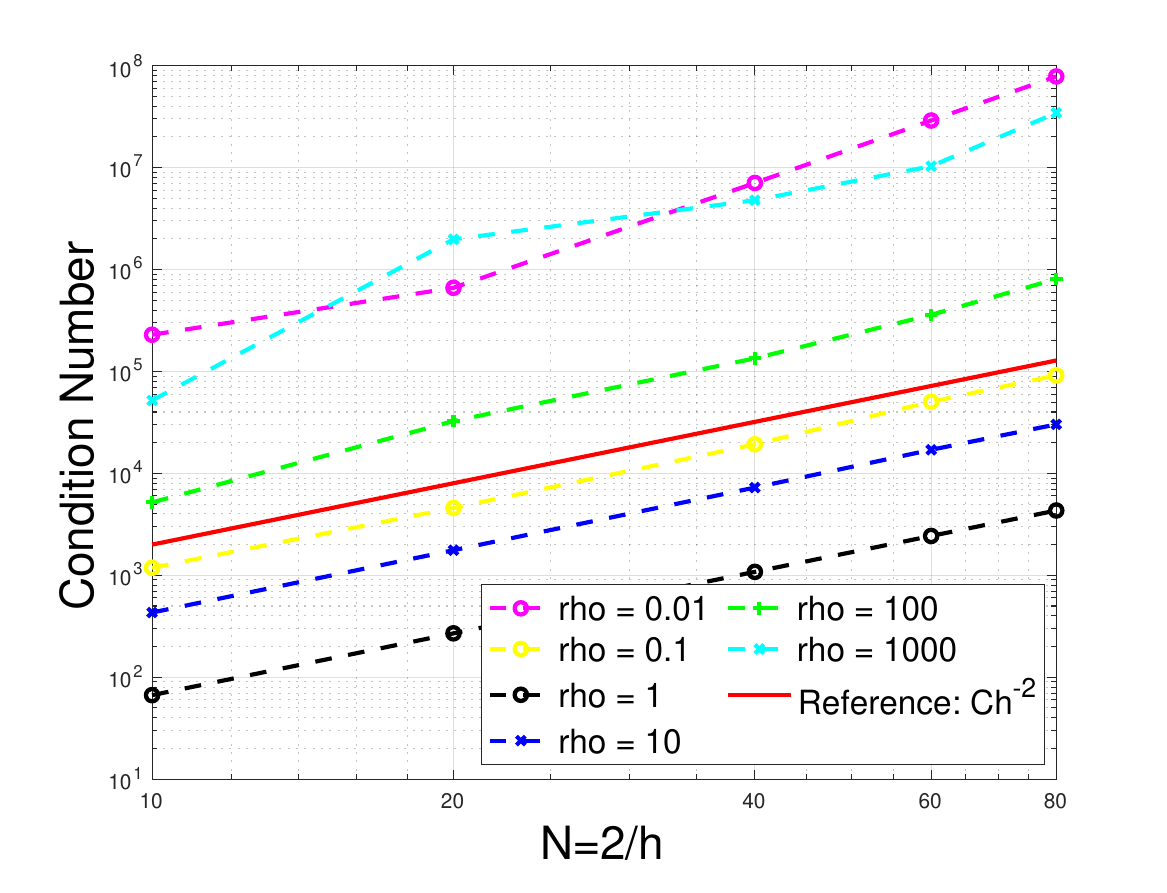}\hspace{-4mm}
\includegraphics[width=.48\textwidth]{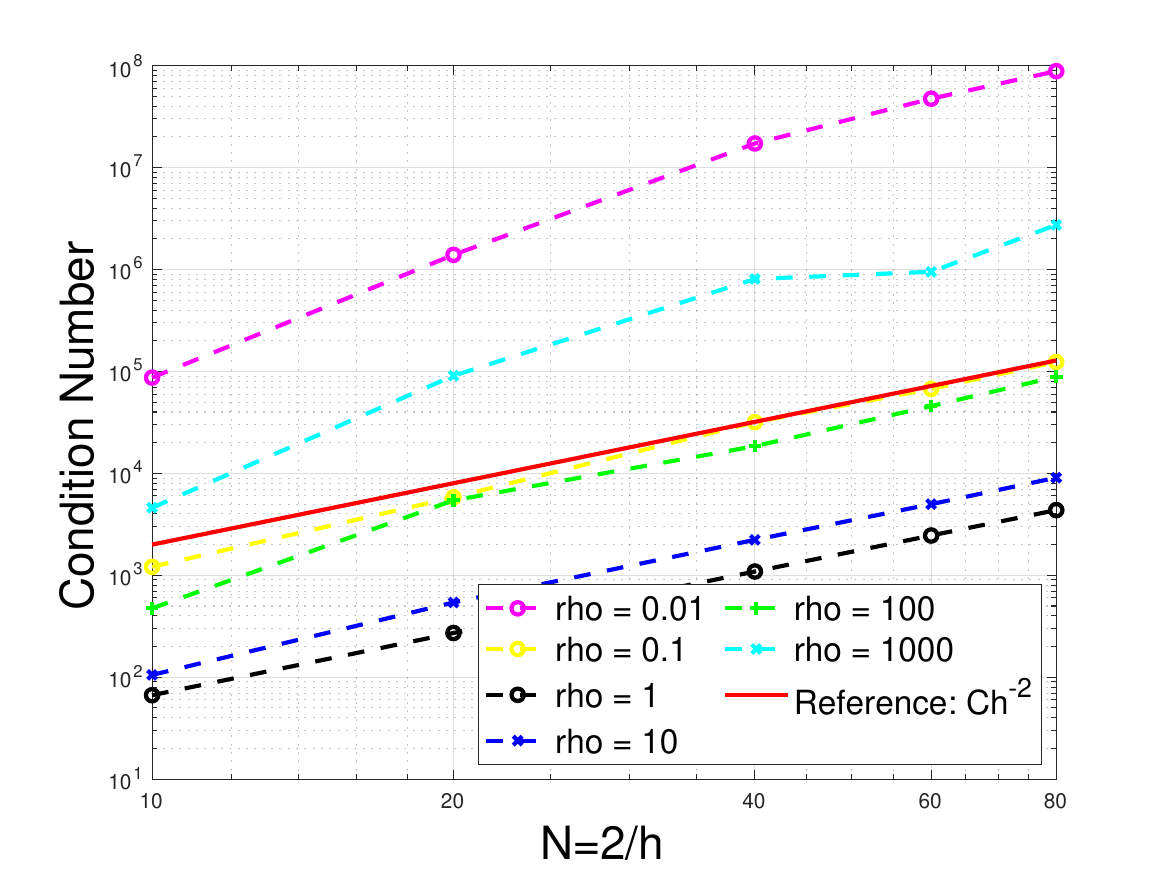}
\end{center}
\caption{Condition numbers for the spherical interface (left) and orthocircle interface (right).}
\label{fig: cond}
\end{figure}

In addition, we investigated how the coefficient jump affects the conditioning of the stiffness matrix.  Fixing $\beta^- = 1$, we vary the jump ratio
\[
\rho = \frac{\beta^{+}}{\beta^{-}}\in\{0.01,\,0.1,\,1,\,10,\,100,\,1000\}.
\]
Figure \ref{fig: cond} plots the resulting condition numbers for both interface geometries. The condition number grows monotonically with $\rho$. The larger jumps in the diffusion coefficient lead to increasingly ill-conditioned systems. 

\subsubsection*{Example 3: Robustness regarding Interface Locations}
In this example, we test the robustness of our enriched IFE method in terms of interface locations. We consider the same example designed in \cite{2024Guo}. Let $\Omega =[-1,1]^3$, and let the interface $\gamma$ be a  ``squirecle":  \begin{equation}
\gamma: x^4+y^4+z^4-r_0^4 = 0,~~~~~\text{with}~~~ r_0 = 0.75-\epsilon,
\end{equation}
whose shape is close to a square but has rounded corners. Here $\epsilon>0$ is a parameter to control the interface location relative to the mesh. Since $x = 0.75$, $y = 0.75$, and $z=0.75$ are aligned with the mesh, the smaller value of $\epsilon$ means the interface is closer to the boundary of the interface element. The exact solution is given by 
\begin{equation}
u(x,y,z) = 
\left\{
\begin{split}
&\dfrac{1}{\beta^-}(x^4+y^4+z^4)^\alpha~~~& \text{in}~\Omega^- := \{(x,y,z)\in\Omega:\gamma(x,y,z)<0\},\\
&\dfrac{1}{\beta^+}(x^4+y^4+z^4)^\alpha + (\dfrac{1}{\beta^-}-\dfrac{1}{\beta^+})r_0^{4\alpha}~~~& \text{in}~\Omega^+ := \{(x,y,z)\in\Omega:\gamma(x,y,z)>0\}.
\end{split}
\right.
\end{equation}
where $\alpha = 1/2$, $\beta^- = 1$, and $\beta^+=100$.

\begin{figure}
\begin{center}
\includegraphics[width=.48\textwidth]{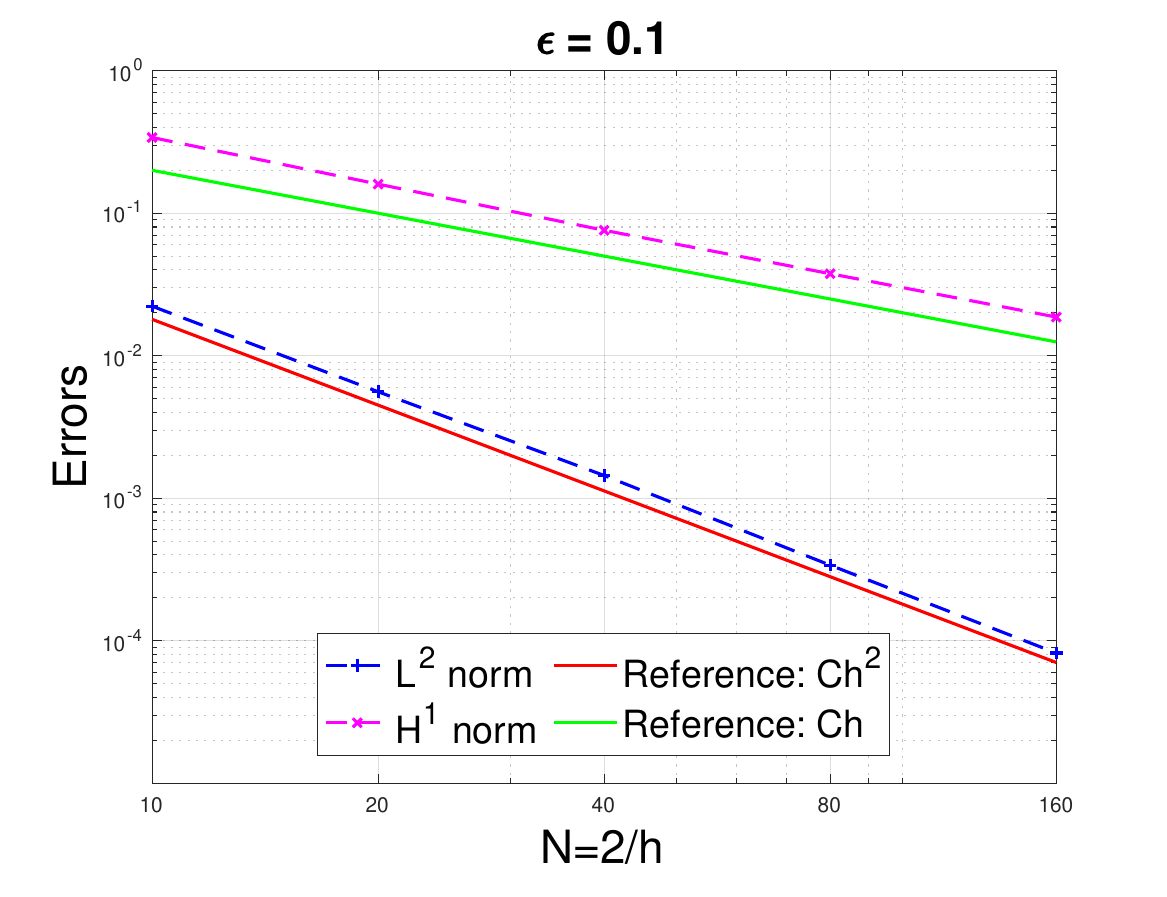}\hspace{-4mm}
\includegraphics[width=.48\textwidth]{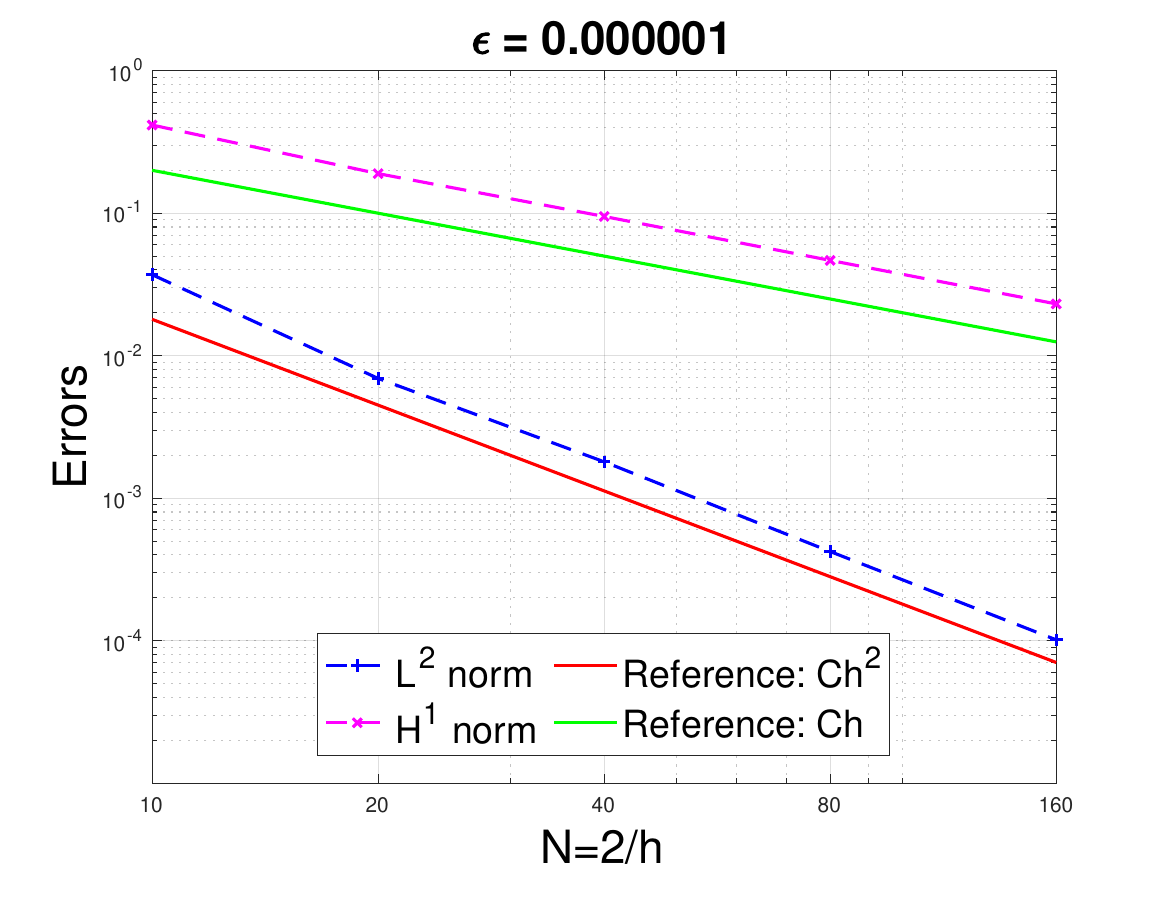}
\end{center}
\caption{The convergence for IFE solutions with $\epsilon = 10^{-1}$ (left) and $\epsilon = 10^{-6}$ (right) for Example 3.}
\label{fig: epsilon impact}
\end{figure}

Figure \ref{fig: epsilon impact} shows the errors of IFE solutions for $\epsilon = 10^{-1}$ and $\epsilon = 10^{-6}$. In both cases, the $L^2$- and $H^1$-norm errors decay at the optimal rate. Table \ref{tab:epsilon impact} list the AMG iteration counts required to reach a residual tolerance of $10^{-8}$. These results demonstrate that the IFE method and AMG solver remain robust, even in the case of small-cut interface elements.

\begin{table}[h!]
  \centering
  \caption{Iterations of AMG with two types of interface locations}
  \label{tab:epsilon impact}
  \begin{tabular}{|c|c|c|c|c|c|c|c|c|c|c|}
    \hline
    N  & \multicolumn{2}{c}{10}& \multicolumn{2}{|c}{20}   & \multicolumn{2}{|c}{40} & \multicolumn{2}{|c}{80} & \multicolumn{2}{|c|}{160}\\ \hline
    $\epsilon$ &$10^{-1}$ & $10^{-6}$  &$10^{-1}$ & $10^{-6}$&$10^{-1}$ & $10^{-6}$&$10^{-1}$ & $10^{-6}$ &$10^{-1}$ & $10^{-6}$\\
       \hline
     AMG iterations & 11 & 9 & 11&  11 & 11 & 11 &12 & 12 &13 &13\\
     \hline
  \end{tabular}
\end{table}

\subsubsection*{Example 4: Interface without analytical equation}
In this example, we apply the enriched IFE method to real-world interfaces for which no analytic solutions exist. These interfaces are available as cloud points and we refer to \cite{2011RouhaniSappa} for the availability of data of cloud points. 
For each set of cloud points, a smooth interface surface $\gamma$ is reconstructed using the {implicit B-spline} technique introduced in \cite{2015RouhaniSappaBoyer}. A signed-distance function is then formed by computing the distance from each mesh nodes to cloud data points. The zero level-set of the signed-distance function defines the computational interface. Each interface is embedded in a bounded cuboid domain $\Omega$. The exact solution is set to be
\begin{equation}
u(x,y,z) = 
\left\{
\begin{split}
&\sin(x^2+y^2+z^2)~~~& \text{in}~\Omega^- := \{(x,y,z)\in\Omega:\gamma(x,y,z)<0\},\\
&\cos(x^2+y^2+z^2)~~~& \text{in}~\Omega^+ := \{(x,y,z)\in\Omega:\gamma(x,y,z)>0\}.
\end{split}
\right.
\end{equation}
with the diffusion coefficients $\beta^- = 1$ and $\beta^+=10$. The jump date $q_1$ and $q_2$ follow directly. 
 
\begin{figure}
\includegraphics[height=1.9in]{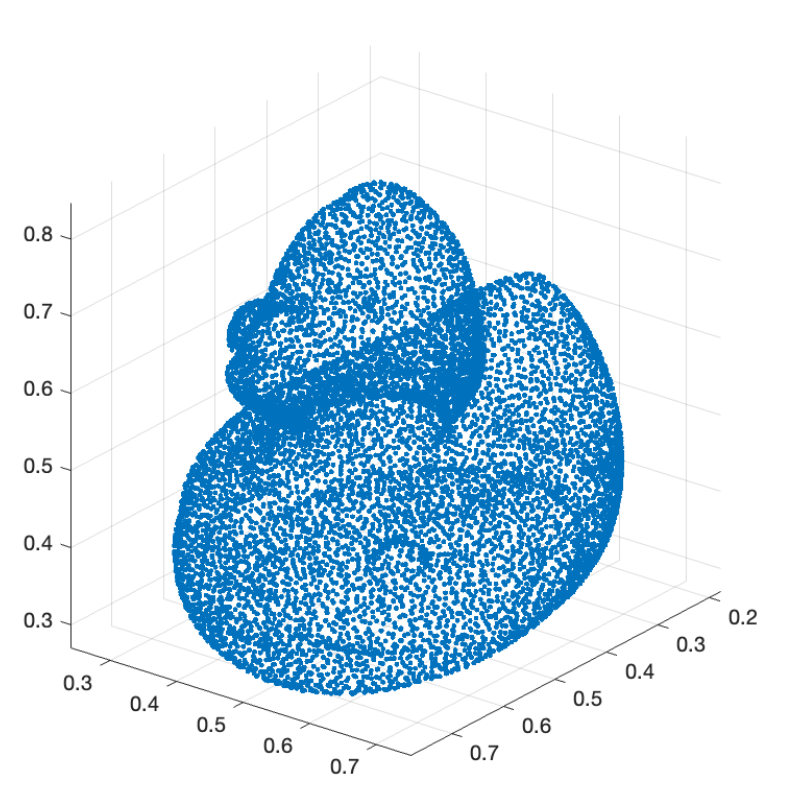}
\includegraphics[height=1.9in]{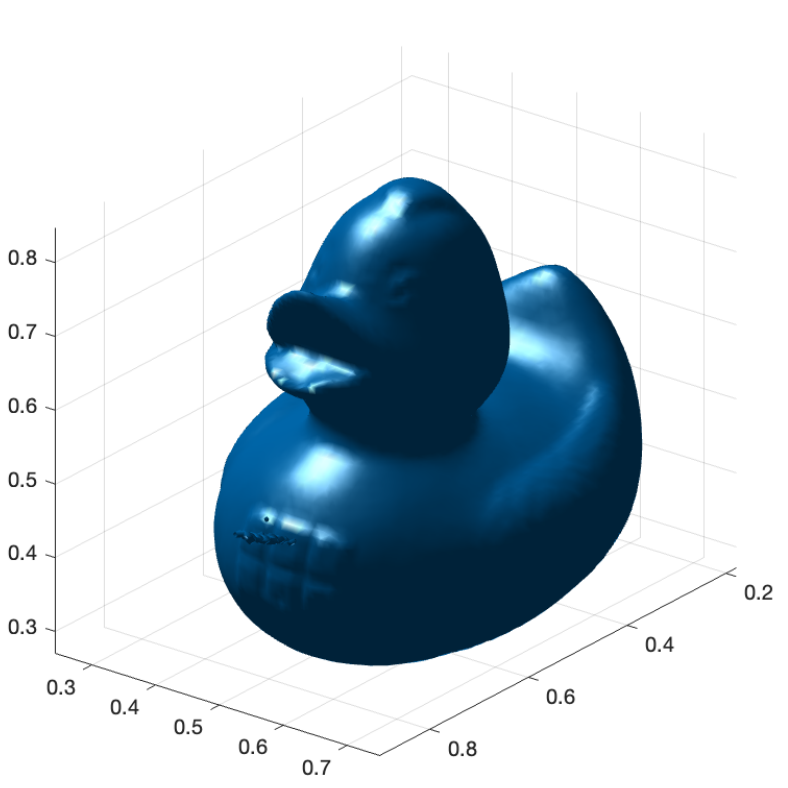}
\includegraphics[height=1.9in]{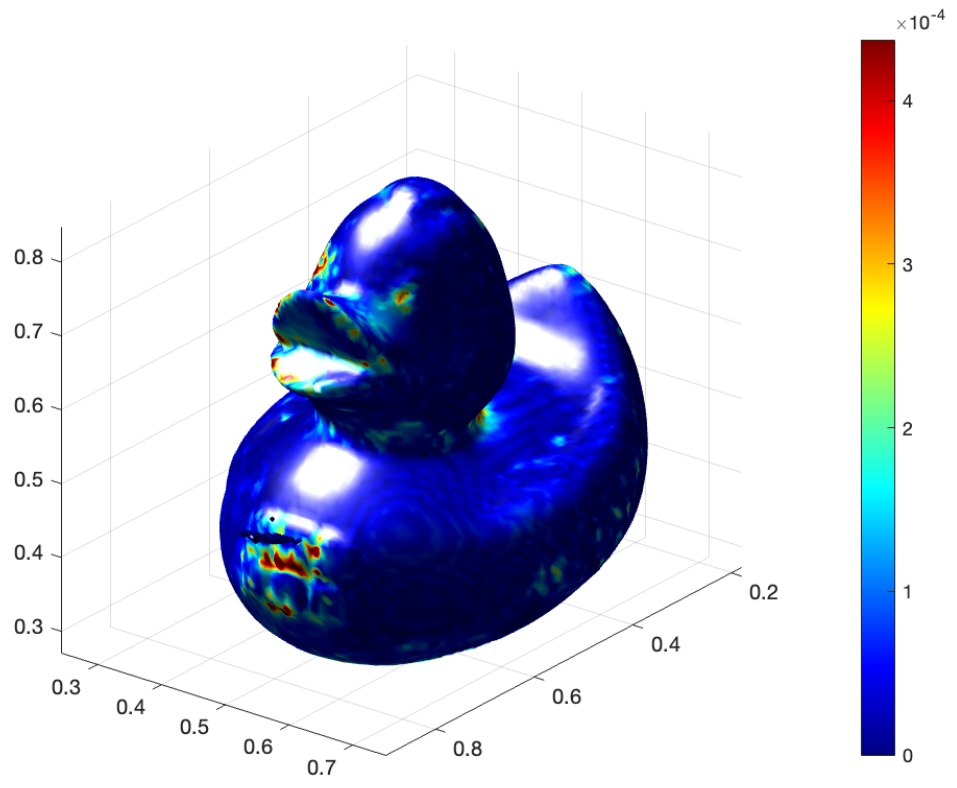}
\caption{From left: cloud points of duck interface, reconstructed computational interface, error of IFE solution on the interface.}
\label{fig: duck}
\end{figure}

\begin{figure}
\begin{center}
\includegraphics[height=2in]{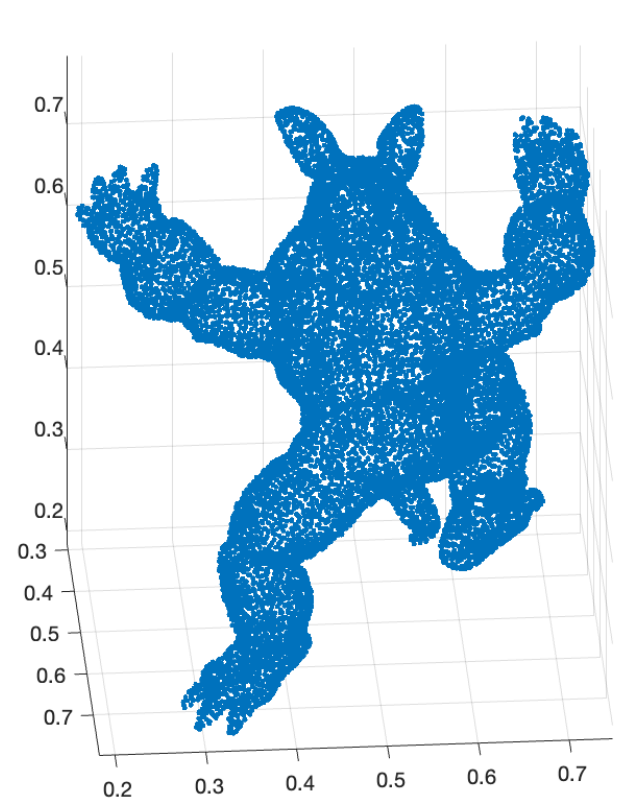}~~~~
\includegraphics[height=2in]{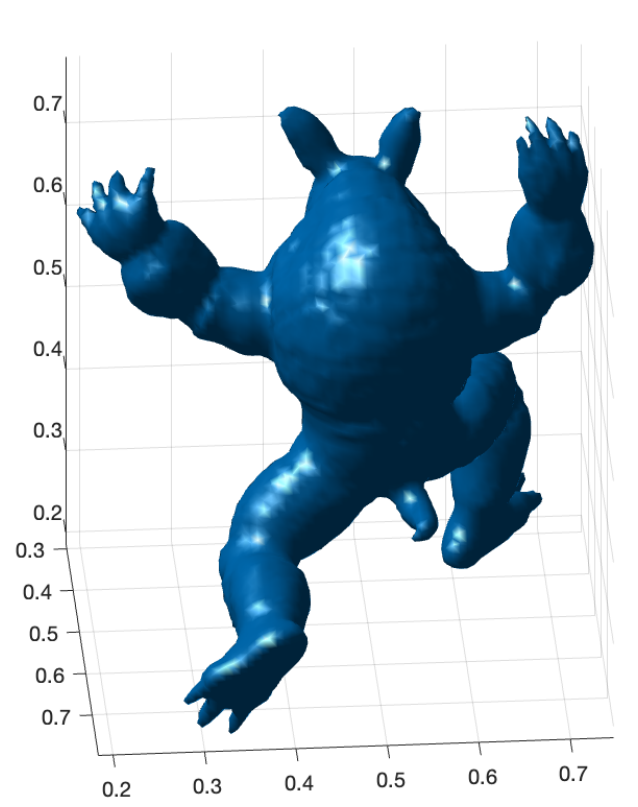}~~~~
\includegraphics[height=2in]{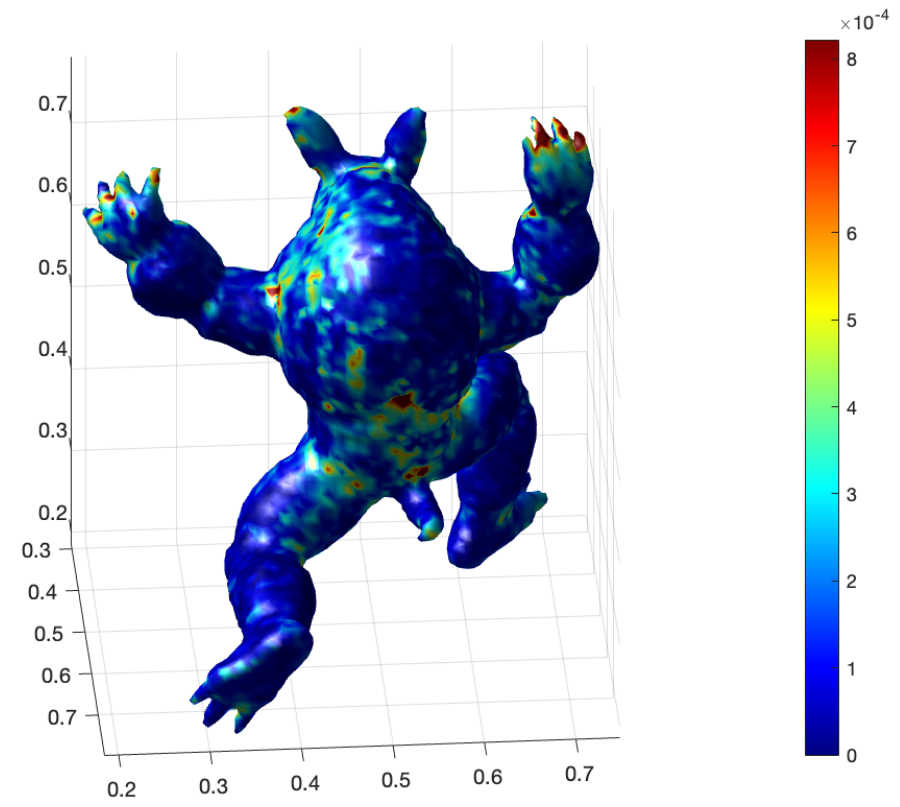}
\end{center}
\caption{From left: cloud points of armadillo interface, reconstructed computational interface, error of IFE solution on the interface.}
\label{fig: armadillo}
\end{figure}

\begin{figure}
\begin{center}
\includegraphics[height=2.2in]{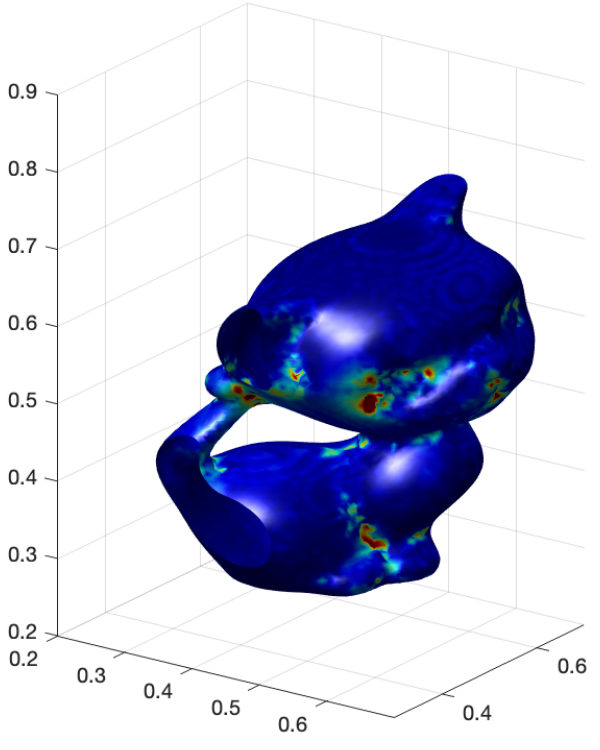}~
\includegraphics[height=2.2in]{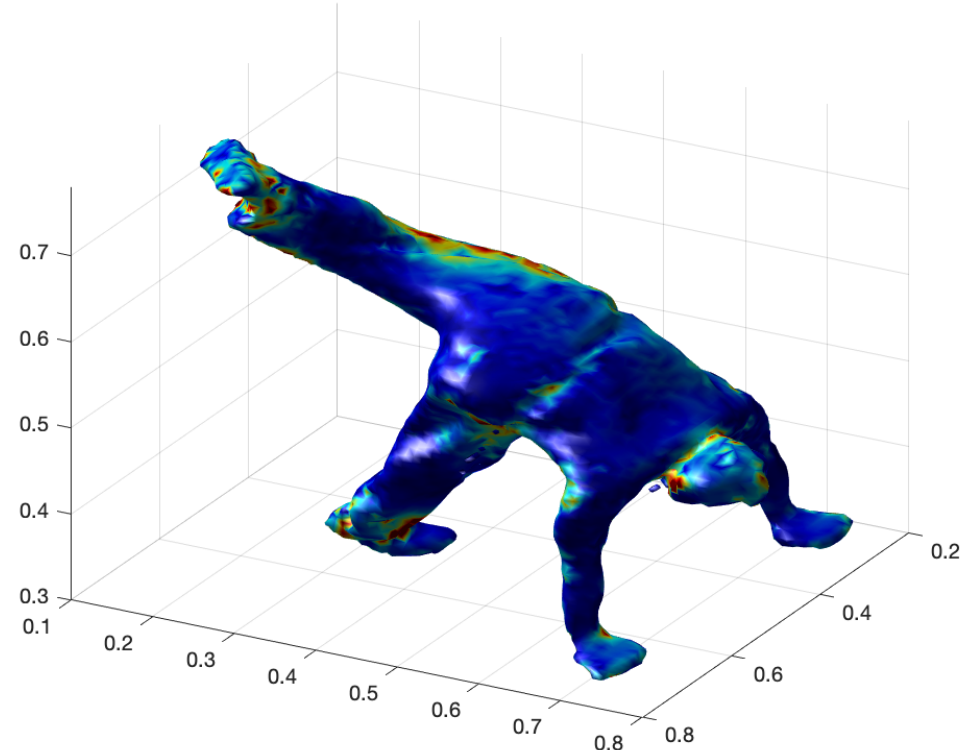}~
\includegraphics[height=2.2in]{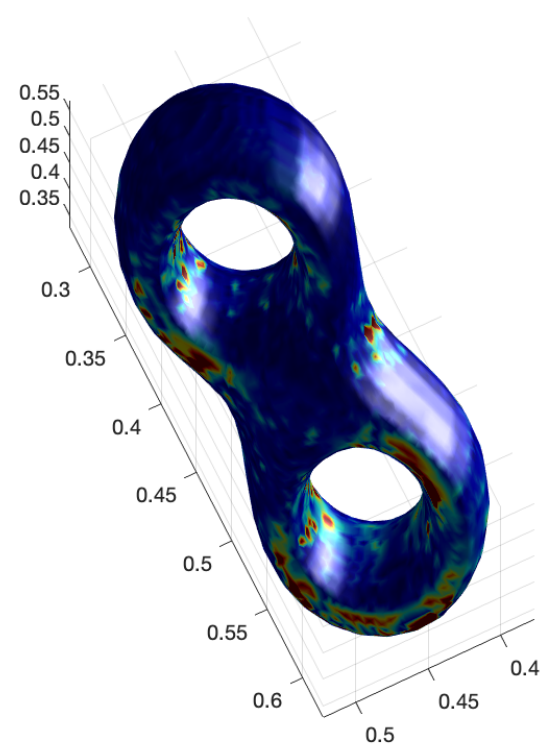}
\end{center}
\caption{From left: IFE solutions for a kitten shaped interface; a human shaped interface; and an eight shaped interface.}
\label{fig: more interfaces}
\end{figure}

The five test geometries ---Duck, Armadillo, Kitten, Human, and the numeral Eight--- are summarized in Table \ref{tab:numerical_summary}. Interface problems for the Duck and Armadillo are solved on a structured mesh with $N=100$. The corresponding surface error maps (right plots of Figures \ref{fig: duck} and \ref{fig: armadillo}) show larger errors on the Armadillo, reflecting its higher geometric complexity indicated by more cloud points in the Armadillo shape. The error maps for the other three interfaces are depicted in Figure \ref{fig: more interfaces}.

\begin{figure}
\begin{center}
\includegraphics[width=.49\textwidth]{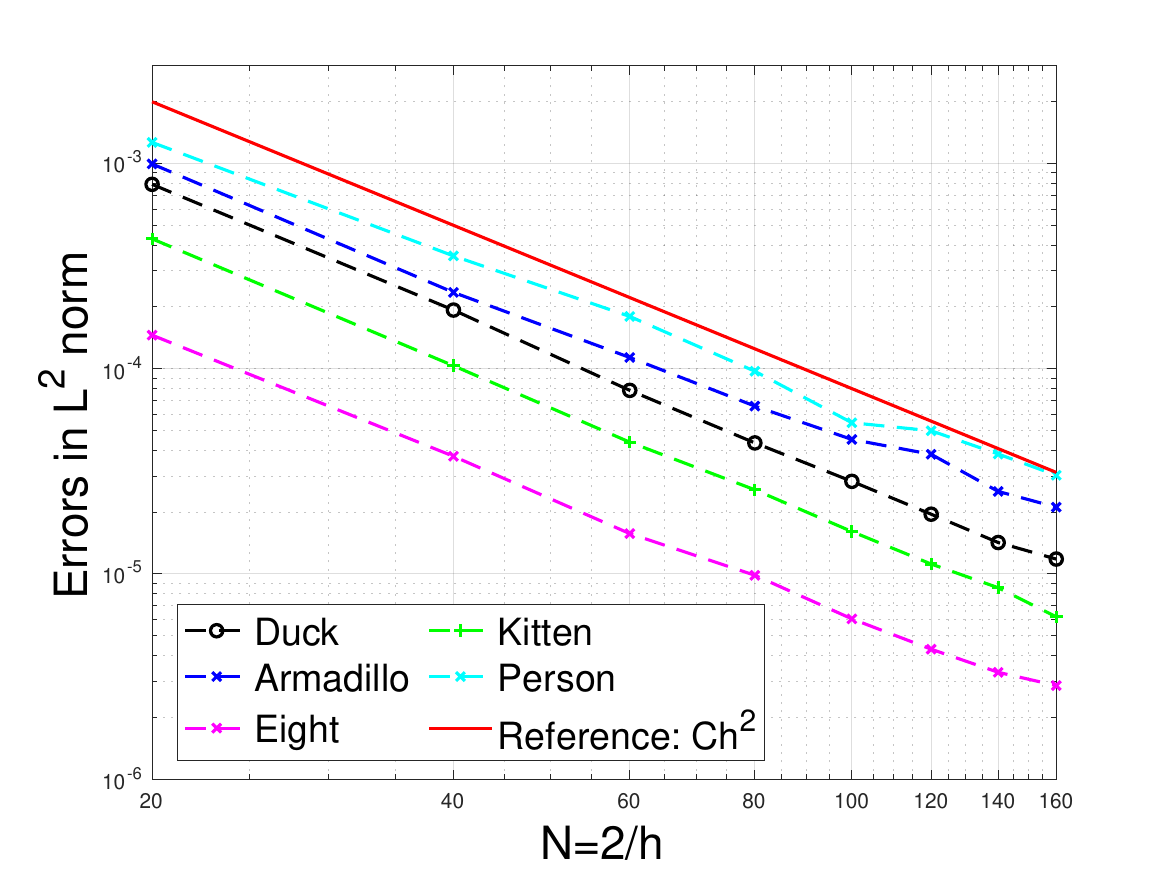}
\includegraphics[width=.49\textwidth]{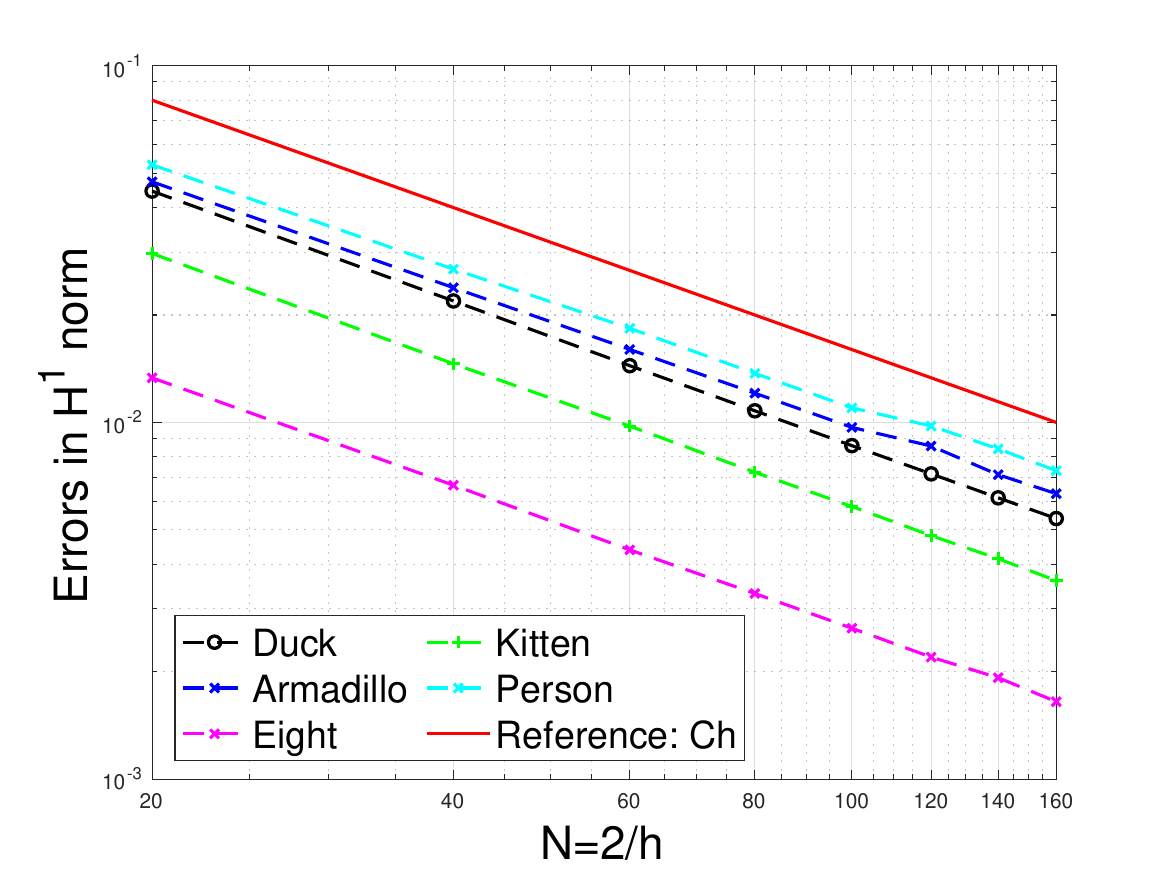}
\end{center}
\caption{Convergence in $L^2$ and $H^1$ norms of IFE solutions for the five types of interface geometry.}
\label{fig: convergence real world}
\end{figure}

For convergence study, each domain is uniformly partitioned to $N=20,40,60,\cdots, 160$. The $L^2$- and $H^1$-norm errors plotted in Figure \ref{fig: convergence real world} decreases at the optimal rates, demonstrating the accuracy and robustness of the enriched IFE method on challenging real‑world interfaces.

\begin{table}[h!]
  \centering
  \caption{Summary of data information for 3D Interfaces}
  \label{tab:numerical_summary}
  \begin{tabular}{|l|c|c|c|c|}
    \hline
    Interface Shape & \# Cloud Points  & Domain $\Omega$ & Figure Reference \\ \hline
    Duck              & 9640     & $[0.1,0.9]\times[0.1,0.9]\times [0.1,0.9] $  & Figure \ref{fig: duck} \\
    Armadillo       & 16519   &  $[0.1,0.9]\times[0.1,0.9]\times [0.1,0.9]$ & Figure \ref{fig: armadillo} \\
    Kitten             & 11039    & $[0.2,0.8]\times [0.1,0.9]\times[0.2,0.8]$   & Figure \ref{fig: more interfaces} \\
    Human           & 10002   &  $[0,1]\times[0,1]\times[0.2,0.8]$ &Figure  \ref{fig: more interfaces} \\
    Eight                & 766      & $[0.2,0.7]\times[0.3,0.6]\times[0.1,0.8]$ & Figure \ref{fig: more interfaces}  \\ \hline
  \end{tabular}
\end{table}

\section{Conclusion}
We developed an enriched immersed finite element (IFE) method for solving three-dimensional interface problems with non-homogeneous jump conditions. By applying a homogenization strategy, we construct enrichment function directly from the jump data, so the resulting homogeneous IFE space remains isomorphic, uniformly in the interface position, to the standard finite element space on the same mesh. We proved optimal error estimates and uniformly bounded condition numbers. Extensive numerical experiments confirm these theoretical results. 

\begin{appendix}

\section{Proof Lemma \ref{q_est} for $t=0$.}
\label{append1}
For each $\tilde{q}\in H^{s+1/2}(\omega_T)$, we introduce $\tilde{q}^{\bot}\in H^s(\Gamma^{\omega_T}_h)$ as the projection of $\tilde{q}$ onto $\Gamma^{\omega_T}_h$: for every $\bfx\in \Gamma^{\omega_T}_h$, define $\tilde{q}^{\bot}(\bfx) := q(\tilde{\bfx})$ where $\tilde{\bfx}\in\Gamma^{\omega_T}$ such that its projection onto $\Gamma^{\omega_T}_h$ is exactly $\bfx$. Use the plane $\Gamma^{\omega_T}_h$ to construct a local coordinate system denoted as $(\hat{x}_1,\hat{x}_2,\hat{x}_3)$, and assume that the surface $\Gamma^{\omega_T}$ can be described by a function $\hat{x}_3 = f(\hat{x}_1,\hat{x}_2)$. Then, we have $\tilde{q}^{\bot}(\hat{x}_1,\hat{x}_2)=\tilde{q}(\hat{x}_1,\hat{x}_2,f(\hat{x}_1,\hat{x}_2))$. Now, we consider the following decomposition
\begin{equation}
\label{q_est_eq1_1}
\| \tilde{q} - \Pi^k_{\Gamma^{\omega_T}_h} \tilde{q} \|_{L^2(\Gamma^{\omega_T}_h)} \le  \underbrace{\| \tilde{q} -  \tilde{q}^{\bot} \|_{L^2(\Gamma^{\omega_T}_h)}}_{(I)} + \underbrace{ \| \tilde{q}^{\bot} -   \Pi^k_{\Gamma^{\omega_T}_h} \tilde{q}^{\bot} \|_{L^2(\Gamma^{\omega_T}_h)} }_{(II)} + \underbrace{ \|  \Pi^k_{\Gamma^{\omega_T}_h}( \tilde{q} - \tilde{q}^{\bot}) \|_{L^2(\Gamma^{\omega_T}_h)} }_{(III)}.
\end{equation}
For $(I)$, we note that
\begin{equation}
\label{q_est_eq1_2}
\tilde{q}^{\bot}(\hat{x}_1,\hat{x}_2) = \tilde{q}(\hat{x}_1,\hat{x}_2,f(\hat{x}_1,\hat{x}_2)) = \tilde{q}(\hat{x}_1,\hat{x}_2,0) + \int_{0}^{f(\hat{x}_1,\hat{x}_2)} \partial_{\hat{x}_3}\tilde{q}(\hat{x}_1,\hat{x}_2,\hat{x}_3) \dd \hat{x}_3.
\end{equation}
It then implies, with the H\"older's inequality, that
\begin{equation}
\begin{split}
\label{q_est_eq1_3}
\| \tilde{q} -  \tilde{q}^{\bot} \|^2_{L^2(\Gamma^{\omega_T}_h)} &  = \int_{\Gamma^{\omega_T}_h} \left(  \int_{0}^{f(\hat{x}_1,\hat{x}_2)} \partial_{\hat{x}_3}\tilde{q}(\hat{x}_1,\hat{x}_2,\hat{x}_3) \dd \hat{x}_3 \right)^2 \dd\hat{x}_1\dd\hat{x}_2 \\
& \le   \int_{\Gamma^{\omega_T}_h} |f(\hat{x}_1,\hat{x}_2)|  \int_{0}^{f(\hat{x}_1,\hat{x}_2)} |\partial_{\hat{x}_3}\tilde{q}(\hat{x}_1,\hat{x}_2,\hat{x}_3)|^2 \dd \hat{x}_3  \dd\hat{x}_1\dd\hat{x}_2 \lesssim h_T^2   \| \tilde{q} \|^2_{H^1(\omega^{\text{int}}_T)}
\end{split}
\end{equation}
where we have used \eqref{lem_geo_gamma_eq01} in the last inequality. 
As for $(II)$, the well-known projection property on convex region by Proposition 6.1 \cite{1980DupontScott} immediately yields
\begin{equation}
\begin{split}
\label{q_est_eq1}
\| \tilde{q}^{\bot} - \Pi^k_{\Gamma^{\omega_T}_h} \tilde{q}^{\bot} \|^2_{L^2(\Gamma^{\omega_T}_h )} \lesssim h^{2\min\{k+1,s\}} | \tilde{q}^{\bot} |^2_{H^s(\Gamma^{\omega_T}_h)}.
\end{split}
\end{equation}
Given any two points $\bfx,\bfy\in \Gamma^{\omega_T}_h$, consider the corresponding points $\tilde{\bfx},\tilde{\bfy}\in \Gamma^{\omega_T}$, and let $\theta$ be the angle between $\bfx\bfy$ and $\tilde{\bfx}\tilde{\bfy}$. By elementary geometry, the angle between $\bfx\bfy$ and $\tilde{\bfx}\tilde{\bfy}$ denoted by $\theta$ is bounded by the maximum angle between the normal vectors to $\Gamma^T$ and $\Gamma^T_h$. Then, \eqref{lem_geo_gamma_eq02} implies that 
\begin{equation}
\label{q_est_eq2}
\| \tilde{\bfx} - \tilde{\bfy} \| = \| \bfx - \bfy \|/\cos(\theta) = \| \bfx - \bfy \|/(1-  c^2_{\Gamma,2}h^2_T/2 ) \lesssim \| \bfx - \bfy \|.
\end{equation}
For each double-index $\bfalpha=(\alpha_1,\alpha_2)$ with $\alpha_i\ge 0$, we can use chain rule to obtain 
\begin{equation}
\label{q_est_eq3}
  \partial_{\bfalpha}\tilde{q}^{\bot}(\bfx) =  \partial_{\bfalpha}q(\hat{x}_1,\hat{x}_2,f(\hat{x}_1,\hat{x}_2))  = \sum_{|\bfalpha'|=|\alpha|} C^f_{\alpha'}(\tilde{\bfx}) \partial_{\alpha'} \tilde{q}(\tilde{\bfx})
\end{equation}
where $\bfalpha'$ is a triple-index, and $C^f_{\alpha'}$ is a function only depending on $\bfalpha'$ and $f$.
Let $\bar{s}$ be the largest integer not greater than $s$, and define $s'=s-\bar{s}$. By definition, applying \eqref{q_est_eq2} and \eqref{q_est_eq3} we arrive at
\begin{equation}
\begin{split}
\label{q_est_eq4}
 | \tilde{q}^{\bot} |^2_{H^s(\Gamma^{\omega_T}_h)} & = \sum_{|\bfalpha |=\bar{s}} \int_{\Gamma^{\omega_T}_h} \int_{\Gamma^{\omega_T}_h} \frac{| \partial_{\bfalpha} \tilde{q}^{\bot}(\bfx) - \partial_{\bfalpha} \tilde{q}^{\bot}(\bfy)|^2}{\| \bfx- \bfy \|^{3+2s'}} \dd \bfx \dd \bfy \\
 & \lesssim \sum_{|\bfalpha' |=\bar{s}} \int_{\Gamma^{\omega_T}_h} \int_{\Gamma^{\omega_T}_h} \frac{|C^f_{\alpha'}(\tilde{\bfx}) \partial_{\bfalpha'}\tilde{q}(\tilde{\bfx}) - C^f_{\alpha'}(\tilde{\bfx}) \partial_{\bfalpha'}\tilde{q}(\tilde{\bfy})|^2}{\| \tilde{\bfx}- \tilde{\bfy} \|^{3+2s'}}  \sqrt{1+ \| \nabla f(\bfx)\|^2}  \sqrt{1+ \| \nabla f(\bfy)\|^2 }  \dd \bfx \dd \bfy \\
& = \sum_{|\bfalpha' |=\bar{s}} | C^f_{\alpha'} \partial_{\bfalpha'}\tilde{q} |_{H^{s'}(\Gamma)} \lesssim  \sum_{|\bfalpha' |=\bar{s}} |\partial_{\bfalpha'} \tilde{q} |_{H^{s'}(\Gamma^T)} \lesssim |\tilde{q}|_{H^s(\Gamma^T)}
 \end{split}
\end{equation}
where in the second last inequality we have used the multiplication theorem for fractional order Sobolev spaces by the assumption of the interface being sufficiently smooth and the hidden constant only depends on $f$. Putting \eqref{q_est_eq4} into \eqref{q_est_eq1} yields the estimate for $(II)$. Finally, $(III)$ follows from the stability of $\Pi^k_{\Gamma^{\omega_T}_h}$ and the estimate for $(I)$. These estimates together yield the desired result.

\end{appendix}


\end{document}